\documentclass[aap,MSNbibl,seceqn,dvips]{arximspdf}
\usepackage{algorithm}
\usepackage{algpseudocode}

\doi{10.1214/13-AAP953} 
\volume{25}
\issue{1}
\pubyear{2015}
\firstpage{1}
\lastpage{45}

\makeatletter

\newcommand{\rright}{\right}
\newcommand{\lleft}{\left}
\newcommand{\rrvert}{\vert}
\newcommand{\llvert}{\vert}

\newtheorem{cor}{Corollary}
\newtheorem{lem}{Lemma}
\newtheorem{prop}{Proposition}
\newproclaim{rem}{Remark}

\newcommand{\sgn}{\operatorname{sgn}}
\newcommand{\eqref}[1]{(\ref{#1})}

\renewcommand{\emptyset}{\varnothing}
\def\sfrac#1#2{#1/#2}

\def\vafrac#1#2{(#1)/(#2)}

\makeatother

\begin{document}
\begin{frontmatter}

\title{On the stability of some controlled Markov chains and its applications
to stochastic approximation with Markovian dynamic}
\runtitle{On the stability of controlled Markov chains}

\begin{aug}
\author[A]{\fnms{Christophe}~\snm{Andrieu}\corref{}\thanksref{T1}\ead[label=e1]{c.andrieu@bristol.ac.uk}},
\author[A]{\fnms{Vladislav B.}~\snm{Tadi\'{c}}\ead[label=e2]{v.b.tadic@bristol.ac.uk}}
\and
\author[B]{\fnms{Matti}~\snm{Vihola}\thanksref{T2}\ead[label=e3]{matti.vihola@iki.fi}}
\runauthor{C. Andrieu, V. B. Tadi\'{c} and M. Vihola}
\affiliation{University of Bristol, University of Bristol and University
of Jyv\"{a}skyl\"{a}}
\address[A]{C. Andrieu\\
V. B. Tadi\'{c}\\
School of Mathematics\\
University of Bristol\\
BS8 1TW\\
United Kingdom\\
\printead{e1}\\
\printead{e2}}

\address[B]{M. Vihola\\
Department of Mathematics and Statistics\\
University
of Jyv\"{a}skyl\"{a}\\
P.O. Box 35\\
FI-40014\\
Finland\\
\printead{e3}}
\end{aug}
\thankstext{T1}{Supported in part by an EPSRC Advance Research
Fellowship and a Winton Capital research prize for this research.}
\thankstext{T2}{Supported in part by the Academy of Finland (project
250575) and by the Finnish Academy of Science and Letters, Vilho,
Yrj\"{o} and Kalle V\"{a}is\"{a}l\"{a} Foundation.}

\received{\smonth{5} \syear{2012}}
\revised{\smonth{5} \syear{2013}}

%
\begin{abstract}
We develop a practical approach to establish the stability, that is,
the recurrence in a given set, of a large class of controlled Markov
chains. These processes arise in various areas of applied science
and encompass important numerical methods. We show in particular how
individual Lyapunov functions and associated drift conditions for
the parametrized family of Markov transition probabilities and the
parameter update can be combined to form Lyapunov functions for the
joint process, leading to the proof of the desired stability property.
Of particular interest is the fact that the approach applies even
in situations where the two components of the process present a time-scale
separation, which is a crucial feature of practical situations. We
then move on to show how such a recurrence property can be used in
the context of stochastic approximation in order to prove the convergence
of the parameter sequence, including in the situation where the so-called
stepsize is adaptively tuned. We finally show that the results apply
to various algorithms of interest in computational statistics and
cognate areas.
\end{abstract}

%

\end{frontmatter}

\section{Introduction: Recurrence of controlled MC and compound
drifts}
\label{121212121}

The class of controlled Markov chain processes underpins numerous
models or algorithms encountered in various areas of engineering or
science (e.g., control, EM algorithm, adaptive MCMC). Consider the space
$(\mathsf{X},\mathcal{B}(\mathsf{X}))$ where $\mathsf{X}\subset
\mathbb{R}^{n_{x}}$
for some $n_{x}\geq1$, a parametrized family of Markov transition
probabilities $\{P_{\theta},\theta\in\Theta\}$ (for some set
$\Theta\subset\mathbb{R}^{n_{\theta}}$)
such that for any $\theta,x\in\Theta\times\mathsf{X}$, $P_{\theta
}(x,\cdot)$
is a probability distribution on $(\mathsf{X},\mathcal{B}(\mathsf{X}))$.
The class of controlled Markov chains we consider in this paper consists
of the class of processes defined on $ ((\Theta\times\mathsf
{X})^{\mathbb{N}},(\mathcal{B}(\Theta)\otimes\mathcal{B}(\mathsf
{X)})^{\otimes\mathbb{N}} )$
initialized at some $(\theta_{0},X_{0})=(\theta,x)\in\Theta\times
\mathsf{X}$,
with probability distribution denoted $\mathbb{P}_{\theta,x}(\cdot)$
[and associated expectation $\mathbb{E}_{\theta,x}(\cdot)$] and defined
recursively for $i\geq0$ as follows:
%
\begin{eqnarray}\label{eqdefinitioncontrolledmc}
&X_{i+1}|(\theta_{0},X_{0},X_{1},
\ldots,X_{i})\sim P_{\theta
_{i}}(X_{i},\cdot),&
\nonumber
\\[-8pt]\\[-8pt]
&\theta_{i+1}:=\phi_{i+1}(\theta_{0},X_{0},X_{1},
\ldots,X_{i+1}),&\nonumber
\end{eqnarray}
for a family of mappings $\{\phi_{i}\dvtx \Theta\times\mathsf
{X}^{i+1}\rightarrow\Theta\}$.
The present paper is concerned with the stability of the sequence
$\{\theta_{i},X_{i}\}$, or more precisely, the recurrence of such
a process in a set $\mathcal{C}\subset\mathsf{X}\times\Theta$;
that is,
we aim to develop practically relevant tools to establish that $\{
\theta_{i},X_{i}\}$
visits $\mathcal{C}$ infinitely often $\mathbb{P}_{\theta,x}$-a.s.
Such a form of stability is central to establish important properties
of the process which, depending on the context, range from the existence
of an invariant distribution for the process or its marginals to the
convergence of the parameter sequence $\{\theta_{i}\}$ to a set of
values of particular interest. This is largely an open problem despite
its practical relevance as illustrated and discussed later in the
paper. The following toy example illustrates the potential difficulties
one may face. Let $\mathsf{X}=\{0,1\}$ and consider the transition
matrix\looseness=-1
\[
P_{\theta}=\lleft[ %
\begin{array} {c@{\quad}c} 1-\exp\bigl(-|\theta|\bigr) &
\exp\bigl(-|\theta|\bigr)
\\
\exp\bigl(-|\theta|\bigr) & 1-\exp\bigl(-|\theta|\bigr) \end{array} %
\rright],
\]
with $\Theta=\mathbb{R}$. This transition matrix has $\pi=(1/2,1/2)$
as invariant distribution, and its second eigenvalue is $\lambda
=1-2\exp(-|\theta|)$.
Set $\theta_{i+1}=\theta_{i}+a/i[1/2-X_{i+1}]$ for some $a>0$. One
could expect $\{\theta_{i}\}$ to converge to a finite value, but
following the argument in \cite{younes}, Section~6.3, one can in fact
show that for some values of $a$, with positive probability, $\{X_{i}\}$
may get stuck in either states while $\{\theta_{i}\}$ diverges. Ergodicity
is lost here due to the fact that $\mathcal{C}_{0}=\Theta\times\{0\}$
or $\mathcal{C}_{1}=\Theta\times\{1\}$ is not visited infinitely
often with probability one.\looseness=-1

The remainder of the paper is organized as follows.
In Section~\ref{secCompound-Lyapunov-functions}
we introduce our methodology, which relies on a classical Lyapunov
function/drift argument to establish recurrence of the joint process
$\{\theta_{i},X_{i}\}$ to a set $\mathcal{C}$ (Lemma~\ref
{lemgeneral-drift-stability}).
Our main result in this section is Theorem~\ref{thmcombined-recurrence}
where it is shown how individual drift conditions of type (\ref{eqw-ineq})
and (\ref{eqV-ineq}) characterizing the evolution from $X_{i}$
to $X_{i+1}$ and $\theta_{i}$ to $\theta_{i+1}$ in (\ref
{eqdefinitioncontrolledmc}),
respectively, can be combined into a joint drift condition in order
to characterize the joint dynamic and establish recurrence to a set
$\mathcal{C}$. It is worth pointing out that the result applies even
in situations where this dynamic exhibits a time-scale separation,
which, as we shall see in the application section, is of practical interest.
This result captures the main ideas behind our general strategy to
establish recurrence, but for simplicity and clarity, remain unspecific
about how the abstract conditions may be relevant in practice.

Section~\ref{1313131313} contains the main practical results
of the paper, Theorem~\ref{propestablishcompositedrift} and its
corollary, where we show how familiar (e.g., \cite{meyn:tweedie:1993}),
but $\theta$-dependent, drift conditions characterizing the evolution
of homogeneous Markov chains with transitions $\{P_{\theta},\theta\in
\Theta\}$
[see assumption \ref{hypfixedthetaconditiononPtheta}] can be
combined with a class of drift conditions characterizing the evolution
of the parameter $\theta$ [see assumption \ref{hypwdriftSA1}]
in order to apply our earlier abstract results for the stability of
the joint process.

In Section~\ref{secmotivation} we focus on a practically important
class of updates for the parameter $\theta$, known as stochastic
approximation \cite{benveniste}, which covers all our subsequent
applications. The corresponding processes aim to find the zeroes of
a function of the parameter $\theta$ and can be seen as noisy gradient
algorithms. The aims of the section are to introduce sufficient background
for the application section and to establish Theorem~\ref
{thlocalergodicityimpliesstability}.
The result of this theorem highlights the central role played by recurrence
in an appropriate set $\mathcal{C}$ in this scenario, in order to
ensure that such numerical methods are stable and that they achieve
their goal.

Finally in Section~\ref{171717171717} we show how the results
established earlier apply in the context of adaptive Markov chain
Monte Carlo (MCMC) algorithms, a particular type of MCMC algorithms
which aim to optimize their performance ``on the fly.'' More specifically,
we show that our general results apply to both the AM algorithm of
\cite{haario:saksman:tamminen:2001} but also the coerced acceptance
probability algorithm \cite{andrieu:robert:2001,atchade:rosenthal:2006}
and a novel variation.

\section{Compound Lyapunov functions for some two timescale controlled\break Markov~chains}
\label{secCompound-Lyapunov-functions}

The approach we adopt throughout this paper relies on a classical
Lyapunov function and drift argument commonly used in the (homogeneous)
Markov chain setting \cite{meyn:tweedie:1993}. Due to the potential
time inhomogeneity of the process above, it is useful to consider a
sequence of Lyapunov functions $\{W_{i}\}$ satisfying a sequence
of drift conditions and leading to the following classical result,
provided here together with its proof (in Appendix~\ref
{secAppendix-for-Section1})
for completeness only. Hereafter, for any $i\geq0$ we let $\mathcal
{F}_{i}:=\sigma (\theta_{0},X_{0},X_{1},\ldots,X_{i} )$
and for any $u,v\in\mathbb{R}^{2}$ we define $u\vee v:=\max\{u,v\}$
and $u\wedge v:=\min\{u,v\}$.

%
\begin{lem}
\label{lemgeneral-drift-stability} Let $\{W_{i}\}$ be a sequence
of functions $W_{i}\dvtx \Theta\times\mathsf{X}\to[0,\infty)$ such that
for the controlled Markov chain defined in (\ref{eqdefinitioncontrolledmc})
for all $\theta,x\in\Theta\times\mathsf{X}$:
\begin{longlist}[(2)]
\item[(1)] for all $i\geq0$, $\mathbb{E}_{\theta,x}[W_{i}(\theta
_{i},X_{i})]<\infty$,
\item[(2)] there exist $\mathcal{C}\subset\Theta\times\mathsf{X}$, a sequence
$\{\delta_{i},i\geq1\}$ of nonnegative scalars such that
$\sum_{i=1}^{\infty}\delta_{i}=\infty$ and an integer $i_{w}<\infty$
such that for all $i\ge i_{w}$, and whenever $(\theta
_{i},X_{i})\notin\mathcal{C}$,
$\mathbb{P}_{\theta,x}$-a.s.
%
\begin{equation}\label{eqgeneric-W-drift}
\mathbb{E}_{\theta,x}\bigl[W_{i+1}(\theta_{i+1},X_{i+1})
\mid\mathcal {F}_{i}\bigr]\le W_{i}(\theta_{i},X_{i})-
\delta_{i+1}.
\end{equation}
\end{longlist}

Then $\sum_{i=1}^{n}\mathbb{I}\{(\theta_{i},X_{i})\in\mathcal{C}\}
=\infty$,
$\mathbb{P}_{\theta,x}$-a.s.\vadjust{\goodbreak}
\end{lem}

The main result of this section consists of showing that it is possible
to construct joint Lyapunov function sequences $\{W_{i}\}$ which
satisfy drifts to a set $\mathcal{C}$, such that the conditions of
Lemma~\ref{lemgeneral-drift-stability} hold, from two separate Lyapunov
functions $w(\theta)$ and $V(x)$ each satisfying an individual drift
condition characterizing the two respective updates involved in the
definition of $\{\theta_{i},X_{i}\}$ in \eqref{eqdefinitioncontrolledmc}.
The form of these individual drifts is given below in \eqref{eqw-ineq}
and \eqref{eqV-ineq}: it is worth pointing out that we allow the
drift on $w(\theta)$ to vanish with time since $\{\gamma_{i}\}$
may be allowed to vanish. This is practically very relevant since
in many situations of interest the ``size'' of the increments $|\theta
_{i+1}-\theta_{i}|$
may vanish as $i\rightarrow\infty$ while that of $|X_{i+1}-X_{i}|$
may not. The role of the sequence $\{\gamma_{i}\}$ is to accommodate
the possibility of two distinct timescales for the two updates in
\eqref{eqdefinitioncontrolledmc}---examples are numerous and some
will be presented later in Sections~\ref{secmotivation} and \ref
{171717171717}.
We will consider two scenarios which share very similar assumptions,
and will be labeled with $s\in\{0,1\}$.
\renewcommand\thelonglist{(A1)}
\renewcommand\labellonglist{\thelonglist}
\begin{longlist}[(A1)]
\item \label{aw-v-bounds} Suppose
$V\dvtx \mathsf{X}\to[1,\infty)$
and $w\dvtx \Theta\to[1,\infty)$ are two functions such that there exist
functions $\Delta_{w},\Delta_{V}\dvtx \Theta\times\mathsf{X}\to\mathbb{R}$,
a set $\mathcal{C}\subset\Theta\times\mathsf{X}$, a sequence of strictly
positive integers $\{\gamma_{i},i\geq1\}$ such that:
\begin{longlist}[(2)]
\item[(1)]$\{\gamma_{i}\}$ is bounded,
\item[(2)] for some integer $i_{0}\geq0$, $\mathbb{P}_{\theta,x}$-a.s. the
following individual drifts hold for all $i\ge i_{0}$:
%
\begin{eqnarray}
\label
{eqw-ineq}\mathbb{E}_{\theta,x}\bigl[w(\theta_{i+1})\mid\mathcal{F}_{i}
\bigr]&\le& w(\theta_{i})-\gamma_{i+1}\Delta_{w}(
\theta_{i},X_{i}),
\\
\label{eqV-ineq}\mathbb{E}_{\theta,x}\bigl[V(X_{i+1})\mid\mathcal{F}_{i}
\bigr]&\le& V(X_{i})-\Delta_{V}(\theta_{i},X_{i})
\end{eqnarray}
and $\mathbb{E}_{\theta,x}[w(\theta_{i})]<\infty$ and $\mathbb
{E}_{\theta,x}[V(X_{i})]<\infty$,
\item[(3)] there exist constants $\delta\in(0,\infty)$ and $\upsilon
_{v},\upsilon_{w}\in(0,1]$
such that
%
\begin{equation}\label{eqpsi-notin-c}
\upsilon_{w}\frac{\Delta_{w}(\theta,x)}{w^{1-\upsilon_{w}}(\theta
)}+\upsilon_{v}
\frac{\Delta_{V}(\theta,x)}{V^{1-\upsilon
_{v}}(x)}\geq\delta w^{s\times\upsilon_{w}}(\theta) \qquad \mbox{for } (\theta,x)\notin
\mathcal{C}
\end{equation}
and
\[
\sup_{(\theta,x)\in\mathcal{C}}\bigl|\Delta_{w}(\theta,x)\bigr|\vee\bigl|\Delta
_{V}(\theta,x)\bigr|<\infty.
\]
\end{longlist}
\end{longlist}

The following theorem establishes two recurrence results for $\{\theta
_{i},X_{i}\}$
to $\mathcal{C}$. The first result requires the strongest set of
assumptions but also establishes a stronger result, namely that the
first moment of the return times to $\mathcal{C}$ are uniformly bounded
in time. The second result requires weaker assumptions but does not
guarantee the existence of a uniform in time upper bound on characteristics
of the return times. A~particular contribution here is the rescaling
of either the Lyapunov function $w(\theta)$ or $V(x)$ in order to
allow for their respective drift terms to be compared on the same
time scale.

%
\begin{thm}
\label{thmcombined-recurrence}Consider the controlled Markov chain
defined in (\ref{eqdefinitioncontrolledmc}). Define the sequences
of functions $\{W_{i}\dvtx \Theta\times\mathsf{X}\rightarrow[1,\infty)\}$
and $\{U_{i}\dvtx \Theta\times\mathsf{X}\rightarrow[1,\infty)\}$ for
$i\geq1$
and $\theta,x\in\Theta\times\mathsf{X}$ as follows:
\[
W_{i}(\theta,x):=V^{\upsilon_{v}}(x)+w^{\upsilon_{w}}(\theta )/
\gamma_{i} \quad \mbox{and}\quad  U_{i}(\theta,x):=
\gamma_{i}W_{i}(\theta,x),
\]
where $\{\gamma_{i}\}$, $w(\cdot)$, $V(\cdot)$, $\upsilon_{v}$
and $\upsilon_{w}$ are as in \textup{\ref{aw-v-bounds}}, which is assumed
to hold. Then:
\begin{longlist}[(1)]
\item[(1)] if $s=1$ and $\bar{\ell}:=\limsup_{i\to\infty}(\gamma
_{i+1}^{-1}-\gamma_{i}^{-1})<\delta$,
then for any $\delta_{W}\in(0,\delta-\bar{\ell})$ there exists
$i_{W}\geq i_{0}$
such that for any $i\geq i_{W}$, whenever $(\theta_{i},X_{i})\notin
\mathcal{C}$,
$\mathbb{P}_{\theta,x}$-a.s.
%
\begin{equation}\label{eqWdrift}
\mathbb{E}_{\theta,x}\bigl[W_{i+1}(\theta_{i+1},X_{i+1})
\mid\mathcal {F}_{i}\bigr]\le W_{i}(\theta_{i},X_{i})-
\delta_{W},
\end{equation}
and $\mathbb{E}_{\theta,x}[W_{i}(\theta_{i},X_{i})]<\infty$, and
$\sum_{i=1}^{\infty}\mathbb{I}\{(\theta_{i},X_{i})\in\mathcal{C}\}
=\infty$,
$\mathbb{P}_{\theta,x}$-almost surely,
\item[(2)] if $s=0$, $\{\gamma_{i}\}$ is nonincreasing then for any $i\geq i_{0}$,
whenever $(\theta_{i},X_{i})\notin\mathcal{C}$, $\mathbb{P}_{\theta,x}$-a.s.
%
\begin{equation}\label{eqUdrift}
\mathbb{E}_{\theta,x}\bigl[U_{i+1}(\theta_{i+1},X_{i+1})
\mid\mathcal {F}_{i}\bigr]\le U_{i}(\theta_{i},X_{i})-
\delta\gamma_{i+1}
\end{equation}
and moreover $\mathbb{E}_{\theta,x}[U_{i}(\theta_{i},X_{i})]<\infty$.
If in addition $\sum_{i=1}^{\infty}\gamma_{i}=\infty$, then\break $\sum_{i=1}^{\infty}\mathbb{I}\{(\theta_{i},\allowbreak X_{i})\in\mathcal{C}\}
=\infty$,
$\mathbb{P}_{\theta,x}$-almost surely.
\end{longlist}
\end{thm}

\begin{pf}
By \ref{aw-v-bounds}, Jensen's inequality and the classical concavity
identity $(1+x)^{\upsilon}\leq1+\upsilon x$ for $x\in[-1,\infty)$
and $\upsilon\in(0,1]$, we have for any $i\ge i_{0}$ and $\mathbb
{P}_{\theta,x}$-a.s.
%
\begin{eqnarray}\label{eqfirstinequaforWi}
&&\mathbb{E}_{\theta,x}\bigl[  W_{i+1}(\theta_{i+1},X_{i+1})
\mid\mathcal {F}_{i}\bigr]\nonumber
\\
&&\qquad \le V^{\upsilon_{v}}(X_{i}) \biggl(1-\frac{\Delta_{V}(\theta
_{i},X_{i})}{V(X_{i})}
\biggr)^{\upsilon_{v}}+\gamma _{i+1}^{-1}w^{\upsilon_{w}}(
\theta_{i}) \biggl(1-\gamma_{i+1}\frac
{\Delta_{w}(\theta_{i},X_{i})}{w(\theta_{i})}
\biggr)^{\upsilon
_{w}}
\nonumber
\\[-8pt]\\[-8pt]
&&\qquad \leq V^{\upsilon_{v}}(X_{i}) \biggl(1-\upsilon_{v}
\frac{\Delta
_{V}(\theta_{i},X_{i})}{V(X_{i})} \biggr)+\gamma _{i+1}^{-1}w^{\upsilon_{w}}(
\theta_{i}) \biggl(1-\gamma_{i+1}\upsilon _{w}
\frac{\Delta_{w}(\theta_{i},X_{i})}{w(\theta_{i})} \biggr)
\nonumber
\\
&&\qquad =  W_{i}(\theta_{i},X_{i})+ \bigl(
\gamma_{i+1}^{-1}-\gamma _{i}^{-1}
\bigr)w^{\upsilon_{w}}(\theta_{i})- \biggl(\upsilon _{v}
\frac{\Delta_{V}(\theta_{i},X_{i})}{V^{1-\upsilon
_{v}}(X_{i})}+\upsilon_{w}\frac{\Delta_{w}(\theta
_{i},X_{i})}{w^{1-\upsilon_{w}}(\theta_{i})} \biggr).\hspace*{-20pt}
\nonumber
\end{eqnarray}
Now consider the scenario where $s=1$. Let $\delta_{W}\in(0,\delta
-\bar{\ell})$
and $i_{W}\ge i_{0}$ be such that $\sup_{i\ge i_{w}}(\gamma
_{i+1}^{-1}-\gamma_{i}^{-1})<\delta-\delta_{W}$.
Then, for all $i\ge i_{W}$ and $(\theta_{i},X_{i})\notin\mathcal{C}$,
$\mathbb{P}_{\theta,x}$-a.s.
\[
\mathbb{E}_{\theta,x}\bigl[W_{i+1}(\theta_{i+1},X_{i+1})
\mid\mathcal {F}_{i}\bigr]\le W_{i}(\theta_{i},X_{i})-
\delta_{W}.
\]
Let $C:=[\sup_{i\geq i_{0}}\gamma_{i}(\gamma_{i+1}^{-1}-\gamma
_{i}^{-1})]\vee[\sup_{(\theta,x)\in\mathcal{C}}|\bar{\gamma
}_{i_{0}}\Delta_{w}(\theta,x)|\vee|\Delta_{V}(\theta,x)|]$
with $\bar{\gamma}_{i_{0}}=\sup_{i\geq i_{0}}\gamma_{i}$. Now for
any $i\ge i_{0}$ and $(\theta_{i},X_{i})\in\Theta\times\mathsf{X}$
we have, starting with the first inequality in (\ref{eqfirstinequaforWi}),
%
\begin{eqnarray*}
&&\mathbb{E}_{\theta,x}\bigl[W_{i+1}(\theta_{i+1},X_{i+1})\mid\mathcal{F}_{i}\bigr]
\\
&&\qquad  \le(1+C)^{\upsilon_{v}}V^{\upsilon_{v}}(X_{i})+(1+C)^{\upsilon
_{w}}
\bigl[\gamma_{i}\bigl(\gamma_{i+1}^{-1}-
\gamma_{i}^{-1}\bigr)+1\bigr]w^{\upsilon
_{w}}(
\theta_{i})/\gamma_{i}
\\
&&\qquad  \leq(1+C)^{2}W_{i}(\theta_{i},X_{i})
\end{eqnarray*}
for $s\in\{0,1\}$. From these inequalities we therefore deduce that
for any $i\geq i_{0}$, $\mathbb{E}_{\theta,x}[W_{i}(\theta
_{i},X_{i})]\le(1+C)^{2(i-i_{0})}\mathbb{E}_{\theta,x}[W_{i_{0}}(\theta_{i_{0}},X_{i_{0}})]<\infty$
where the last inequality follows from our assumptions. For the scenario
where $s=0$ with $U_{i}(\theta,x)=\gamma_{i}W_{i}(\theta,x)$, we
obtain from (\ref{eqfirstinequaforWi})
\begin{eqnarray*}
\mathbb{E}_{\theta,x}\bigl[U_{i+1}(\theta_{i+1},X_{i+1})
\mid\mathcal {F}_{i}\bigr]&\le& U_{i}(\theta_{i},X_{i})+(
\gamma_{i+1}-\gamma _{i})W_{i}(
\theta_{i},X_{i})
\\
&&{}+\gamma_{i+1} \bigl(\gamma_{i+1}^{-1}-
\gamma_{i}^{-1} \bigr)w^{\upsilon_{w}}(\theta_{i})\\
&&{}-
\gamma_{i+1} \biggl(\upsilon_{v}\frac
{\Delta_{V}(\theta_{i},X_{i})}{V^{1-\upsilon_{v}}(X_{i})}+\upsilon
_{w}\frac{\Delta_{w}(\theta_{i},X_{i})}{w^{1-\upsilon_{w}}(\theta
_{i})} \biggr)
\end{eqnarray*}
and since
\[
(\gamma_{i+1}-\gamma_{i})W_{i}(\theta,x)+
\gamma_{i+1} \bigl(\gamma _{i+1}^{-1}-
\gamma_{i}^{-1} \bigr)w^{\upsilon_{w}}(\theta)=(\gamma
_{i+1}-\gamma_{i})V(x)
\]
and $\{\gamma_{i}\}$ is nonincreasing, we conclude (\ref{eqUdrift})
for $(\theta_{i},X_{i})\in\mathcal{C}^{c}$. Notice further that
$U_{i}(\theta,x)\leq\gamma_{1}W_{i}(\theta,x)$,
implying $\mathbb{E}_{\theta,x}[U_{i}(\theta_{i},X_{i})]<\infty$
for any $i\geq i_{0}$. We now conclude in both scenarios with Lemma~\ref{lemgeneral-drift-stability}.
\end{pf}

Some comments are in order concerning the choice of the Lyapunov
functions and the assumptions. First we clarify the role of $\upsilon_{v}$
and $\upsilon_{w}$, which are additional degrees of freedom one may
find helpful to establish \eqref{eqpsi-notin-c} in regions of $\Theta
\times\mathsf{X}$
where $\Delta_{V}(\theta,x)$ [resp., $\Delta_{w}(\theta,x)$] is negative
and of large magnitude, but $V$ (resp.,~$w$) is itself large. Notice
also that more general concave transformations of $V$ and $w$ could
be considered for the definition of $W_{i}$ and $U_{i}$, but we
do not pursue this here. We would also like to point out that other
Lyapunov functions of the form $U_{i}^{\alpha}(\theta,x):=\gamma
_{i}^{\alpha}W_{i}(\theta,x)$
for $\alpha\geq0$ may be considered but we have found the scenarios
$\alpha=0$ and $\alpha=1$ to be of interest only. Finally whereas
it is clear that \eqref{eqpsi-notin-c} is stronger for $s=1$ than
$s=0$, we also note that $\limsup_{i\to\infty}(\gamma
_{i+1}^{-1}-\gamma_{i}^{-1})<\delta-\delta_{W}$
implies $\sum_{i=1}^{\infty}\gamma_{i}=\infty$.

In the next section we consider a practically relevant scenario encountered
in practice, for which we identify $\mathcal{C}$ and $\{W_{i}\}$,
but also establish an even stronger drift than in \eqref{eqWdrift}.
We will show in Section~\ref{171717171717} that such results are
satisfied in realistic scenarios.

\section{Simultaneous \texorpdfstring{$\theta$}{$theta$}-dependent drift conditions and
stability}
\label{1313131313}

The results presented in the previous section are rather abstract
since $\Delta_{w},\Delta_{V}$ and $\mathcal{C}$ are not specified.
Here we add some structure, and in particular, show how a simultaneous
drift condition on the family of Markov transition probabilities
$P_{\theta}$
for $\theta\in\Theta$, where the dependence on $\theta$ is explicit,
can be used to prove the stability of the sequence $\{\theta
_{i},X_{i}\}$
to a well-identified set $\mathcal{C}\subset\Theta\times\mathsf{X}$.
For ease of exposition we focus throughout this section on the situation
where $\phi_{i}:=\phi_{\gamma_{i}}$ for some family of updates $\{
\phi_{\gamma}\dvtx \Theta\times\mathsf{X}\rightarrow\Theta,\gamma\in
(0,\gamma^{+}]\}$
and a positive sequence $\{\gamma_{i}\}\subset(0,\gamma
^{+}]^{\mathbb{N}}$,
allowing us to define the update $\theta_{i+1}=\phi_{\gamma
_{i+1}}(\theta_{i},X_{i+1})$
for $i\geq0$. This directly covers most relevant applications in
computational statistics and can be easily generalized. As we shall
see, the realistic assumptions we use lead, in fact, to stronger results
than those of the previous section. For any $f\dvtx \mathsf{X}\rightarrow
\mathbb{R}$
we use the standard notation $P_{\theta}f(x):=\int_{\mathsf
{X}}P_{\theta}(x,\mathrm{d}y)f(y)$.
The $\theta$-dependent simultaneous drift conditions we consider
here are as follows:
\renewcommand\thelonglist{(A2)}
\begin{longlist}[(A2)]
\item \label{hypfixedthetaconditiononPtheta}The family
or Markov transition probabilities $\{P_{\theta},\theta\in\Theta\}$
is such that there exist:
\begin{longlist}[(2)]
\item[(1)]$V\dvtx \mathsf{X}\rightarrow{}[1,+\infty)$ and $\mathsf
{C}\subset\mathsf{X}$
such that $\sup_{x\in\mathsf{C}}V(x)<+\infty$,
\item[(2)]$a(\cdot),b(\cdot)\dvtx \Theta\rightarrow[0,+\infty)$ and $\iota
\in[0,1]$
\end{longlist}
such that for any $\theta,x\in\Theta\times\mathsf{X}$,
\[
P_{\theta}V(x)\leq \bigl[V(x)-a^{-1}(\theta)V^{\iota}(x)
\bigr]\mathbb{I}\{x\notin\mathsf{C}\}+b(\theta)\mathbb{I}\{x\in\mathsf {C}\}.
\]
\end{longlist}
For functions $v(\cdot)\dvtx \Theta\rightarrow\mathbb{R}$ we
define the level sets $\mathcal{V}_{M}:=\{\theta\in\Theta\dvtx v(\theta
)\leq M\}$
for all $M\geq0$, and for any set $A$ we will denote $A^{c}$ the
complement of $A$ in either $\Theta$ or $\mathsf{X}$. Notice that
assumption \ref{hypfixedthetaconditiononPtheta} implies that
$\inf_{\theta\in\Theta}a(\theta)>0$. The situations we are interested
in are those for which $\mathcal{A}_{M}^{c}\neq\emptyset$ for any
$M>0$. Hereafter it will be convenient to denote for any $\theta,x\in
\Theta\times\mathsf{X}$,
$\gamma\in(0,\gamma^{+}]$ and $f\dvtx \Theta\times\mathsf
{X}\rightarrow\mathbb{R}^{n_{f}}$
\[
P_{\theta,\gamma}f(\theta,x):=\int_{\mathsf{X}}P_{\theta
}(x,
\mathrm{d}y)f\bigl(\phi_{\gamma}(\theta,y),y\bigr).
\]

\renewcommand\thelonglist{(A3)}\begin{longlist}[(A3)]
\item\label{hypwdriftSA1}The family of
mappings $\{\phi_{\gamma}\dvtx \Theta\times\mathsf{X}\rightarrow\Theta,\gamma\in(0,\gamma^{+}]\}$
is such that there exists a Lyapunov function $w(\cdot)\dvtx \Theta
\rightarrow[1,\infty)$
such that [with $\{P_{\theta},\theta\in\Theta\}$, $V(\cdot
),\mathsf{C},a(\cdot),b(\cdot)$
and $\iota$ as in \ref{hypfixedthetaconditiononPtheta}]:
{\renewcommand\thelonglist{(\arabic{longlist})}
\renewcommand\labellonglist{\thelonglist}
\begin{longlist}[(5)]
\item\label{hypsomepropertiesofw}for any $M>0$, $\sup_{\theta\in
\mathcal{W}_{M}}a(\theta)<\infty$
and $\lim_{M\rightarrow\infty}\sup_{\theta\in\mathcal
{W}_{M}^{c}}b(\theta)/w(\theta)=0$,
\item there exists $\beta\in[0,1]$, $c(\cdot),d(\cdot)\dvtx \Theta
\rightarrow[0,\infty)$
and $e(\cdot)\dvtx \Theta\rightarrow[1,\infty)$ such that for all
$\gamma\in(0,\gamma^{+}]$
and $\theta,x\in\Theta\times\mathsf{X}$,
\[
P_{\theta,\gamma}w(\theta,x)\leq w(\theta)-\gamma w(\theta)\Delta
\bigl(V_{\beta}(\theta,x)\mathbb{I}\{x\notin\mathsf{C}\}+d(\theta )
\mathbb{I}\{x\in\mathsf{C}\} \bigr),
\]
where $V_{\beta}(\theta,x):=c(\theta)+V^{\beta}(x)/e(\theta)$ and
\item\label{hypenuwdivergeswhenadoes}$c(\cdot),d(\cdot)\dvtx \Theta
\rightarrow[0,\infty)$
are bounded, $\lim_{M\rightarrow\infty}\sup_{\theta\in\mathcal
{W}_{M}^{c}}[c(\theta)\vee d(\theta)]=0$
and for any $M>0$, $\sup_{\theta\in\mathcal{W}_{M}}e(\theta
)<\infty$,
\item\label{4}$\Delta(\cdot)\dvtx [0,\infty)\rightarrow\mathbb{R}$ is such that:
{\renewcommand\thelonglist{(\alph{longlist})}
\renewcommand\labellonglist{\thelonglist}
\begin{longlist}[(b)]
\item$\Delta(0)>0$ and it is continuous in a neighborhood of $0$,
\item\label{hypenugrowthofDelta}there exists $p_{\Delta}\in
(0,\iota/\beta]$
such that for all $M>0$ there exists $C_{\Delta,M}>0$, such that
for all $z\geq M$
\[
\bigl|\Delta(z)\bigr|\leq C_{\Delta,M}\times z{}^{p_{\Delta}},
\]
\end{longlist}}
\setcounter{longlist}{4}
\item\label{hypenurelationawV}for any $\epsilon>0$,
\[
\sup_{\theta,x\in\tilde{\mathcal{V}}_{\epsilon}}\frac{a(\theta
)w(\theta)e^{-p_{\Delta}}(\theta)}{V^{\iota-p_{\Delta}\beta
}(x)}<\infty,
\]
where $\tilde{\mathcal{V}}_{\epsilon}:=\{\theta,x\dvtx V^{\beta
}(x)/e(\theta)\geq\epsilon\}$.
\end{longlist}}
\end{longlist}

\begin{rem}
\label{remremarksabouthypmainresult}The conditions above may appear
abstract, but are motivated by the following concrete situations:
\begin{longlist}[(2)]
\item[(1)] The simultaneous fixed-$\theta$ drift conditions \ref{hypfixedthetaconditiononPtheta}
can be established in numerous situations of practical interest. Examples
are given in Section~\ref{171717171717}, where the transition
probabilities share the same invariant distribution, but it should
be pointed out that such drift conditions can also be established
in situations of interest where each transition kernel $P_{\theta}$
has its own invariant distribution $\pi_{\theta}$; this is the case
for example in the context of the stochastic approximation implementation
of the EM algorithm in \cite{andrieu:vihola:2011}. Other examples
can be found in \cite{benveniste} for algorithms used in the area
of digital communications, although the dependence on $\theta$ is
never used.
\item[(2)] Typically the function $\Delta(\cdot)$ in \ref{hypwdriftSA1}
will take the form of a polynomial, as a byproduct of a tractable
approximation of $w(\phi_{\gamma}(\theta,y))$ in terms of $w(\theta)$.
For example, in the situation where $\vartheta=\phi_{\gamma}(\theta,y)=\theta+\gamma H(\theta,y)$,
which corresponds to the standard stochastic approximation framework
(see Section~\ref{secmotivation}), a Taylor expansion of $w(\vartheta)$
around $\theta$ will lead to
\[
w(\vartheta)  \leq w(\theta)+\gamma \bigl\langle H(\theta,y),\nabla w(\theta)
\bigr\rangle+\tfrac{1}{2}\gamma^{2}\bar {w}''
\times\bigl|H(\theta,y)\bigr|^{2}
\]
whenever $\bar{w}'':=\sup_{\theta\in\Theta}|\nabla^{2}w(\theta
)|<\infty$.
With appropriate assumptions on $H(\theta,\allowbreak y)$, one can apply
$P_{\theta}$
to both sides of this inequality and hence obtain a drift condition
on $w(\cdot)$ of the form given in \ref{hypwdriftSA1}.
\item[(3)] The condition required on $p_{\Delta}\in(0,\iota/\beta${]}
can be
understood as being a tradeoff between the strength of the drift in
\ref{hypfixedthetaconditiononPtheta} and the strength of unfavorable
updates $\theta_{+}=\phi_{\gamma}(\theta,x_{+})$ such that
$w(\theta_{+})\gg w(\theta)$.
\end{longlist}

Hereafter for any $\varepsilon\in(0,\Delta(0))$ we will denote
\[
\Gamma_{\varepsilon}:= \bigl\{ \gamma,\bar{\gamma}\in\bigl(0,\gamma
^{+}\bigr]\dvtx \gamma^{-1}-\bar{\gamma}^{-1}<\Delta(0)-
\varepsilon \bigr\},
\]
where we omit the dependence on $\gamma^{+}$ for simplicity.
\end{rem}

%
\begin{thm}
\label{propestablishcompositedrift}Assume that $\{P_{\theta
},\theta\in\Theta\}$
and $\{\phi_{\gamma},\gamma\in(0,\gamma^{+}]\}$ satisfy \ref
{hypfixedthetaconditiononPtheta}
and~\ref{hypwdriftSA1}. Then for any $\varepsilon\in(0,\Delta(0))$
there exist $\lambda_{\ast}\in[1,\infty)$, $\delta,M_{*}\in
(0,+\infty)$
such that for any $\gamma,\bar{\gamma}\in\Gamma_{\varepsilon}$ and
$\theta,x\notin\mathcal{W}_{M_{\ast}}\times\mathsf{C}$,
%
\begin{equation}\label{eqmagicdrift}
\hspace*{20pt}P_{\theta,\gamma}\{\lambda_{*}V+w/\gamma\}(\theta,x)\leq\lambda
_{*}V(x)+w(\theta)/\bar{\gamma}-\delta\bigl[V^{\iota}(x)/a(
\theta )+w(\theta)\bigr].
\end{equation}
\end{thm}

%
\begin{cor}
\label{corcorollaryrecurrence}Let $\{\theta_{i},X_{i}\}$ be the
controlled Markov chain process as described in equation (\ref
{eqdefinitioncontrolledmc})
with for any $i\geq1$ $\phi_{i}(\theta_{0},x_{0},x_{1},\ldots,x_{i}):=\break \phi_{\gamma_{i}}(\theta_{i-1},x_{i})$
for a family $\{\phi_{\gamma}\dvtx \Theta\times\mathsf{X}\rightarrow
\Theta,\gamma\in(0,\gamma^{+}]\}$
and some real positive sequence $\{\gamma_{i}\}$. Assume further
that $\{P_{\theta},\theta\in\Theta\}$ and $\{\phi_{\gamma},\gamma
\in(0,\gamma^{+}]\}$
satisfy \ref{hypfixedthetaconditiononPtheta}, \ref{hypwdriftSA1}
and that $\{\gamma_{i}\}$ is such that
%
\begin{equation}\label{eqconditionongamma}
\bar{\varepsilon}:=\Delta(0)-\limsup_{i\rightarrow\infty}\bigl(\gamma
_{i+1}^{-1}-\gamma_{i}^{-1}\bigr)>0.
\end{equation}
Then, for any $\varepsilon\in(0,\bar{\varepsilon})$ there exists
$M_{*}$ as in Theorem~\ref{propestablishcompositedrift} such
that the set $\mathcal{W}_{M_{\ast}}\times\mathsf{C}$ is visited
infinitely often $\mathbb{P}_{\theta,x}$-a.s. by $\{\theta
_{i},X_{i}\}$.
\end{cor}

\begin{pf} Let $\varepsilon
\in(0,\bar{\varepsilon})$,
$\delta\in(0,1]$, $\lambda_{\ast}>1$ and $M_{\ast}>0$ be as in
Theorem~\ref{propestablishcompositedrift}, and define the family
of (Lyapunov) functions $\{W_{i}(\theta,x):=\lambda_{\ast
}V(x)+w(\theta)/\gamma_{i}\}$.
From the assumption on $\{\gamma_{i}\}$ there exists $i_{0}\in\mathbb{N}$
such that for any $i\geq i_{0}$ and $\theta,x\notin\mathcal
{W}_{M_{\ast}}\times\mathsf{C}$
\[
P_{\theta,\gamma_{i}}W_{i}(\theta,x)\leq W_{i-1}(\theta,x)-
\delta \bigl[V^{\iota}(x)/a(\theta)+w(\theta)\bigr].
\]
The result follows from Lemma~\ref{lemgeneral-drift-stability} since
$\inf_{\theta\in\Theta}w(\theta)>0$.
\end{pf}

\begin{rem}
One can notice that:
\begin{longlist}[(2)]
\item[(1)] in the case where $\gamma_{i}=c_{0}/(c_{1}+i)^{a}$ (\ref
{eqconditionongamma})
is satisfied for any $c_{0}>0$ and $a\in(0,1)$, and for $c_{0}<\Delta(0)$
when $a=1$;
\item[(2)] in the case where $\{\gamma_{i}=\gamma\leq\gamma^{+}\}$ is constant,
$W_{0}(\theta,x)=\lambda_{*}V(x)+w(\theta)$ and for any $\theta\in
\Theta$,
$a(\theta)\leq Cw^{\varkappa}(\theta)$ for $\varkappa>0$, then one
may show that for any $i\geq i_{0}$ and $\theta,x\notin\mathcal
{W}_{M_{\ast}}\times\mathsf{C}$
\[
P_{\theta,\gamma}W_{0}(\theta,x)\leq W_{0}(\theta,x)-
\delta 'W_{0}^{\iota/(1+\varkappa)}(\theta,x).
\]
Indeed, from a standard convexity inequality, for any $l\in(0,1]$,
%
\begin{equation}\label{eqmoreintermediateforcii}
l\frac{V^{\iota}(x)}{w(\theta)^{\varkappa}}+(1-l)w(\theta)\geq \biggl(\frac{V^{\iota}(x)}{w(\theta)^{\varkappa}}
\biggr)^{l}w^{1-l}(\theta),
\end{equation}
which, with the choice $\bar{l}=1/(1+\varkappa)$, leads to
\[
V^{\iota}(x)/w^{\varkappa}(\theta)+w(\theta)\geq\bar{l}V^{\iota
}(x)/w^{\varkappa}(
\theta)+(1-\bar{l})w(\theta)\geq V^{\iota
/(1+\varkappa)}(x).
\]
As a result we obtain [noting that, without loss of generality, one
can always take $C\geq1$ in the upper bound of $a(\cdot)$ above]
\begin{eqnarray*}
V^{\iota}(x)/a(\theta)+w(\theta) & \geq& C^{-1}2
\frac{1}{2}\bigl[V^{\iota
}(x)/w^{\varkappa}(\theta)+w(\theta)
\bigr]
\\
& \geq&\frac{1}{2C}\bigl[V^{\iota/(1+\varkappa)}(x)+w(\theta)\bigr]
\\
& \geq&\frac{1}{2C}\bigl[V^{\iota/(1+\varkappa)}(x)+w^{\iota
/(1+\varkappa)}(\theta)\bigr]
\\
& \geq& C^{-1}2^{-1-\iota/(1+\varkappa)} \bigl(V(x)+w(\theta)
\bigr)^{\iota/(1+\varkappa)},
\end{eqnarray*}
and we conclude. This suggests the possibility of precisely characterizing
the return times to $\mathcal{W}_{M_{\ast}}\times\mathsf{C}$, as this
form of drift condition is known to lead to the existence of polynomial
moments of return times.
\end{longlist}
\end{rem}

\begin{pf*}{Proof of Theorem~\ref{propestablishcompositedrift}} Choose $\varepsilon\in
(0,\Delta(0))$
and $\epsilon_{-}>0$ such that for any $|z|\leq\epsilon_{-}$,
$|\Delta(0)-\Delta (z )|\leq\varepsilon/2$.
This implies
%
\begin{equation}\label{eqchoiceM0implies}
\hspace*{21pt}\sup_{\gamma,\bar{\gamma}\in\Gamma_{\varepsilon}}\bigl(\gamma ^{-1}-\bar{
\gamma}^{-1}\bigr)-\inf_{\{z:|z|\leq\epsilon_{-}\}}\Delta (z )\leq
\Delta(0)-\varepsilon+\varepsilon/2-\Delta (0)=-\varepsilon/2.
\end{equation}
Now let $M_{0}\geq0$ be such that $\sup_{\theta\in\mathcal
{W}_{M_{0}}^{c}}d(\theta)\leq\epsilon_{-}$
and $\sup_{\theta\in\mathcal{W}_{M_{0}}^{c}}c(\theta)\leq\epsilon_{-}/2$.
From \ref{hypfixedthetaconditiononPtheta} and \ref{hypwdriftSA1}
we have for $(\theta,x)\in\Theta\times\mathsf{X}$ and $\lambda\in
(0,\infty)$
\begin{eqnarray*}
&&P_{\theta,\gamma}\{\lambda V+w/\gamma\}(\theta,x)\\
&&\qquad \leq\lambda
\bigl[V(x)-a^{-1}(\theta)V^{\iota}(x) \bigr]\mathbb{I}\{x\notin
\mathsf {C}\}+\lambda b(\theta)\mathbb{I}\{x\in\mathsf{C}\}
\\
&&\qquad \quad {}+w(\theta)/\gamma-w(\theta)\Delta \bigl(V_{\beta}(\theta,x)\mathbb{I}\{x
\notin\mathsf{C}\}+d(\theta)\mathbb{I}\{x\in \mathsf{C}\} \bigr).
\end{eqnarray*}
Note that for $(\theta,x)\in\mathcal{W}_{M_{0}}\times\mathsf{C}^{c}$,
$V_{\beta}(\theta,x)\geq\bar{M}_{0}^{-1}:=1/\sup_{\theta\in
\mathcal{W}_{M_{0}}}e(\theta)>0$
from \ref{hypwdriftSA1}\ref{hypenuwdivergeswhenadoes},
and therefore from \ref{hypwdriftSA1}\ref{4}\ref{hypenugrowthofDelta}
\[
w(\theta)\bigl\llvert \Delta\bigl(V_{\beta}(\theta,x)\bigr)\bigr\rrvert
\leq C_{\Delta,\bar{M}_{0}^{-1}}M_{0}\sup_{\theta\in\Theta}
\bigl(e^{-1}(\theta )+c(\theta)\bigr)^{p_{\Delta}}\times
V^{p_{\Delta}\beta}(x),
\]
and $\sup_{\theta\in\mathcal{W}_{M_{0}}}\Delta(d(\theta))<\infty$
as $\Delta(\cdot)$ is bounded on compact sets. In addition  $\sup_{\theta\in
\Theta}d(\theta)<\infty$.
Let now
\[
C'_{\Delta,M_{0}}:=\Bigl[M_{0}C_{\Delta,\bar{M}_{0}^{-1}}\sup
_{\theta
\in\Theta}\bigl(e^{-1}(\theta)+c(\theta)
\bigr)^{p_{\Delta}}\Bigr]\vee\Bigl[M_{0}\sup_{\theta\in\mathcal{W}_{M_{0}}}
\Delta\bigl(d(\theta)\bigr)\Bigr]<\infty,
\]
notice that $\iota\geq p_{\Delta}\beta$, and recall that $V\geq1$.
Then we have for $\gamma,\bar{\gamma}\in\Gamma_{\varepsilon}$
\[
P_{\theta,\gamma}\{\lambda V+w/\gamma\}(\theta,x)\leq\lambda V(x)+w(\theta)/
\bar{\gamma}+\Lambda(\theta,x),
\]
with
%
\begin{eqnarray}\label{eqdefofLambda}
\Lambda(\theta,x)&:=&-\lambda V(x)+\lambda \bigl[V(x)-a^{-1}(\theta
)V^{\iota}(x) \bigr]\mathbb{I}\{x\notin\mathsf{C}\}\nonumber\\
&&{}+\lambda b(\theta)
\mathbb{I}\{x\in\mathsf{C}\}
+\bigl(\Delta(0)-\varepsilon\bigr)w(\theta)\nonumber\\[-8pt]\\[-8pt]
&&{}-\Delta \bigl(V_{\beta}(\theta,x)\mathbb{I}\{x\notin\mathsf{C}\}
+d(\theta)\mathbb{I}\{x\in \mathsf{C}\}
\bigr)w(\theta) \mathbb{I}\bigl\{\theta\in\mathcal {W}_{M_{0}}^{c}
\bigr\}\nonumber\\
&&{}+C'_{\Delta,M_{0}}V^{\iota}(x) \mathbb{I}\{\theta\in
\mathcal {W}_{M_{0}}\}.\nonumber
\end{eqnarray}
It will be convenient below to refer to the following inequality:
%
\begin{equation}
\Lambda(\theta,x)\leq-\tilde{\delta}\bigl[V^{\iota}(x)/a(\theta )+w(
\theta)\bigr],\label{eqconditioncasenotinC}
\end{equation}
for $(\theta,x)\notin\mathcal{W}_{\tilde{M}}\times\mathsf{S}$ and
various instantiations of $\tilde{\delta},\tilde{M},\lambda>0$ and
$\mathsf{S}\subset\mathsf{X}$. Our ultimate aim is to prove that
under the stated assumptions there exist $\delta,\lambda_{\ast}\in
(0,+\infty)$
and $M_{\ast}\geq M_{0}$ such that \eqref{eqconditioncasenotinC}
holds for $(\theta,x)\in(\mathcal{W}_{M_{*}}\times\mathsf{C})^{c}$.
For any $M\geq M_{0}$, we use the following partition:
\[
(\mathcal{W}_{M}\times\mathsf{C})^{c}=\bigl(
\mathcal{W}_{M_{0}}\times \mathsf{C}^{c}\bigr)\cup\bigl(
\mathcal{W}_{M_{0}}^{c}\times\mathsf {C}^{c}\bigr)
\cup\bigl(\mathcal{W}_{M}^{c}\times\mathsf{C}\bigr),
\]
which leads us to consider three cases, (a), (b) and (c), from left
to right.

(a) For $(\theta,x)\in\mathcal{W}_{M_{0}}\times\mathsf{C}^{c}$
and any $\lambda>0$, we have
\begin{eqnarray*}
\Lambda(\theta,x) & \leq&\bigl[\Delta(0)-\varepsilon\bigr]w(\theta)-\lambda
V^{\iota}(x)/a(\theta)+C'_{\Delta,M_{0}}V^{\iota}(x)
\\
& \leq&\bigl[\Delta(0)-\varepsilon\bigr]\sup_{\theta\in\mathcal
{W}_{M_{0}}}w(
\theta)+V^{\iota}(x) \Bigl[C'_{\Delta,M_{0}}-\lambda \bigl/\sup
_{\vartheta\in\mathcal{W}_{M_{0}}}a(\vartheta) \Bigr],
\end{eqnarray*}
where we note that $\sup_{\vartheta\in\mathcal
{W}_{M_{0}}}a(\vartheta)<\infty$
from \ref{hypwdriftSA1}\ref{hypsomepropertiesofw}. Now,
from our choice of $M_{0}$ and since $V\geq1$ and $\inf_{\vartheta
\in\Theta}a(\vartheta)>0$,
we conclude about the existence of $\lambda_{a},\delta_{a}>0$ such
that for all $\lambda\geq\lambda_{a}$, $(\theta,x)\in\mathcal
{W}_{M_{0}}\times\mathcal{\mathsf{C}}^{c}$
\begin{eqnarray*}
\Lambda(\theta,x)&\leq&\bigl[\Delta(0)-\varepsilon\bigr]M_{0}+V^{\iota
}(x)
\Bigl[C'_{\Delta,M_{0}}-\lambda\bigl/\sup_{\vartheta\in\mathcal
{W}_{M_{0}}}a(
\vartheta) \Bigr]\\
&\leq&-\delta_{a}\bigl[V^{\iota
}(x)/a(\theta)+w(
\theta)\bigr].
\end{eqnarray*}
Therefore (\ref{eqconditioncasenotinC}) is satisfied with $\tilde
{M}=M_{0}$,
any $\lambda\geq\lambda_{a}$ and $\tilde{\delta}=\delta_{a}$.

(b) For $(\theta,x)\in\mathcal{W}_{M_{0}}^{c}\times
\mathsf{\mathsf{C}}^{c}$
and any $\lambda>0$, we have
\[
\Lambda(\theta,x)\leq-\lambda V^{\iota}(x)/a(\theta)+ \bigl[\Delta (0)-
\varepsilon-\Delta \bigl(V_{\beta}(\theta,x) \bigr) \bigr]w(\theta),
\]
and we seek to show that there exists $\lambda_{b}=\lambda_{i}\vee
\lambda_{ii}$
and $\delta_{b}=\delta_{i}\wedge\delta_{ii}>0$ (where $\lambda
_{i},\lambda_{ii}>0$
and $\delta_{i},\delta_{ii}>0$ are given below in the proof) such
that for all $\lambda\geq\lambda_{b}$ and $(\theta,x)\in\mathcal
{W}_{M_{0}}^{c}\times\mathsf{C}^{c}$,
(\ref{eqconditioncasenotinC}) is satisfied with $\tilde{\delta
}=\delta_{b}$.
In what follows we will use the following intermediate results. From
\ref{hypwdriftSA1}\ref{hypenurelationawV} we have that
for our earlier choice of $\epsilon_{-}$ and any $(\theta,x)\in
\mathcal{W}_{M_{0}}^{c}\times\mathsf{C}^{c}$,
the condition
\[
V_{\beta}(\theta,x)=\frac{V^{\beta}(x)}{e(\theta)}+c(\theta)\geq
\epsilon_{-}\quad  \mbox{implies}\quad  \frac{V^{\beta}(x)}{e(\theta)}\geq
\epsilon_{-}-\sup_{\vartheta\in\mathcal{W}_{M_{0}^{c}}}c(\vartheta )\geq
\epsilon_{-}/2,
\]
and therefore that for $q\in\{0,p_{\Delta}\}$
%
\begin{equation}\label{eqintermediateresultsc-i}
\sup_{(\mathcal{W}_{M_{0}}^{c}\times\mathsf{C}^{c})\cap\{\theta,x:V_{\beta}(\theta,x)\geq\epsilon_{-}\}}\frac{a(\theta)w(\theta
)e^{-q}(\theta)}{V^{\iota-q\beta}(x)}\leq C_{\epsilon_{-}}<\infty.
\end{equation}
Indeed the case $q=p_{\Delta}$ is true by assumption and for $V^{\beta
}(x)/e(\theta)\geq\epsilon_{-}/2$
\[
\frac{a(\theta)w(\theta)e^{-p_{\Delta}}(\theta)}{V^{\iota
-p_{\Delta}\beta}(x)}=\frac{a(\theta)w(\theta)}{V^{\iota
}(x)} \biggl(\frac{V^{\beta}(x)}{e(\theta)}
\biggr)^{p_{\Delta}}\geq \frac{a(\theta)w(\theta)}{V^{\iota}(x)} (\epsilon_{-}/2
)^{p_{\Delta}}
\]
from which we conclude. We now partition $\mathcal
{W}_{M_{0}}^{c}\times\mathsf{C}^{c}$
by considering the two following subsets:

(i) From our choice of $M_{0}$ and $\epsilon_{-}$
and (\ref{eqchoiceM0implies}), we deduce that on the subset
$(\mathcal{W}_{M_{0}}^{c}\times\mathsf{C}^{c})\cap\{\theta,x\dvtx V_{\beta}(\theta,x)<\epsilon_{-}\}$
\[
\Delta(0)-\varepsilon-\Delta \bigl(V_{\beta}(\theta,x) \bigr)\leq -
\varepsilon/2,
\]
and consequently
\[
\Lambda(\theta,x)\leq-\lambda V^{\iota}(x)/a(\theta)-w(\theta )
\varepsilon/2,
\]
and we conclude about the existence of $\lambda_{i}, \delta_{i}>0$
such that (\ref{eqconditioncasenotinC}) holds for any $\lambda
\geq\lambda_{i}$
and $\tilde{\delta}=\delta_{i}$.

(ii) By our assumption on $\Delta(\cdot)$ there exists
$C'_{\Delta,\epsilon_{-}}>0$
such that for any $z\geq\epsilon_{-}$
\[
\Delta(0)-\varepsilon-\Delta(z)\leq C'_{\Delta,\epsilon
_{-}}z^{p_{\Delta}}.
\]
Consequently we deduce that on $(\mathcal{W}_{M_{0}}^{c}\times\mathsf
{C}^{c})\cap\{\theta,x\dvtx V_{\beta}(\theta,x)\geq\epsilon_{-}\}$
\begin{eqnarray*}
\Lambda(\theta,x) & \leq&\frac{V^{\iota}(x)}{a(\theta)} \biggl[C'_{\Delta,\epsilon_{-}}
\frac{a(\theta)w(\theta)}{V^{\iota
}(x)} \biggl(\frac{V^{\beta}(x)}{e(\theta)}+c(\theta) \biggr)^{p_{\Delta}}-
\lambda \biggr]
\\
& \leq&\frac{V^{\iota}(x)}{a(\theta)} \biggl[2^{p_{\Delta
}}C'_{\Delta,\epsilon_{-}}
\frac{a(\theta)w(\theta)e(\theta
)^{-p_{\Delta}}}{V^{\iota-p_{\Delta}\beta}(x)}-\lambda \biggr].
\end{eqnarray*}
We now choose $\lambda>2^{p_{\Delta}}C'_{\Delta,\epsilon
_{-}}C_{\epsilon_{-}}$,
and from (\ref{eqintermediateresultsc-i}) with $q=0$ we have
$V^{\iota}(x)/\allowbreak a(\theta)\geq (\epsilon_{-}/2 )^{p_{\Delta
}}C_{\epsilon_{-}}^{-1}w(\theta)$
and therefore
\[
\Lambda(\theta,x)\leq\frac{1}{2} \biggl[\frac{V^{\iota
}(x)}{a(\theta)}+w(\theta) (
\epsilon_{-}/2 )^{p_{\Delta
}}C_{\epsilon_{-}}^{-1}
\biggr] \bigl(2^{p_{\Delta}}C'_{\Delta,\epsilon_{-}}C_{\epsilon_{-}}-
\lambda \bigr).
\]
We conclude about the existence of $\lambda_{ii},\delta_{ii}>0$
such that (\ref{eqconditioncasenotinC}) holds for any $\lambda
\geq\lambda_{ii}$
and $\tilde{\delta}=\delta_{ii}$.

(c) First, we note from our choice of $M_{0}$, (\ref
{eqchoiceM0implies})
and $(\ref{eqdefofLambda})$ that for any $\theta,x\in\mathcal
{W}_{M_{0}}^{c}\times\mathsf{C}$
and $\lambda>0$ the function $\Lambda(\theta,x)$ is upper bounded
by
\begin{eqnarray*}
\Lambda(\theta,x) & \leq&-\lambda V(x)+\lambda b(\theta)+ \bigl[\Delta(0)-
\varepsilon-\Delta \bigl(d(\theta) \bigr) \bigr]w(\theta)
\\
& \leq&-\lambda V(x)+\bigl[\lambda b(\theta)/w(\theta)-\varepsilon /2\bigr]w(
\theta).
\end{eqnarray*}
We now show that for any $\lambda\in(0,+\infty)$ there exist
$M_{\lambda}\in{}[ M_{0},+\infty)$
and $\delta_{\lambda}\in(0,+\infty)$ such that for all $\theta,x\in\mathcal{W}_{M_{\lambda}}^{c}\times\mathsf{C}$,
(\ref{eqconditioncasenotinC}) is satisfied with $\tilde{\delta
}=\delta_{\lambda}$,
a function of $\lambda$. From our last inequality and since $\iota\in[0,1]$
and $V\geq1$, for any $M\geq M_{0}$ and $\theta,x\in\mathcal
{W}_{M}^{c}\times\mathsf{C}$
\[
\Lambda(\theta,x)  \leq-\lambda{ \Bigl[\inf_{\vartheta
\in\Theta}a(
\vartheta) \Bigr]} V^{\iota}(x)/a(\theta)+\Bigl[\lambda \sup
_{\vartheta\in\mathcal{W}_{M}^{c}}b(\vartheta)/w(\vartheta )-\varepsilon/2\Bigr]w(\theta).
\]
We conclude about the existence of $M_{\lambda}$ and $\delta_{\lambda}$
as above from our assumption on $b(\cdot)$.

We now conclude by letting $\lambda_{\ast}\geq\lambda_{a,b}=\lambda
_{a}\vee\lambda_{b}$, $M_{\ast}\geq M_{\lambda_{a}\vee\lambda_{b}}$
and $\delta=\delta_{\lambda_{a}\vee\lambda_{b}}\wedge\delta
_{a}\wedge\delta_{b}$.
\end{pf*}

\section{The central role of stability in the context of stochastic
approximation
with Markovian dynamic}
\label{secmotivation}

In this section we illustrate the central role played by the form
of stability considered in this paper to establish that some controlled
Markov chains of practical relevance possess some desired properties.
We focus on a particular class of controlled Markov chains driven
by a so-called stochastic approximation recursion (also known as the
Robbins--Monro algorithm). The motivation for such algorithms, described
below, is to find the roots of the function $h(\cdot)\dvtx \Theta
\rightarrow\mathbb{R}^{n_{\theta}}$
\[
h(\theta):=\int_{\mathsf{X}}H(\theta,x)\pi_{\theta}(
\mathrm{d}x),
\]
for a family of functions $\{H(\theta,x)\dvtx \Theta\times\mathsf
{X}\rightarrow\Theta\}$
and a family of probability distributions $\{\pi_{\theta},\theta\in
\Theta\}$
defined on some space $\mathsf{X}\times\mathcal{B}(\mathsf{X})$.
This is a ubiquitous problem in statistics, engineering and computer
science. These roots are rarely available analytically and a way of
finding them numerically consists of considering the following controlled
Markov chain on $ ((\Theta\times\mathsf{X})^{\mathbb
{N}},(\mathcal{B}(\Theta)\otimes\mathcal{B}(\mathsf{X)})^{\otimes
\mathbb{N}} )$,
initialized at some $(\theta_{0},X_{0})=(\theta,x)\in\Theta\times
\mathsf{X}$
and defined recursively for a sequence of stepsizes $\{\gamma_{i}\}$
for $i\geq0$,
%
\begin{eqnarray}\label
{eqSAwithMarkovianDynamic}
&X_{i+1}|\mathcal{F}_{i}\sim P_{\theta_{i}}(X_{i},
\cdot),&
\nonumber
\\[-8pt]\\[-8pt]
&\theta_{i+1}=\theta_{i}+\gamma_{i+1}H(
\theta_{i},X_{i+1}),&\nonumber
\end{eqnarray}
where $\{P_{\theta},\theta\in\Theta\}$ (for some set $\Theta
\subset\mathbb{R}^{n_{\theta}}$)
is a family of Markov transition probabilities such that for each
$\theta\in\Theta$, $P_{\theta}$ leaves $\pi_{\theta}$ invariant,
that is, is such that $\pi_{\theta}P_{\theta}=\pi_{\theta}$. The
rational for this recursion is as follows. Let us first rewrite the
Robbins--Monro recursion as
\[
\theta_{i+1}=\theta_{i}+\gamma_{i+1} \bigl[h(
\theta_{i})+\xi _{i+1} \bigr],
\]
where $\{\xi_{i+1}=H(\theta_{i},X_{i+1})-h(\theta_{i})\}$, which
is traditionally refered to as the ``noise.'' Then
$\{\theta_{i}\}$ can be thought as being a noisy version of the sequence
$\{\bar{\theta}_{i}\}$ defined as $\bar{\theta}_{i+1}=\bar{\theta
}_{i}+\gamma_{i+1}h(\bar{\theta}_{i})$,
and it is believable that the properties of $\{\theta_{i}\}$ are
closely related to those of the noiseless sequence $\{\bar{\theta
}_{i}\}$
provided the average effect of the noise on this sequence is negligible.
This requires some form of averaging, or ergodicity, property on $\{\xi
_{i}\}$.

The convergence of such sequences has been well studied by various
authors, starting with the seminal work of \cite{benveniste}, under
various assumptions on all the quantities involved. A crucial step
of such convergence analyses, however, consists of assuming that the
sequence $\{\theta_{i}\}$ remains bounded in a compact set of $\Theta$
with probability one. This problem has traditionally been either ignored
or circumvented by means of modifications of the recursion (\ref
{eqSAwithMarkovianDynamic}).
Indeed, one of the major difficulties specific to the Markovian dynamic
scenario is that $\{\theta_{i}\}$ governs the ergodicity of $\{X_{i}\}$
(and hence $\{\xi_{i}\}$) and that stability properties of $\{\theta
_{i}\}$
relying on those of $\{\bar{\theta}_{i}\}$ require good ergodicity
properties which might vanish whenever $\{\theta_{i}\}$ approaches
a set $\partial\Theta$ away from the zeroes of $h(\theta)$, resulting
in instability. Most existing results rely on modifications of the
updates $\{\phi_{\gamma}\}$ designed to ensure a form of ergodicity
of $\{\xi_{i}\}$ which in turn ensures that $\{\theta_{i}\}$ inherits
the stability properties of $\{\bar{\theta}_{i}\}$. The only known
results we are aware of where stability is established for (\ref
{eqSAwithMarkovianDynamic})
without any modification are \cite{benveniste}, Part II, Section~1.9,
where assumption (1.9.3) may not be satisfied in numerous cases of
interest or directly verifiable, and \cite{saksman:vihola:2008} in
a particular scenario.

The approach we follow here is significantly different from that developed
in the aforementioned works and consists of dividing the difficult
problem of proving boundedness away from $\partial\Theta$ into two
simpler tasks. First using the results of Sections~\ref{secCompound-Lyapunov-functions}
or \ref{1313131313}, one may establish that the sequence $\{
\theta_{i},X_{i}\}$
visits some set $\mathcal{W}\times C\subset\Theta\times\mathsf{X}$
infinitely often $\mathbb{P}_{\theta,x}$-a.s., a set which has the
particularity that transition probabilities $\{P_{\theta},\theta\in
\mathcal{W}\}$
have uniformly good ergodicity properties. Then, using these facts,
one can show that $\{\theta_{i}\}$ follows the trajectories of the
deterministic recursion $\{\bar{\theta}_{i}\}$ more and more accurately
at each visit of $\mathcal{W}\times C$, and eventually remains in
a set only slightly larger than $\mathcal{W}$ provided the deterministic
sequence is itself stable. The advantage of our approach is that instead
of aiming to establish ergodicity properties of $\{\xi_{i}\}$ in
worse case scenarios for the sequence $\{\theta_{i}\}$, we decouple
the analysis of the behavior of $\{\theta_{i}\}$ when it approaches
$\partial\Theta$ from the study of the ergodicity properties of $\{
\xi_{i}\}$,
which need to be studied for ``reasonable'' values of $\theta$
only. Before stating the main result of this section we state our
assumptions.
\renewcommand\thelonglist{(A4)}
\begin{longlist}[(A4)]
\item \label{hypthetastaysinTheta}Let $\{
H(\theta,x)\}$,
$\{\gamma_{i}\}$, $\{\pi_{\theta}\}$ and $\{P_{\theta}\}$ be as
above. We assume that:
\begin{enumerate}[(3)]
\item[(1)] there exists $\gamma_{+}>0$ such that:

\begin{longlist}[(b)]
\item[(a)]$\gamma:=\{\gamma_{i}\}\subset{}[0,\gamma^{+}]^{\mathbb{N}}$,
\item[(b)] for any $\theta,x\in\Theta\times\mathsf{X}$ and $\gamma\in
{}[0,\gamma^{+}]$
\[
\theta+\gamma H(\theta,x)\in\Theta,
\]
\end{longlist}
\item[(2)]$H\dvtx \Theta\times\mathsf{X}\rightarrow\mathbb{R}^{n_{\theta
}}$ is such
that for any $\theta\in\Theta$, $\int_{\mathsf{X}}|H(\theta,x)|\pi_{\theta}(\mathrm{d}x)<+\infty$,
\item[(3)] and for any $\theta\in\Theta$, $\pi_{\theta}P_{\theta}=\pi
_{\theta}$.
\end{enumerate}
\end{longlist}

A practical technique to prove the boundedness of the noiseless sequence
consists, whenever possible, of determining a Lyapunov function
$w\dvtx \Theta\rightarrow{}[0,\infty)$
such that $ \langle\nabla w(\theta),h(\theta) \rangle
\leq0$
away from the roots of $h(\theta)$, where $\nabla w$ denotes the
gradient of $w$ with respect to $\theta$, and for $u,v\in\mathbb{R}^{n}$,
$ \langle\cdot,\cdot \rangle$ is their Euclidean inner
product (we will later on also use the notation $|v|=\sqrt{
\langle v,v \rangle}$
to denote the Euclidean norm of $v$). Note that although we use here
the same symbol $w$ as in Sections~\ref{121212121} and \ref
{1313131313},
the Lyapunov function below might be different.
\renewcommand\thelonglist{(A5)}
\begin{longlist}[(A5)]
\item \label{hypALyapunovFunction} $\Theta
$ is an
open subset of $\mathbb{R}^{n_{\theta}}$, $h\dvtx \Theta\rightarrow
\mathbb{R}^{n_{\theta}}$
is continuous and there exists a continuously differentiable function
$w\dvtx \Theta\rightarrow{}[0,\infty)$ such that:
\renewcommand\theenumi{(\arabic{enumi})}
\renewcommand\labelenumi{\theenumi}
\begin{enumerate}[(3)]
\item\label{enumcompacitycalL} there exists $M_{0}>0$ such that
\[
\mathcal{L}:= \bigl\{ \theta\in\Theta, \bigl\langle\nabla w(\theta ),h(\theta)
\bigr\rangle=0 \bigr\} \subset \bigl\{ \theta\in \Theta,w(\theta)<M_{0}
\bigr\},
\]

\item\label{enumcompacitylevelsets} there exists $M_{1}\in
(M_{0},\infty]$
such that $\mathcal{W}_{M_{1}}$ is a compact set,
\item\label{enumgloballyapunovfunction} for any $\theta\in\Theta
\setminus\mathcal{L}$,
$ \langle\nabla w(\theta),h(\theta) \rangle<0$.
\end{enumerate}
\end{longlist}

We now introduce some additional notation needed to describe the ergodicity
properties of $\{\xi_{i}\}$ every time the sequence $\{\theta
_{i},X_{i}\}$
visits some set $\mathcal{W}\times C$. More precisely, consider the
stochastic processes $\{\vartheta_{i},\mathfrak{X}_{i}\}$ defined
on $ ((\Theta\times\mathsf{X})^{\mathbb{N}},(\mathcal
{B}(\Theta)\otimes\mathcal{B}(\mathsf{X)})^{\otimes\mathbb
{N}} )$
which use the stepsize sequence $\boldsymbol{\gamma}^{\leftarrow l}:=\{
\gamma_{i+l},i\geq0\}$
for some $l\geq0$, initialized with $\vartheta_{0},\mathfrak
{X}_{0}\in\Theta\times\mathsf{X}$
and such that for $i\geq0$,\looseness=-1
%
\begin{eqnarray}
&\mathfrak{X}_{i+1}|(\vartheta_{0},\mathfrak{X}_{0},
\mathfrak {X}_{1},\ldots,\mathfrak{X}_{i})\sim
P_{\vartheta_{i}}(\mathfrak {X}_{i},\cdot),&
\nonumber
\\[-8pt]\\[-8pt]
&\vartheta_{i+1}=\vartheta_{i}+\gamma_{i+1+l}H(
\vartheta _{i},\mathfrak{X}_{i+1}).&\nonumber
\end{eqnarray}
In order to take the shift in the stepsize sequence into account, we
denote by $\bar{\mathbb{P}}_{\theta,x}^{\boldsymbol{\gamma}^{\leftarrow l}}$ and
$\bar{\mathbb{E}}_{\theta,x}^{\boldsymbol{\gamma}^{\leftarrow l}}$ the associated
probability distribution and expectation operator for $\vartheta
_{0}=\theta\in\Theta$
and $\mathfrak{X}_{0}=x\in\mathsf{X}$, and point out that in
contrast to $\mathbb{P}_{\theta,x}$ defined earlier for $\{\theta
_{i},X_{i}\}$,
the notational dependence on $\boldsymbol{\gamma}^{\leftarrow l}$ for
$l\geq0$ is here crucial. For any $M>0$ we define the exit time
from the level set $\mathcal{W}_{M}$, $\sigma(\mathcal{W}_{M}):=\inf
\{ k\geq1\dvtx \vartheta_{k}\notin\mathcal{W}_{M} \} $
(with the standard convention that $\inf\{\emptyset\}=+\infty$), and
for any $j\geq1$ we define $\varsigma_{j}:=H(\vartheta
_{j-1},\mathfrak{X}_{j})-h(\vartheta_{j-1})$.
We extend the result of \cite{andrieu:moulines:priouret:2005}, Proposition~5.2
(see also \cite{andrieu:moulines:2006}) in order to establish the
following result.

%
\begin{thm}
\label{thlocalergodicityimpliesstability}Assume \textup{\ref
{hypthetastaysinTheta}}
and \textup{\ref{hypALyapunovFunction}}, that $\{\gamma_{i}\}$ is such
that $\limsup_{i\rightarrow\infty}\gamma_{i}=0$, and let $M\in
(M_{0},M_{1}]$.
Assume that there exists $\mathsf{C}\subset\mathsf{X}$ such that:
\begin{longlist}[(2)]
\item[(1)] for any $\epsilon>0$,
%
\begin{equation}\label{eqconditiononcumalivenoisesums}
\quad\limsup_{l\rightarrow\infty}\sup_{\theta,x\in\mathcal
{W}_{M_{0}}\times\mathsf{C}}\bar{
\mathbb{P}}_{\theta,x}^{\boldsymbol{\gamma}
^{\leftarrow l}} \Biggl(\sup_{k\geq1}
\mathbb{I}\bigl\{\sigma(\mathcal {W}_{M})\geq k\bigr\}\Biggl\llvert
\sum_{j=1}^{k}\gamma_{j+l}
\varsigma _{j}\Biggr\rrvert >\epsilon \Biggr)<1;
\end{equation}

\item[(2)] for any $\theta,x\in\Theta\times\mathsf{X}$,
\[
\mathbb{P}_{\theta,x} \Biggl(\bigcap_{k=1}^{\infty}
\bigcup_{i=k}^{\infty
}\bigl\{(\theta_{i},X_{i})
\in\mathcal{W}_{M_{0}}\times\mathsf{C}\bigr\} \Biggr)=1,
\]
that is, $\{\theta_{i},X_{i}\}$ defined by equation (\ref
{eqSAwithMarkovianDynamic})
visits $\mathcal{W}_{M_{0}}\times\mathsf{C}$ infinitely often
$\mathbb{P}_{\theta,x}$-a.s.
\end{longlist}
Then the sequence $\{\theta_{i}\}$, as defined by equation (\ref
{eqSAwithMarkovianDynamic}),
is such that
\[
\mathbb{P}_{\theta,x} \Biggl(\bigcup_{k=1}^{\infty}
\bigcap_{i=k}^{\infty
}\{\theta_{i}\in
\mathcal{W}_{M}\} \Biggr)=1,
\]
that is, $\{\theta_{i}\}$ eventually remains in $\mathcal{W}_{M}$,
$\mathbb{P}_{\theta,x}$-a.s.
\end{thm}

\begin{rem}
Proving equation (\ref{eqconditiononcumalivenoisesums}) is now rather
well understood in general scenarios as soon as some form of local
(in $\theta$) uniform ergodicity of $\{P_{\theta}\}$ is satisfied
and can be checked in practice; see \cite{andrieu:moulines:2006}
and \cite{andrieu:moulines:priouret:2005}, for example, and the recent
results in \cite{atchade:fort:2008}. In the present paper, we rather
focus on finding verifiable conditions on $\{\gamma_{i}\},\{H(\theta,x)\}$
and $\{P_{\theta}\}$ which ensure that $\{\theta_{i},X_{i}\}$, as
defined by equation (\ref{eqSAwithMarkovianDynamic}), visits
$\mathcal{W}_{M_{0}}\times C$
infinitely often $\mathbb{P}_{\theta,x}$-a.s., which in combination
with the aforementioned existing results will allow us to conclude
about the stability of a large class of controlled MCMC algorithms.
\end{rem}

\begin{pf*}{Proof of Theorem~\ref{thlocalergodicityimpliesstability}}
For $M\in(M_{0},M_{1}]$ we let $\delta_{0}>0$ and $\lambda_{0}>0$
be as in Theorem~\ref{theokeyproposition} from \cite{andrieu:moulines:priouret:2005}, Proposition~5.2,
given in  Appendix \ref{secAppendix-for-Section3} for convenience. We consider the sequence $\{
T_{i},i\geq1\}$
of successive return times to $\mathcal{W}_{M_{0}}\times C$ ``separated
by at least an exit from $\mathcal{W}_{M}$,'' formally
defined for $i\geq0$ as
\[
T_{i+1}=\inf\bigl\{j\geq T_{i}+1\dvtx \exists l\in
\{T_{i}+1,\ldots,j\}/\theta _{l}\notin\mathcal{W}_{M}
\mbox{ and } (\theta_{j},X_{j})\in \mathcal{W}_{M_{0}}
\times C\bigr\},
\]
with the conventions $T_{0}=0$ and $\inf\{\emptyset\}=+\infty$.
It will be useful below to note that for any $i\geq1$, $T_{i}\geq i$.
Let $n_{0}\in\mathbb{N}$ be such that $\sup_{k\geq n_{0}}\gamma
_{k}\leq\lambda_{0}$.
We first show that for any $\theta,x\in\Theta\times\mathsf{X}$,
%
\begin{equation}\label{eqsupTkinfinite}
\mathbb{P}_{\theta,x} \biggl({\bigcup_{k\geq1}}
\{ T_{k}=+\infty \} \biggr)=1,
\end{equation}
and to achieve this, we establish a bound on $\sup_{\theta,x\in
\Theta\times\mathsf{X}}\mathbb{P}_{\theta,x} (T_{n}<+\infty
)$
for $n\geq n_{0}$. Notice that from the strong Markov property, for
any $\theta,x\in\Theta\times\mathsf{X}$ and $l\geq n_{0}$
\[
\mathbb{P}_{\theta,x} (T_{l+1}<+\infty )=\mathbb
{E}_{\theta,x} \bigl(\mathbb{I}\{T_{l}<+\infty\}
\mathbb{P}_{\theta
_{T_{l}},X_{T_{l}}} (T_{1}<+\infty ) \bigr).
\]
In addition, for any $\theta,x\in\Theta\times\mathsf{X}$, we have
\begin{eqnarray*}
&&\mathbb{I}\{T_{l}<+\infty\}\mathbb{P}_{\theta
_{T_{l}},X_{T_{l}}}
(T_{1}<+\infty )\\
&&\qquad \leq\mathbb{I}\{ T_{l}<+\infty\}\bar{
\mathbb{P}}_{\theta_{T_{l}},X_{T_{l}}}^{\boldsymbol{\gamma}^{\leftarrow T_{l}}} \bigl(\sigma(\mathcal{W}_{M})<+
\infty \bigr),\qquad  \mathbb{P}_{\theta,x}\mbox{-a.s.}
\end{eqnarray*}
and for any $q\geq0$,
\[
\bar{\mathbb{P}}_{\theta,x}^{\boldsymbol{\gamma}^{\leftarrow q}} \bigl(\sigma(
\mathcal{W}_{M})<+\infty \bigr)=\bar{\mathbb{P}}_{\theta,x}^{\boldsymbol{\gamma}^{\leftarrow q}}
\biggl({\bigcup_{k\geq1}}\bigl\{\sigma(
\mathcal{W}_{M})=k\bigr\} \biggr).
\]
From Theorem~\ref{theokeyproposition} we deduce that for any $q\geq n_{0}$,
\[
{\bigcup_{k\geq1}}\bigl\{\sigma(\mathcal{W}_{M})=k
\bigr\}\subset \Biggl\{ \sup_{k\geq1}\mathbb{I}\bigl\{\sigma(
\mathcal{W}_{M})\geq k\bigr\} \Biggl\llvert \sum
_{j=1}^{k}\gamma_{j+q}
\varsigma_{j}\Biggr\rrvert >\delta_{0} \Biggr\},
\]
which implies that for any $T_{l}\geq l\geq n_{0}$,
\begin{eqnarray*}
&&\mathbb{I}\{T_{l}  <+\infty\}\mathbb{P}_{\theta
_{T_{l}},X_{T_{l}}}
(T_{1}<+\infty )
\\
&&\qquad  \leq\sup_{q\geq l}\sup_{\theta,x\in\mathcal{W}_{M_{0}}\times
C}\bar{
\mathbb{P}}_{\theta,x}^{\boldsymbol{\gamma}^{\leftarrow
q}} \Biggl(\sup_{k\geq1}
\mathbb{I}\bigl\{\sigma(\mathcal{W}_{M})\geq k\bigr\} \Biggl\llvert
\sum_{j=1}^{k}\gamma_{j+q}
\varsigma_{j}\Biggr\rrvert >\delta_{0} \Biggr).
\end{eqnarray*}
Consequently by induction one obtains that for any $n>n_{0}$,
\begin{eqnarray*}
&&\mathbb{P}_{\theta,x}  (T_{n}<+\infty )
\\
&&\qquad  \leq\prod_{l=n_{0}}^{n-1}\sup
_{q\geq l}\sup_{\theta,x\in
\mathcal{W}_{M_{0}}\times C}\bar{\mathbb{P}}_{\theta,x}^{\boldsymbol{\gamma}^{\leftarrow q}}
\Biggl(\sup_{k\geq1}\mathbb{I}\bigl\{\sigma (
\mathcal{W}_{M})\geq k\bigr\}\Biggl\llvert \sum
_{j=1}^{k}\gamma _{j+q}
\varsigma_{j}\Biggr\rrvert >\delta_{0} \Biggr).
\end{eqnarray*}
Result (\ref{eqsupTkinfinite}) then follows by a standard Borel--Cantelli
argument under the condition of the theorem. We now prove that $\{
\theta_{k}\}$
eventually remains in $\mathcal{W}_{M}$, $\mathbb{P}_{\theta,x}$-a.s.
First notice that by construction of $\{T_{i}\}$,
\begin{eqnarray*}
&&{\bigcup_{k\geq1}}  \{ T_{k}=+\infty \}
\\
&&\qquad =  {\bigcup_{k\geq1}} \{ T_{k-1}<+\infty,T_{k}=+\infty \}
\\
&&\qquad =  {\bigcup_{m\geq0}} {\bigcup
_{k\geq1}} \{ T_{k-1}=m,T_{k}=+\infty \}
\\
&&\qquad =  {\bigcup_{m\geq0}} \bigl\{ (\theta_{m},X_{m})
\in \mathcal{W}_{M_{0}}\times C \bigr\} \\
&&\hphantom{\qquad =  {\bigcup_{m\geq0}}}{}\cap \bigl(\{\theta_{l}
\in \mathcal{W}_{M},\forall l\geq m+1\}\cup\{\theta_{l}
\notin\mathcal {W}_{M_{0}},\forall l\geq m+1\} \bigr).
\end{eqnarray*}
Now, since by assumption $\mathbb{P}_{\theta,x} (\bigcap
_{k=1}^{\infty}\bigcup_{i=k}^{\infty}\{(\theta_{i},X_{i})\in\mathcal
{W}_{M_{0}}\times\mathsf{C}\} )=1$,
we deduce that for any $\theta,x\in\Theta\times\mathsf{X}$
\[
\mathbb{P}_{\theta,x} \biggl({\bigcup_{m\geq0}}
\bigl\{(\theta _{m},X_{m})\in\mathcal{W}_{M_{0}}
\times C,\theta_{l}\in\mathcal {W}_{M},\forall l\geq m+1
\bigr\} \biggr)=1,
\]
and we conclude.
\end{pf*}

We briefly discuss here other applications of our stability results,
particularly in the situation where the step-size sequence is held
constant. Such fixed stepsize algorithms have been popular in engineering
since they provide the algorithms with both some form of robustness
and a ``tracking'' ability. The analysis of these algorithms naturally
requires one to establish stability first \cite{benveniste}. We would
like, however, to point out another important application in the context
of adaptive step-size algorithms. Indeed, the choice of $\{\gamma_{i}\}$
is known to have an important impact on the convergence properties
of $\{\theta_{i}\}$. In particular it is well known that if $\{\gamma
_{i}\}$
vanishes too quickly in the early iterations of \eqref
{eqSAwithMarkovianDynamic},
convergence may be very slow. A natural way to address this problem
consists of adaptively selecting the sequence of stepsizes $\{\gamma
_{i}\}$.
A strategy due to Kesten and further extended by Delyon and Juditsky
in \cite{delyon:juditsky:1993} is as follows. Given a nonincreasing
function $\gamma(\cdot):\rightarrow(0,\infty)$ and $s_{0}=0$, consider
 Algorithm \ref{algadaptivestepsizes} which is a modification of~\eqref
{eqSAwithMarkovianDynamic}.
\begin{algorithm}[b]
\caption{Adaptive step-size algorithm}
\label{algadaptivestepsizes}
\begin{itemize}
\item$X_{i+1}\sim P_{\theta_{i}}(X_{i},\cdot)$
\item$\theta_{i+1}=\theta_{i}+\gamma(s_{i})H(\theta_{i},X_{i+1})$
\item$s_{i+1}=s_{i}+\mathbb{I}\{ \langle H(\theta
_{i-1},X_{i}),H(\theta_{i},X_{i+1}) \rangle<0\}$
\end{itemize}
\end{algorithm}

Assume for brevity that the root of $h(\theta)=0$, $\theta_{*}$,
is unique. The rationale behind this recursion is that for sufficiently
regular scenarios, one may expect the event $\{ \langle H(\theta
_{i-1},X_{i}),H(\theta_{i},X_{i+1}) \rangle<0\}$
to occur with higher probability when $\theta_{i}$ is in a neighborhood
$B(\theta_{*},\epsilon)$ of $\theta_{*}$ than when outside this
neighborhood. As a result $\gamma_{i}:=\gamma(s_{i})$ decreases
slowly as long as $\{\theta_{i}\}$ is outside this neighborhood
of $\theta_{*}$, and decreases much faster whenever $\{\theta_{i}\}$
approaches $\theta^{*}$. Convergence to $\theta_{*}$ requires that
$\gamma_{i}\rightarrow0$ or equivalently that $s_{i}\rightarrow
\infty$
with probability one. This means that one should show that for any
$\gamma\in\{\gamma(0),\gamma(1),\gamma(2),\ldots\}$, the fixed stepsize
sequence $\vartheta_{i}^{\gamma}=\vartheta_{i-1}^{\gamma}+\gamma
H(\vartheta_{i-1}^{\gamma},X_{i}^{\gamma})$
for $i\geq1$ is recurrent in the aforementioned neighborhood, which
is the essence of the proof of \cite{delyon:juditsky:1993}. Our results
allow one to establish that the homogeneous Markov chain $\{\vartheta
_{i}^{\gamma},X_{i}^{\gamma}\}$
is recurrent in some set $\mathcal{W}_{M}\times\mathsf{C}$, the  first
crucial step of the proof of \cite{delyon:juditsky:1993}. A detailed 
analysis of such a result is, however, beyond the scope of the present
paper.

\section{Examples: Some adaptive MCMC algorithms}
\label{171717171717}

In this section we illustrate how the results established in Sections~\ref{1313131313}--\ref{secmotivation} can be straightforwardly
applied to a variety of adaptive Markov chain Monte Carlo (MCMC) algorithms,
where the aim is to automatically optimally choose the parameter
$\theta$
of a family of Markov chain transition probabilities $\{P_{\theta
},\theta\in\Theta\}$,
defined on some $\mathsf{X}\subset\mathbb{R}$ sharing a common invariant
distribution with density $\pi(\cdot)$ with respect to the Lebesgue
measure. More specifically, we focus here on the symmetric random
walk Metropolis (SRWM) algorithm with transition probability defined
for $(\theta,x,A)\in\Theta\times\mathsf{X}\times\mathcal
{B}(\mathsf{X})$
as
%
\begin{eqnarray}\label{eqtransitionSRWM}
P_{\theta}(x,A)&=&\int_{A-x}\alpha(x,x+z)q_{\theta}(z)
\,\mathrm {d}z\nonumber\\[-8pt]\\[-8pt]
&&{}+\mathbb{I}\{x\in A\}\int_{\mathsf{X}-x} \bigl(1-\alpha
(x,x+z) \bigr)q_{\theta}(z)\, \mathrm{d}z, \nonumber
\end{eqnarray}
where for any $x,y\in\mathsf{X}^{2}$, $\alpha(x,y):=1\wedge\pi
(y)/\pi(x)$
and $\{q_{\theta}(\cdot),\theta\in\Theta\}$ is a family of symmetric
increment probability densities with respect to the Lebesgue measure
defined on $\mathsf{Z}\times\mathcal{B}(\mathsf{Z})$ for some
$\mathsf{Z}\subset\mathsf{X}$.
Various choices for $q_{\theta}(\cdot)$ are possible.

\begin{algorithm}[b]
\caption{AM algorithm, iteration $i+1$}
\label{algAM}
\begin{itemize}
\item Sample $X_{i+1}\sim P_{\theta_{i}}(X_{i},\cdot)$
\item Update of the tuning parameter
%
\begin{eqnarray}\label
{eqalgorithmhaariosaksman}
\mu_{i+1} & =&\mu_{i}+\gamma_{i+1}(X_{i+1}-
\mu_{i}),
\nonumber
\\[-8pt]\\[-8pt]
\Gamma_{i+1} & =&\Gamma_{i}+\gamma_{i+1}
\bigl((X_{i+1}-\mu _{i}) (X_{i+1}-
\mu_{i})^{\mathsf{T}}-\Gamma_{i}\bigr). \nonumber
\end{eqnarray}
\end{itemize}
\end{algorithm}

The AM (adaptive Metropolis) algorithm \cite{haario:saksman:tamminen:1999,haario:saksman:tamminen:2001}
is concerned with the situation where $\mathsf{X}=\mathbb{R}^{n_{x}}$
for some $n_{x}\geq1$ and $\theta=[\mu|\Gamma]\in\Theta=\mathbb
{R}^{n_{x}}\times\mathcal{C}$
where $\mathcal{C}\subset\mathbb{R}^{n_{x}\times n_{x}}$ is the cone
of symmetric positive definite matrices and $q_{\theta}(z):=\det
^{-1/2}((2.38^{2}/n_{x})(\Gamma+\epsilon_{\mathrm{AM}}I_{n_{x}\times
n_{x}}))\times q (((2.38^{2}/n_{x})(\Gamma+\epsilon
_{\mathrm{AM}}I_{n_{x}\times n_{x}}))^{-1/2}z) $
for $q(z)=\mathcal{N}(z;0,I_{n_{x}\times n_{x}})$ and some $\epsilon
_{\mathrm{AM}}\in(0,1)$.
In fact other choices for $q(\cdot)$ are possible as long as it is
symmetric, that is, $q(z)=q(-z)$ for all $z\in\mathsf{Z}$. In \cite
{gelman:roberts:gilks:1995}
it is shown that in some circumstances the ``optimal''
covariance matrix for the Normal-SRWM is $\Gamma_{\pi}$, where
$\Gamma_{\pi}$
is the true covariance matrix of the target distribution $\pi(\cdot)$,
assumed here to exist. The AM algorithm of \cite{haario:saksman:tamminen:2001}
essentially implements the following algorithm to estimate $\Gamma$
on the fly. Let $\epsilon_{\mathrm{AM}}>0$ and let $X_{0}=x\in\mathsf{X}$,
then for $i\geq0$ and with $\theta_{i}:=[\mu_{i}|\Gamma_{i}]$ here one can consider Algorithm~\ref{algAM}.

It was realized in \cite{andrieu:robert:2001} that this algorithm
is a particular instance of (\ref{eqSAwithMarkovianDynamic}) where
$H\dvtx \Theta\times\mathsf{X}\rightarrow\Theta$ is
%
\begin{equation}\label{eqdefinitionchamphaariosaksman}
H(\theta,x):= \bigl[x-\mu|(x-\mu) (x-\mu)^{\mathsf{T}}-\Gamma
\bigr]{}^{\mathsf{T}},
\end{equation}
and the corresponding mean field is
\[
h(\theta)= \bigl[\mu_{\pi}-\mu|(\mu_{\pi}-\mu) (
\mu_{\pi}-\mu )^{\mathsf{T}}+\Gamma_{\pi}-\Gamma
\bigr]{}^{\mathsf{T}}.
\]
We show in Section~\ref{subStability-of-AM} (Theorem~\ref
{thmCVAMalgorithms})
that the stability of these recursions is a direct consequence of
the result of Section~\ref{1313131313} and a result from \cite
{saksman:vihola:2008},
which establishes \ref{hypfixedthetaconditiononPtheta} for a
class of target distribution densities $\pi(\cdot)$. The result of
Section~\ref{secmotivation} then directly applies to the AM algorithm,
leading to the conclusion that $\{\theta_{i}\}$ eventually remains
in a compact set with probability one. While the boundedness of $\{
\theta_{i}\}$
has already been established in \cite{saksman:vihola:2008} using
different arguments, our results are more general in several ways.
For example, Theorem~\ref{thmCVAMalgorithms} shows that the AM
algorithm is stable when the sequence of stepsizes $\{\gamma_{i}\}$
is constant, which opens up the way for the analysis of more sophisticated
and robust versions of the AM algorithm. Theorem~\ref{thmCVAMalgorithms}
also shows that the AM algorithm is also stable for heavier tailed
distributions than in \cite{saksman:vihola:2008}, in the situation
where $\mathsf{X}=\mathbb{R}$, for both decreasing or constant stepsize
sequences. As should be clear from our current analysis, a full study
of the multivariate scenario is a different (and significant) research
project.

\begin{algorithm}[b]
\caption{Coerced acceptance probability RWM, iteration $i+1$}
\label{algstandardcoerced}
\begin{itemize}
\item Update the state $X_{i},Y_{i}$, with $Z_{i+1}\sim q_{\theta
_{i}}(\cdot)$
\begin{eqnarray*}
Y_{i+1} & = & X_{i}+Z_{i+1},
\\
X_{i+1} & = & %
\cases{ Y_{i+1} & \quad \mbox{with
probability \ensuremath{\alpha (X_{i},Y_{i+1})}},
\cr
X_{i} & \quad \mbox{otherwise} } %
\end{eqnarray*}

\item Update the scaling parameter
\[
\theta_{i+1}=\theta_{i}+\gamma_{i+1}\bigl(
\alpha(X_{i},Y_{i+1})-\alpha _{\ast}\bigr).
\]
\end{itemize}
\end{algorithm}

We now turn to another type of popular adaptive scheme for the SRWM.
Let $\mathsf{X}=\mathbb{R}^{n_{x}}$ and $\Theta=\mathbb{R}$. Suppose
$q(\cdot)$ is a symmetric probability density on $\mathsf{X}$, and
define the family of proposal distributions $\{q_{\theta}(\cdot
),\theta\in\Theta\}$
as $q_{\theta}(z):=\exp(-\theta)q(\exp(-\theta)z)$. A possible increment
probability density is again $q_{\theta}(z)=\mathcal{N}(z,\sigma:=\exp(\theta))$.
Let $\alpha_{*}\in(0,1)$ be a desired mean acceptance probability
for the SRWM. The following algorithm aims to optimize $\theta^{*}$
in order to achieve an expected acceptance rate of $\alpha_{*}$ \cite
{andrieu:robert:2001}
and is often used as one of the components of more sophisticated schemes.
The resulting procedure is displayed in Algorithm \ref{algstandardcoerced}.

In Section~\ref{subdriftforSRWM} we prove the stability of $\{
\theta_{i},X_{i}\}$
for a broad class of probability densities $\pi(\cdot)$, including
a heavy-tailed scenario and situations where the stepsize sequence
is constant (Theorem~\ref{thmCVcoercedalgorithms}). It should be
pointed out that in this case, in contrast with the AM algorithm scenario,
we do not require a lower bound on the scaling factor $\exp(\theta)$,
which requires establishing \ref{hypfixedthetaconditiononPtheta}
for both arbitrarily large and small values of $\exp(\theta)$ and
leads us to proving the new result Theorem~\ref
{thmdriftconditionsimpletailcondition}
(a stability result has been proved in \cite{vihola:2011}, but in
a less general scenario). In fact the theory we have developed suggests
improvements on this standard algorithm whose stability can be easily
established thanks to the theory developed earlier in the paper. An
example is given below: the rationale behind the algorithm is that
for very poor initializations, the increments on the parameter are
initially large, while still leading to a stable dynamic.
See Algorithm \ref{algfastcoerced}.
\begin{algorithm}[t]
\caption{Fast coerced acceptance probability RWM, iteration $i+1$}
\label{algfastcoerced}
\begin{itemize}
\item Update the state $X_{i},Y_{i}$, with $Z_{i+1}\sim q_{\theta
_{i}}(\cdot)$
\begin{eqnarray*}
Y_{i+1} & = & X_{i}+Z_{i+1},
\\
X_{i+1} & = & %
\cases{ Y_{i+1} & \quad \mbox{with
probability \ensuremath{\alpha (X_{i},Y_{i+1})}},
\cr
X_{i} &\quad  \mbox{otherwise} } %
\end{eqnarray*}

\item Update the scaling parameter
\[
\theta_{i+1}=\theta_{i}+\gamma_{i+1}\bigl(|
\theta_{i}|+1\bigr) \bigl(\alpha (X_{i},Y_{i+1})-
\alpha_{\ast}\bigr).
\]
\end{itemize}
\end{algorithm}

The proofs of stability of the three algorithms above rely on common
key intermediate results. In Section~\ref{subdriftforSRWM} we
establish \ref{hypfixedthetaconditiononPtheta} for the SRWM
under two different sets of assumptions on $\pi(\cdot)$ and $q(\cdot)$.
In Section~\ref{subStability-of-AM} we establish \ref{hypwdriftSA1}
for the AM algorithm and conclude with Theorem~\ref{thmCVAMalgorithms},
while in Section~\ref{subStability-of-coerced} we establish \ref
{hypwdriftSA1}
for the coerced acceptance algorithms, and conclude with Theorem~\ref
{thmCVcoercedalgorithms}.

\subsection{Establishing \texorpdfstring{\protect\ref{hypfixedthetaconditiononPtheta}}{(A2)} for SRWM algorithms}
\label{subdriftforSRWM}

\subsubsection{The superexponential \texorpdfstring{``$\Gamma+\epsilon_{\mathrm{AM}}I_{n_{x}\times n_{x}}$''}{``$Gamma+epsilon_{AM} I_{n_{x} times n_{x}}$''} scenario}

In this section we state a result which establishes that \ref
{hypfixedthetaconditiononPtheta}
is satisfied for the SRWM transition probability on $\mathsf
{X}=\mathbb{R}^{n_{x}}$
under suitable conditions on $\pi(\cdot)$ and $q(\cdot)$ \cite{jarner:hansen:1998}.
\renewcommand\thelonglist{(A6)}
\begin{longlist}[(A6)]
\item \label{hyppropertiesPI}The
probability distribution
$\pi(\cdot)$ has the following properties:
\begin{longlist}[(3)]
\item[(1)] it is positive on every compact set and continuously differentiable;
\item[(2)] there exists $\rho>1$ such that
%
\begin{equation}
\lim_{R\rightarrow+\infty}\sup_{\{x:|x|\geq R\}}\frac{x}{\llvert  x\rrvert ^{\rho}}
\cdot\nabla\log\pi (x )=-\infty;
\end{equation}

\item[(3)] the contours $\partial A (x )= \{ y\dvtx \pi(y)=\pi
(x) \} $
are asymptotically regular, that is, for some $R>0$,
%
\begin{equation}
\sup_{\{x:|x|\geq R\}}\frac{x}{\llvert  x\rrvert }\cdot\frac
{\nabla\pi (x )}{\llvert \nabla\pi (x
)\rrvert }<0;
\end{equation}

\item[(4)] the proposal distribution density $q$ is that of a standardized Gaussian
or Student's $t$-distribution.
\end{longlist}
\end{longlist}

The following theorem quantifies the way in which ergodicity of the
SRWM vanishes under \ref{hyppropertiesPI} as some of its eigenvalues
become large. The norm used for matrices below is $|A|=\sqrt{\mathsf
{Tr}(AA^{\mathsf{T}})}$
and recall that here $\theta=[\mu|\Gamma]\in\Theta=\mathbb
{R}^{n_{x}}\times\mathcal{C}$.

\begin{prop}
\label{propsimultaneousdriftsuperexponential}Let $\eta\in(0,1)$
and $V(x)\propto\pi^{-\eta}(x)$. Under \textup{\ref{hyppropertiesPI}}
one can choose $V\geq1$ and there exist $a,b\in(0,\infty)$ and
$\mathsf{C}=B(0,R)$
for some $R>0$ such that for any $\theta,x\in\Theta\times\mathsf{X}$,
\[
P_{\theta}V(x)\leq\bigl(1-a/\sqrt{\det(\Gamma+\epsilon
_{\mathrm{AM}}I_{n_{x}\times n_{x}})}\bigr)V(x)+b\mathbb{I}\{x\in\mathsf{C}\},
\]
and we note that for any $\Gamma\in\mathcal{C}$,
\[
\sqrt{\det(\Gamma)}\leq n_{x}^{-n_{x}/4}|\Gamma|^{n_{x}/2}.
\]
\end{prop}

\begin{pf}
The first statement is proved in \cite{saksman:vihola:2008}, Proposition~15,
and the second statement follows from the standard arithmetic/geometric
mean inequality applied to the eigenvalues of $\Gamma^{2}$,
%
\begin{equation}\label{eqsmallresultdettheta}
\det\bigl(\Gamma^{2}\bigr)\leq \biggl(\frac{\mathsf{Tr} (\Gamma
^{2} )}{n_{x}}
\biggr)^{n_{x}}= \biggl(\frac{|\Gamma
|^{2}}{n_{x}} \biggr)^{n_{x}}.
\end{equation}
\upqed
\end{pf}

\subsubsection{Establishing \texorpdfstring{\protect\ref{hypfixedthetaconditiononPtheta}}{(A2)}
for the AM
algorithms with weak tail assumptions}

In this section we prove \ref{hypfixedthetaconditiononPtheta}
for the SRWM algorithm on $\mathsf{X}=\mathbb{R}$ in the situation
where no lower bound on the scaling parameter of the proposal distribution
is assumed and under a weaker assumption on the vanishing rate of
the tails of the target density $\pi(\cdot)$ than in the previous
subsection. More precisely, let $P_{\theta}$ denote here the random-walk
Metropolis kernel with symmetric proposal distribution $q_{\theta
}(z)=\exp(-\theta)q(z/\exp(\theta))$
for $\theta\in\Theta:=(-\infty,\infty)$. For notational simplicity,
in this subsection, we introduce the piece of notation $\sigma=\exp
(\theta)$
and use $P_{\sigma}$ instead of $P_{\log\sigma}$, $q_{\sigma}$
instead of $q_{\log\sigma}$ throughout and say that $\sigma\in\exp
(\Theta)$.
We will also use the piece of notation $\ell(x):=\log\pi(x)$. We
require the following assumptions on $\pi(\cdot)$ and the increment
proposal density $q_{\sigma}$:
\renewcommand\thelonglist{(A7)}
\begin{longlist}[(A7)]
\item\label{hyppismallthetasubgeomdrift}
\mbox{}
\renewcommand\theenumi{(\arabic{enumi})}
\renewcommand\labelenumi{\theenumi}
\begin{enumerate}[(3)]
\item\label{subhypregularityderivatives}The target distribution $\pi
(\cdot)$
on $(\mathsf{X},\mathcal{B}(\mathsf{X}))$ has the following properties:
\begin{longlist}[(c)]
\item[(a)] It has a density $\pi(x)$ with respect to the Lebesgue measure,
\item[(b)]$\pi(x)$ is bounded away from $0$ on any compact set of $\mathbb{R}$,
\item[(c)]$\ell(x)$ is twice differentiable. We denote $\ell'(x):=\nabla
\ell(x)$
and $\ell''(x):=\nabla^{2}\ell(x)$,
\item[(d)] for any $M>0$, defining $\epsilon_{x}:=M/|\ell'(x)|$,
\begin{eqnarray*}
\lim_{R\rightarrow\infty}\sup_{x\in B^{c}(0,R)}\sup
_{|t|\leq
\epsilon_{x}}\frac{|\ell''(x+t)|}{|\ell'(x+t)|^{2}}&=&0,
\\
\lim_{R\rightarrow\infty}\sup_{x\in B^{c}(0,R)}\sup
_{|t|\leq
\epsilon_{x}}\frac{|\ell''(x+t)|}{|\ell'(x)|^{2}}&=&0,
\\
\lim_{R\rightarrow\infty}\sup_{x\in B^{c}(0,R)}\sup
_{0\leq t\leq
\epsilon_{x}}\biggl|\frac{|\ell'(x-t)|}{|\ell'(x+t)|}-1\biggr|&=&0.
\end{eqnarray*}
\end{longlist}
\item \label{subhypsubgeometrictailsofpi}The tails of $\pi(x)$ decay
at a minimum rate characterized as follows: there exist $p\in(0,1)$
such that
\begin{eqnarray*}
\lim_{R\rightarrow\infty}\sup_{x\in B^{c}(0,R)}\frac{\ell
'(x)}{|x|^{p-1}}&<&0
\quad \mbox{and}\\
  \forall\gamma\in(0,1) \qquad \lim
_{R\rightarrow\infty}\inf_{x\in B^{c}(0,R)}\frac
{\pi^{-\gamma}(x)}{|\ell'(x)|}&>&0.
\end{eqnarray*}

\item\label{subhypincrementproposaldensity}The increment proposal density
$q_{\sigma}(z)$ is of the form $q_{\sigma}(z)=\frac{1}{\sigma
}q(z/\sigma)$
for some symmetric probability density $q(z)$, such that
$\operatorname{supp}(q)=[-1,1]=:\mathsf{Z}$
and $q(\cdot)\dvtx \mathsf{Z}\rightarrow[\underline{q},\bar{q}]$ for
$\underline{q},\bar{q}\in(0,\infty)$.
\end{enumerate}
\end{longlist}

%
\begin{rem}
Consider $\ell(x)=C-|x|^{\alpha}$, for $|x|\geq R_{\ell}$ and
$\alpha>0$.
Then for $x\geq R_{\ell}$, $\ell'(x)=-\alpha x^{\alpha-1}$ and $\ell
''(x)=-\alpha(\alpha-1)x^{\alpha-2}$
and all the conditions in \ref{hyppismallthetasubgeomdrift}
are satisfied.
\end{rem}

\begin{rem}
The support condition on $q$ can be removed, but this requires one
to control additional ``tail integral'' terms in the proofs of this
section, which would add further to already long arguments. We have
opted for this presentation for brevity and clarity since it is the
terms that we handle which are both crucial and difficult to control.
\end{rem}

The following theorem establishes \ref{hypfixedthetaconditiononPtheta}
for the scalar SRWM with $V(x)\propto\pi^{-\beta}(x)$ for some
$\beta\in(0,1)$
under \ref{hyppismallthetasubgeomdrift}.

%
\begin{thm}
\label{thmdriftconditionsimpletailcondition}Consider the SRWM targetting
$\pi(\cdot)$ and with increment proposal density $q_{\sigma}$. Assume
they satisfy \textup{\ref{hyppismallthetasubgeomdrift}}, and let
$V(x):=c\pi^{-\eta}(x)$
for some $\eta\in(0,1)$ and $c\in(0,\infty)$ such that $V\geq1$.
Then for any $\iota\in(0,1)$ there exists $R\geq0$ and $a_{0}\in
(0,\infty)$
such that:
\begin{longlist}[(2)]
\item[(1)] for any $x\in B^{c}(0,R)$, with $\tilde{a}^{-1}(\sigma
):=a_{0}/(\sigma\vee\sigma^{-2})$,
\[
P_{\sigma}V(x)\leq V(x)-\tilde{a}^{-1}(\sigma)V^{\iota}(x),
\]

\item[(2)] there exists a constant $b\in(0,\infty)$ such that for all
$\sigma,x\in\exp(\Theta)\times B(0,R)$
we have $P_{\sigma}V(x)\leq b$.
\end{longlist}
\end{thm}

\begin{pf}
Let $\iota\in(0,1)$ and $R\geq R_{0}$ such that $\inf_{x\in
B^{c}(0,R)}V^{1-\iota}(x)|\ell'(x)|^{2}>0$
and $\inf_{x\in B^{c}(0,R)}V^{1-\iota}(x)/|\ell'(x)|>0$, where $R_{0}$
is given in Proposition~\ref{propintermediatePVoverVminus1upperbound}.
The existence of $R$ is ensured by \ref
{hyppismallthetasubgeomdrift}\ref{subhypsubgeometrictailsofpi}
and the choice of $V$. Indeed, from the assumption, for $x\in B^{c}(0,R)$,
we have from Lemma~\ref{lempropertiesUpsilonetc} that $V(x)\geq
C_{\Upsilon,2}\exp(\eta C_{\Upsilon,1}^{-1}|x|^{p})$
for some constant $C>0$ and $|\ell'(x)|\geq C_{\ell}|x|^{p-1}$ for
$x\geq R_{\ell}$ for some $C_{\ell},R_{\ell}>0$, and we can conclude.
From Proposition~\ref{propintermediatePVoverVminus1upperbound} below,
for $x\in B^{c}(0,R)$ and $\sigma\leq c_{0}/|\ell'(x)|$, we have
\begin{eqnarray*}
P_{\sigma}V(x) & \leq& V(x)-a'_{0}
\sigma^{2}\bigl|\ell'(x)\bigr|^{2}V(x)
\\
& \leq& V(x)-a'_{0}\inf_{x_{0}\in B^{c}(0,R)}\bigl|\ell
'(x_{0})\bigr|^{2}V^{1-\iota}(x_{0})
\times\sigma^{2}V^{\iota}(x).
\end{eqnarray*}
For $|x|\geq\sigma\geq c_{0}/|\ell'(x)|$, we have
\begin{eqnarray*}
P_{\sigma}V(x) & \leq& V(x)-a'_{0}V(x)/\bigl|\sigma
\ell'(x)\bigr|
\\
& \leq& V(x)-a'_{0}\inf_{x_{0}\in B^{c}(0,R)}V^{1-\iota}(x_{0})/\bigl|
\ell '(x_{0})\bigr|\sigma^{-1}V^{\iota}(x).
\end{eqnarray*}
For $\sigma\geq|x|\geq R$, we have
\begin{eqnarray*}
P_{\sigma}V(x) & \leq& V(x)-a'_{0}|x|\times V(x)/
\sigma
\\
& \leq& V(x)-a'_{0}R\times V(x)/\sigma.
\end{eqnarray*}
Now we can use the trivial inequalities $\sigma^{2}\leq\sigma\leq
1\leq\sigma^{-1}$
(case $\sigma\in(0,1]$) and $\sigma^{-1}\leq1\leq\sigma\leq
\sigma^{2}$
[case $\sigma\in(1,\infty)$] which lead to the following upper bound:
\[
P_{\sigma}V(x)\leq V(x)-a_{0}\bigl(\sigma^{-1}
\wedge\sigma^{2}\bigr)V^{\iota}(x).
\]
The second claim follows immediately from the bound $PV(x)\leq2V(x)$
easily obtained from (\ref{eqPVoverVminus1decompositin}) and the
fact that $\sup_{x\in B(0,R)}V(x)<\infty$.
\end{pf}

%
\begin{prop}
\label{propintermediatePVoverVminus1upperbound}Consider the SRWM
targetting $\pi(\cdot)$ satisfying \ref{hyppismallthetasubgeomdrift}.
Let $V(x):=c\pi^{-\eta}(x)$ for some $\eta\in(0,1)$ and $c$ such
that for all $x\in\mathsf{X}$, $V(x)\geq1$. Then there exist $a'_{0},R_{0}>0$
such that for any $x\in B^{c}(0,R_{0})$ and any $\sigma\in\exp
(\Theta)=(0,\infty)$,
\[
\frac{P_{\sigma}V(x)}{V(x)}-1\leq-a'_{0}\times %
\cases{ \sigma^{2}\bigl|\ell'(x)\bigr|^{2}, &\quad  \mbox{if }$
\sigma\bigl|\ell'(x)\bigr|<c_{0}$,
\cr
1/\bigl|\sigma\ell'(x)\bigr|,
& \quad \mbox{if }$|x|\geq\sigma\geq c_{0}/\bigl|\ell '(x)\bigr|$,
\cr
|x|/\sigma, & \quad \mbox{if }$1>|x|/\sigma$.} %
\]
\end{prop}

\begin{pf}
Without loss of generality we detail the situation where $x>0$ since
the case $x<0$ can be straightforwardly addressed by considering
the density $\pi_{-}(x):=\pi(-x)$ which also satisfies \ref
{hyppismallthetasubgeomdrift},
and hence Lemmata \ref{lemdecompPVoverVsubgeometric}, \ref
{lemsubgeomboundpsithetalessthanx}
and \ref{lemsubgeomboundonTsthetalargerthanx} (given below). In
what follows the terms $T_{i}(\sigma,x)$ for $i=1,2,3,4$ are defined
in Lemma~\ref{lemdecompPVoverVsubgeometric}. Choose $R\geq
R_{PV}\vee R_{\psi}\vee R_{T}$
such that for $x\in B^{c}(0,R)$, $x\geq c_{0}/|\ell'(x)|$, where
$R_{PV},R_{\psi},R_{T}$ and $c_{0}$ are as in Lemmata \ref
{lemdecompPVoverVsubgeometric},
\ref{lemsubgeomboundpsithetalessthanx} and \ref
{lemsubgeomboundonTsthetalargerthanx}.
First from Lemma~\ref{lemsubgeomboundpsithetalessthanx}, we have
for $x\geq c_{0}/|\ell'(x)|$ and any $\sigma\in\exp(\Theta)$,
\begin{eqnarray*}
T_{1}(  \sigma,x)&=&\int_{0}^{\sigma\wedge c_{0}/|\ell'(x)|}\psi
_{x}(z)q_{\sigma}(z)\,\mathrm{d}z\\
&&{}+\mathbb{I}\bigl\{\sigma\geq
c_{0}/\bigl|\ell '(x)\bigr|\bigr\}\int_{c_{0}/|\ell'(x)|}^{\sigma\wedge x}
\psi_{x}(z)q_{\sigma
}(z)\,\mathrm{d}z
\\
&\leq& -\epsilon_{\psi}\bigl|\ell'(x)\bigr|^{2}\int
_{0}^{\sigma\wedge
c_{0}/|\ell'(x)|}z^{2}q_{\sigma}(z)
\,\mathrm{d}z\\
&&{}-\epsilon_{\psi
}\mathbb{I}\bigl\{\sigma\geq
c_{0}/\bigl|\ell'(x)\bigr|\bigr\}\int_{c_{0}/|\ell
'(x)|}^{\sigma\wedge x}q_{\sigma}(z)
\,\mathrm{d}z
\\
&\leq& -\epsilon_{\psi}\underline{q}/3\bigl|\ell'(x)\bigr|^{2}
\bigl[\sigma\wedge c_{0}/\bigl|\ell'(x)\bigr|\bigr]^{3}/
\sigma\\
&&{}-\epsilon_{\psi}\underline{q}\mathbb {I}\bigl\{\sigma\geq
c_{0}/\bigl|\ell'(x)\bigr|\bigr\}\bigl[\sigma\wedge
x-c_{0}/\bigl|\ell '(x)\bigr|\bigr]/\sigma,
\end{eqnarray*}
and therefore for $\sigma\leq x$
\begin{eqnarray*}
T_{1}(\sigma,x) & \leq&-\epsilon_{\psi}\underline{q} \biggl[
\mathbb {I} \biggl\{\sigma<\frac{c_{0}}{|\ell'(x)|} \biggr\}\frac{\sigma
^{2}|\ell'(x)|^{2}}{3}\\
&&\hphantom{-\epsilon_{\psi}\underline{q} \biggl[}{}+
\mathbb{I} \biggl\{\sigma\geq\frac
{c_{0}}{|\ell'(x)|} \biggr\} \biggl(
\frac{c_{0}^{3}}{3\sigma|\ell
'(x)|}+1-\frac{c_{0}}{\sigma|\ell'(x)|} \biggr) \biggr]
\\
& \leq&\frac{-\epsilon_{\psi}\underline{q}}{3} \biggl[\sigma ^{2}\bigl|\ell'(x)\bigr|^{2}
\times\mathbb{I} \biggl\{\sigma<\frac
{c_{0}}{|\ell'(x)|} \biggr\}\\
&&\hphantom{\frac{-\epsilon_{\psi}\underline{q}}{3} \biggl[}{}+\frac{c_{0}^{3}}{\sigma|\ell
'(x)|}
\times\mathbb{I} \biggl\{ x\geq\sigma\geq\frac{c_{0}}{|\ell
'(x)|} \biggr\} \biggr].
\end{eqnarray*}
Now from Lemma~\ref{lemsubgeomboundonTsthetalargerthanx} for
$\sigma\geq x\geq R$
we have
\begin{eqnarray*}
T_{1}(\sigma,x)+T_{2}(\sigma,x)&\leq&-\epsilon_{T}
\times x/\sigma,
\\
T_{3}(\sigma,x)&\leq&0,
\end{eqnarray*}
and for $\sigma\geq x-\Upsilon(x)$ we have
\[
T_{3}(\sigma,x)+T_{4}(\sigma,x)\leq-\epsilon_{T}
\times\bigl(-\Upsilon (x)\bigr)/\sigma
\]
and we conclude with Lemma~\ref{lemdecompPVoverVsubgeometric} and
by treating the case where $x<0$ in a similar fashion.
\end{pf}

Note that, as pointed out in the proof of Proposition~\ref
{propintermediatePVoverVminus1upperbound},
it is sufficient to specialize most of the results of Lemmata \ref
{lempropertiesUpsilonetc},
\ref{lemdecompPVoverVsubgeometric}, \ref
{lemsubgeomboundpsithetalessthanx}
and \ref{lemsubgeomboundonTsthetalargerthanx} (stated below and
proved in Appendix~\ref{secProofs-from-Section}) to the case $x>0$.
The following lemma establishes some key properties implied specifically
by \ref{hyppismallthetasubgeomdrift}\ref
{subhypsubgeometrictailsofpi},
which will also be used in Section~\ref{subStability-of-coerced}.

%
\begin{lem}
\label{lempropertiesUpsilonetc}Assume that $\pi(\cdot)>0$, is differentiable
and satisfies \textup{\ref{hyppismallthetasubgeomdrift}\ref
{subhypsubgeometrictailsofpi}},
define for any $\gamma>0$,
\begin{eqnarray*}
I_{\gamma}(x) &:=&\int_{0}^{\infty} \biggl(
\frac{\pi(x+\sgn
(x)z)}{\pi(x)} \biggr)^{\gamma}\,\mathrm{d}z \quad \mbox{and}\\
J_{\gamma
}(x)&:=&\int_{0}^{|x|} \biggl(
\frac{\pi(x)}{\pi(x-\sgn(x)z)} \biggr)^{\gamma}\,\mathrm{d}z,
\end{eqnarray*}
with $\sgn(x):=x/|x|$ for $x\neq0$ and for any $x>0$, $\Upsilon
(x):=\inf\{y\in\mathsf{X}\dvtx \pi(y)=\pi(x)\}$.
Then:
\begin{enumerate}[(2)]
\item[(1)] the function $\Upsilon(\cdot)$ has the following properties:

\begin{longlist}[(b)]
\item[(a)]$\lim_{x\rightarrow\infty}\Upsilon(x)=-\infty$,
\item[(b)]\emph{there exist constants $C_{\Upsilon,1},C_{\Upsilon,2}\in(0,\infty)$
such that for all $|x|\geq R_{\Upsilon}$
\[
\bigl|\Upsilon(x)\bigr|\vee|x|\leq C_{\Upsilon,1} \bigl(-\log\bigl(\pi
(x)/C_{\Upsilon,2}\bigr) \bigr)^{1/p}
\]
{[}or equivalently $\pi^{-1}(x)\geq C_{\Upsilon,2}\exp
(C_{\Upsilon,1}^{-1}[|\Upsilon(x)|^{p}\vee|x|^{p}] )${]}},
\end{longlist}
\item[(2)] and there exists a constant $C_{\gamma}\in(0,\infty)$ such
that for
any $x\in\mathsf{X}$, $I_{\gamma}(x)\vee J_{\gamma}(x)\leq
C_{\gamma}|x|^{1-p}$.
\end{enumerate}
\end{lem}

We now find a convenient expression for $P_{\sigma}V(x)/V(x)$ valid
for sufficiently large $x$ and all $\sigma$'s.

%
\begin{lem}
\label{lemdecompPVoverVsubgeometric}Assume \ref
{hyppismallthetasubgeomdrift},
and for $x,\eta,s,z\in\mathsf{X}\times(0,1)\times\{-1,1\}\times
\mathsf{Z}$
define $\phi_{x,\eta,s}(z):=[\pi(x+sz)/\pi(x)]^{\eta}$ and
\[
\psi_{x}(z):= \bigl(\phi_{x,-\eta,-1}(z)-1 \bigr)+ \bigl(\phi
_{x,1-\eta,1}(z)-1 \bigr)-\bigl(\phi_{x,1,1}(z)-1 \bigr).
\]
For any $x\geq0$, define $\Upsilon(x):=\inf\{y\in\mathsf{X}:\pi
(x)=\pi(y)\}$,
and let $V(x)\propto\pi^{-\eta}(x)$. Then there exists $R_{PV}>0$
such that for all $x\geq R_{PV}$ and any $\sigma\in\exp(\Theta)$
\begin{eqnarray*}
\frac{P_{\sigma}V(x)}{V(x)}-1 & = & \sum_{i=1}^{4}T_{i}(
\sigma,x)
\end{eqnarray*}
with
\begin{eqnarray*}
T_{1}(\sigma,x) & =&\int_{0}^{\sigma\wedge x}
\psi_{x}(z)q_{\sigma
}(z)\,\mathrm{d}z,
\\
T_{2}(\sigma,x) & =&\mathbb{I}\{\sigma\geq x\}\int_{x}^{\sigma
}
\bigl[\phi_{x,1-\eta,1}(z)-\phi_{x,1,1}(z) \bigr]q_{\sigma
}(z)
\,\mathrm{d}z,
\\
T_{3}(\sigma,x) & =&\mathbb{I}\{\sigma\geq x\}\int
_{\Upsilon
(x)}^{(\sigma-x+\Upsilon(x))\wedge0} \bigl[\phi_{\Upsilon(x),-\eta,-1}(z)-1
\bigr]q_{\sigma}\bigl(z+x-\Upsilon(x)\bigr)\,\mathrm{d}z,
\\
T_{4}(\sigma,x) & =&\mathbb{I}\bigl\{\sigma\geq x-\Upsilon(x)\bigr\}
\\
&&{} \times\int_{0}^{\sigma-(x-\Upsilon(x))} \bigl[
\phi_{\Upsilon
(x),1-\eta,-1}(z)-1+1-\phi_{\Upsilon(x),1,-1}(z) \bigr]\\
&&\hphantom{{} \times\int_{0}^{\sigma-(x-\Upsilon(x))}}{}\times q_{\sigma
}
\bigl(z+x-\Upsilon(x)\bigr)\,\mathrm{d}z.
\end{eqnarray*}
\end{lem}

Here we prove some properties of $\psi_{x}(z)$ which will
allow us to upper bound the term $T_{1}(\sigma,x)$ in the case where
$\sigma\leq x$.

%
\begin{lem}
\label{lemsubgeomboundpsithetalessthanx}Assume \textup{\ref
{hyppismallthetasubgeomdrift}\ref{subhypregularityderivatives}}
and for $\eta\in(0,1)$ let $\psi_{x}(z)$ be as in Lemma~\ref{lemdecompPVoverVsubgeometric}. Then there exist constants
$c_{0},\epsilon_{\psi},R_{\psi}>0$ such that for all $x\geq R_{\psi}$,
$\psi_{x}(z)\leq0$ for $z\in[0,x]$ and $\psi_{x}(z)$ satisfies
the following upper bounds:
%
\begin{equation}\label{eqboundonpsi}
\psi_{x}(z)\leq-\epsilon_{\psi}\times %
\cases{ \bigl|
\ell'(x)\bigr|^{2}z^{2}, &\quad  \mbox{for } $0\leq z\leq
c_{0}/\bigl|\ell'(x)\bigr|$,
\cr
1, & \quad \mbox{for }
$c_{0}/\bigl|\ell'(x)\bigr|\leq z\leq x$. } %
\end{equation}
\end{lem}

Now in the following lemma we address the situation where
$\sigma\geq x$ and require an additional assumption on the vanishing
speed of $\ell'(x)$.

%
\begin{lem}
\label{lemsubgeomboundonTsthetalargerthanx}Assume \textup{\ref
{hyppismallthetasubgeomdrift}},
and let $T_{i}(\sigma,x)$ for $i=1,\ldots,4$ be as defined in Lemma~\ref{lemdecompPVoverVsubgeometric}. Then there exist
$C_{T},R_{T},\epsilon_{T}>0$
such that for $x\geq R_{T}$ and:
\begin{longlist}[(3)]
\item[(1)] for $\sigma\geq x$
\[
T_{1}(\sigma,x)+T_{2}(\sigma,x)\leq-\epsilon_{T}
\times x/\sigma,
\]

\item[(2)] for $\sigma\geq x$
\[
T_{3}(\sigma,x)\leq0,
\]

\item[(3)] and for $\sigma\geq x-\Upsilon(x)$
\begin{eqnarray*}
T_{3}(\sigma,x)  +T_{4}(\sigma,x)\leq-
\epsilon_{T}\bigl(-\Upsilon (x)\bigr)/\sigma.
\end{eqnarray*}
\end{longlist}
\end{lem}

\subsection{Stability of the AM algorithms}
\label{subStability-of-AM}

Thanks to Theorem~\ref{propestablishcompositedrift} and its corollary
we know that recurrence is ensured as soon as \ref
{hypfixedthetaconditiononPtheta}
and \ref{hypwdriftSA1} are satisfied. In the previous section
we established conditions on $\pi(\cdot)$ and $q_{\theta}(\cdot)$
under which \ref{hypfixedthetaconditiononPtheta} is satisfied
for the transition probabilities underpinning Algorithm \ref{algAM}.
We therefore focus on checking that \ref{hypwdriftSA1} is satisfied.
First we start with a result which, together with Theorem~\ref
{thlocalergodicityimpliesstability},
leads to the same conclusions as \cite{saksman:vihola:2008} when
$\{\gamma_{i}\}$ is not constant, but also to the additional stability
of the time-homogeneous Markov chain $\{\theta_{i},X_{i}\}$ when
$\gamma_{i}=\gamma_{0}$ for any $i\geq0$.

%
\begin{thm}
\label{thmCVAMalgorithms}Consider the controlled MC defined by
Algorithm \ref{algAM} for $\mathsf{X}=\mathbb{R}^{n_{x}}$ with
$n_{x}\geq1$
(resp., Algorithm \ref{algAM} for $\mathsf{X}=\mathbb{R}$), assume that
$\pi(\cdot)$ and $q(\cdot)$ satisfy \textup{\ref{hyppropertiesPI}} [resp.,
\textup{\ref{hyppismallthetasubgeomdrift}}] and that $\{\gamma_{i}\}$
is such that $\limsup_{i\rightarrow\infty}\gamma_{i+1}^{-1}-\gamma
_{i}^{-1}<1$.
Then for any $\epsilon>0$ there exist $M,R>0$ such that with $\mathcal
{W}_{M}:=\{\theta\in\Theta\dvtx w(\theta)\leq M\}$
for $w(\theta)=1+|\mu|^{2+\epsilon}+|\Gamma|$
\[
\mathbb{P}_{\theta,x} \Biggl(\bigcap_{k=0}^{\infty}
\bigcup_{i\geq k}\bigl\{ (\theta_{i},X_{i})\in
\mathcal{W}_{M}\times B(0,R)\bigr\} \Biggr)=1.
\]
\end{thm}

The proofs of the theorem for the two sets of assumptions rely on
the following proposition, which establishes \ref{hypwdriftSA1}
for a suitable Lyapunov function $w(\cdot)$. Note that despite its
dependence on \ref{hypfixedthetaconditiononPtheta} the result
does not depend on the expression for $a(\cdot)$.

%
\begin{prop}
\label{propwforAM}Let $\epsilon>0$, and define $w\dvtx \Theta
\rightarrow[0,\infty)$
as
\[
w(\theta):=1+|\mu|^{2+\epsilon}+|\Gamma|,
\]
and assume that \ref{hypfixedthetaconditiononPtheta} holds for
some $V\dvtx \mathsf{X}\rightarrow[1,\infty)$ such that for some $\beta
\in(0,1)$,
we have $V^{\beta}(x)\geq1+|x|^{2+\epsilon}$ for all $x\in\mathsf{X}$.
Let $\gamma^{+}\in(0,1)$. Then there exists $C>0$ such that for
any $\gamma\in(0,\gamma^{+}]$, and any $\theta,x\in\Theta\times
\mathsf{X}$
\begin{eqnarray*}
P_{\theta,\gamma}w(\theta,x)&\leq& w(\theta)\\
&&{}-\gamma w(\theta)\Delta \biggl(w(
\theta)^{-\epsilon/(2+\epsilon)}+\frac{V^{\beta
}(x)\mathbb{I}\{x\notin\mathsf{C}\}+b^{\beta}(\theta)\mathbb{I}\{
x\in\mathsf{C}\}}{w(\theta)} \biggr),
\end{eqnarray*}
where $b\dvtx \Theta\rightarrow[0,\infty)$ is as in \textup{\ref
{hypfixedthetaconditiononPtheta}}
and $\Delta\dvtx [0,\infty)\rightarrow\mathbb{R}$
\[
\Delta(z):=1-C\bigl[z+z^{1/(2+\epsilon)}\bigr].
\]
\end{prop}

\begin{pf}
For $(x_{+},\mu,\Gamma)\in\mathsf{X}\times\mathbb
{R}^{n_{x}}\times\mathcal{C}$
and $\gamma\in[0,1]$, let $\mu_{+}:=\mu+\gamma[x_{+}-\mu]$ and
$\Gamma_{+}:=\Gamma+\gamma[(\mu-x_{+})(\mu-x_{+})^{\mathsf
{T}}-\Gamma]$.
We have the two trivial inequalities
\begin{eqnarray*}
|\mu_{+}| & =&\bigl|\mu+\gamma[x_{+}-\mu]\bigr|\leq(1-\gamma)|\mu|+
\gamma |x_{+}|,
\\
|\Gamma_{+}| & =&\bigl|\Gamma+\gamma\bigl[(\mu-x_{+}) (
\mu-x_{+})^{\mathsf
{T}}-\Gamma\bigr]\bigr|\\
&\leq&(1-\gamma)|\Gamma|+\gamma\bigl|(
\mu-x_{+}) (\mu -x_{+})^{\mathsf{T}}\bigr|
\end{eqnarray*}
which imply that with $w(\theta):=1+|\mu|^{2+\epsilon}+|\Gamma|$,
denoting $\bar{\gamma}:=\gamma/(1-\gamma)<1/(1-\gamma^{+})$,
\begin{eqnarray*}
&&w(\theta_{+})-w(\theta) \\
&&\qquad  \leq-|\mu|^{2+\epsilon}+\gamma \bigl[-|
\Gamma|+|\mu-x_{+}|^{2} \bigr]\\
&&\qquad \quad {}+(1-\gamma)^{2+\epsilon}|\mu
|^{2+\epsilon} \biggl[1+\frac{\gamma}{1-\gamma}\frac{|x_{+}|}{|\mu
|}
\biggr]^{2+\epsilon}
\\
&&\qquad  \leq\gamma \bigl[-|\Gamma|+|\mu-x_{+}|^{2} \bigr]+|\mu
|^{2+\epsilon} \bigl[(1-\gamma) \bigl(1+\bar{\gamma}|x_{+}|/|\mu |
\bigr)^{2+\epsilon}-1 \bigr]
\\
&&\qquad  \leq\gamma \bigl[-w(\theta)+1+|\mu|^{2+\epsilon}+2 \bigl(|\mu
|^{2}+|x_{+}|^{2} \bigr) \bigr]
\\
&&\qquad \quad {} +|\mu|^{2+\epsilon} \bigl[-\gamma \bigl(1+\bar{\gamma }|x_{+}|/|\mu|
\bigr)^{2+\epsilon}+ \bigl(1+\bar{\gamma }|x_{+}|/|\mu| \bigr)^{2+\epsilon}-1
\bigr].
\end{eqnarray*}
By the mean value theorem,
\begin{eqnarray*}
|\mu|^{2+\epsilon} \bigl[ \bigl(1+\bar{\gamma}|x_{+}|/|\mu|
\bigr)^{2+\epsilon}-1 \bigr] & \leq&|\mu|^{2+\epsilon}(2+\epsilon)\bar {
\gamma}|x_{+}|/|\mu|\times \bigl(1+\bar{\gamma}|x_{+}|/|\mu |
\bigr)^{1+\epsilon}
\\
& =&\gamma\frac{2+\epsilon}{1-\gamma}\times|\mu|^{1+\epsilon
}|x_{+}| \bigl(1+\bar{
\gamma}|x_{+}|/|\mu| \bigr)^{1+\epsilon},
\end{eqnarray*}
and since $|\mu|^{2+\epsilon}[1- (1+\bar{\gamma}|x_{+}|/|\mu
| )^{2+\epsilon}]\leq0$
we obtain the following bound:
\begin{eqnarray*}
w(\theta_{+})-w(\theta)&\leq&\gamma w(\theta) \biggl[-1+
\frac
{1}{w(\theta)}+\frac{2|\mu|^{2}}{w(\theta)} \bigl(1+|x_{+}|^{2}/|
\mu|^{2} \bigr)\\
&&\hphantom{\gamma w(\theta) \biggl[}{}+\frac{(2+\epsilon)}{1-\gamma}\frac
{|\mu|^{1+\epsilon}|x_{+}|}{w(\theta)} \bigl(1+\bar{\gamma
}|x_{+}|/|\mu| \bigr)^{1+\epsilon} \biggr]
\\
&\leq&\gamma w(\theta) \bigl(-1+C\Psi(\theta,x_{+}) \bigr),
\end{eqnarray*}
for some $C\in(0,\infty)$ and where
\[
\Psi(\theta,x_{+}):= \biggl(\frac{|\mu|^{2}}{w(\theta)}+\frac
{1+|x_{+}|^{2}}{w(\theta)}
\biggr)+ \biggl(\frac{|\mu|\times
|x_{+}|^{1/(1+\epsilon)}}{w^{1/(1+\epsilon)}(\theta)}+\frac
{|x_{+}|^{1+1/(1+\epsilon)}}{w^{1/(1+\epsilon)}(\theta)} \biggr)^{1+\epsilon}.
\]
Now from Jensen's inequality we have the identity $(a+b)^{1+\epsilon
}\leq2^{\epsilon}(a^{1+\epsilon}+b^{1+\epsilon})$
for $a,b>0$ and the following equalities:
\begin{eqnarray*}
\frac{|\mu|\times|x_{+}|^{1/(1+\epsilon)}}{w^{1/(1+\epsilon
)}(\theta)} & =&\frac{|\mu|}{w^{1/(2+\epsilon)}(\theta)} \biggl(\frac{|x_{+}|^{2+\epsilon}}{w(\theta)}
\biggr)^{1/[(1+\epsilon
)(2+\epsilon)]},
\\
\frac{|x_{+}|^{1+1/(1+\epsilon)}}{w^{1/(1+\epsilon)}(\theta)} & =& \biggl(\frac{|x_{+}|^{2+\epsilon}}{w(\theta)} \biggr)^{1/(1+\epsilon)}
\end{eqnarray*}
yield
\begin{eqnarray*}
&&\biggl(\frac{|\mu|\times|x_{+}|^{1/(1+\epsilon)}}{w^{1/(1+\epsilon
)}(\theta)}+\frac{|x_{+}|^{1+1/(1+\epsilon)}}{w^{1/(1+\epsilon
)}(\theta)} \biggr)^{1+\epsilon}\\
&&\qquad \leq2^{\epsilon} \biggl( \biggl(\frac
{|\mu|^{2+\epsilon}}{w(\theta)} \biggr)^{\vafrac{1+\epsilon
}{2+\epsilon}}
\biggl(\frac{|x_{+}|^{2+\epsilon}}{w(\theta)} \biggr)^{1/(2+\epsilon)}+\frac{|x_{+}|^{2+\epsilon}}{w(\theta)} \biggr).
\end{eqnarray*}
Note also that since $|\mu|^{2+\epsilon}/w(\theta)\leq1$, we have
\[
\frac{|\mu|^{2}}{w(\theta)}=\frac{ (|\mu|^{2+\epsilon}
)^{2/(2+\epsilon)}}{w(\theta)}\leq w(\theta)^{-\epsilon/(2+\epsilon)},
\]
and we therefore deduce that if for any $x\in\mathsf{X}$, $V^{\beta
}(x)\geq1+|x|^{2+\epsilon}$,
then
\begin{eqnarray*}
\Psi(\theta,x_{+}) & \leq& w(\theta)^{-\epsilon/(2+\epsilon)}+
\frac
{1+|x_{+}|^{2}}{w(\theta)}+2^{\epsilon} \biggl( \biggl(\frac
{|x_{+}|^{2+\epsilon}}{w(\theta)}
\biggr)^{1/(2+\epsilon)}+\frac
{|x_{+}|^{2+\epsilon}}{w(\theta)} \biggr)
\\
& \leq& w(\theta)^{-\epsilon/(2+\epsilon)}+2\frac{V^{\beta
}(x_{+})}{w(\theta)}+2^{\epsilon} \biggl(
\biggl(\frac{V^{\beta
}(x_{+})}{w(\theta)} \biggr)^{1/(2+\epsilon)}+\frac{V^{\beta
}(x_{+})}{w(\theta)} \biggr).
\end{eqnarray*}
Now by \ref{hypfixedthetaconditiononPtheta} and Jensen's inequality
we deduce that for some constant $C>1$
\begin{eqnarray*}
&&P_{\gamma,\theta}\Psi(\theta,x)
\\
&&\qquad \leq w(\theta)^{-\epsilon/(2+\epsilon)}+C \biggl[ \biggl(\frac
{V^{\beta}(x)\mathbb{I}\{x\notin\mathsf{C}\}+b^{\beta}(\theta
)\mathbb{I}\{x\in\mathsf{C}\}}{w(\theta)}
\biggr)^{1/(2+\epsilon
)}\\
&&\hspace*{147pt}{}+\frac{V^{\beta}(x)\mathbb{I}\{x\notin\mathsf{C}\}+b^{\beta
}(\theta)\mathbb{I}\{x\in\mathsf{C}\}}{w(\theta)} \biggr],
\end{eqnarray*}
and we conclude.
\end{pf}

\subsubsection{Multivariate case and superexponential tails: \texorpdfstring{\protect\ref{hyppropertiesPI}}{(A6)}}
\mbox{}

\begin{pf*}{Proof of Theorem~\ref{thmCVAMalgorithms} under \ref
{hyppropertiesPI}}
Let $\epsilon>0$, $\beta\in(0,1/(1+n_{x}/2)]$ and $V(x)\propto\pi
^{-\eta}(x)$
for some $\eta\in(0,1)$ where the constant of proportionality is
such that $V^{\beta}(x)\geq1+|x|^{2+\epsilon}$ [which is possible
as $\pi^{-\eta}(x)\geq C_{1}\exp(C_{2}|x|)$ for some $C_{1},C_{2}>0$].
From Proposition~\ref{propsimultaneousdriftsuperexponential} there
exists $a',b'>0$ and $\mathsf{C}:=B(0,R)$ for some $R>0$ such that
for any $x\in\mathsf{X}$ with $w(\theta):=1+|\mu|^{2+\epsilon
}+|\Gamma|$
and appropriate $a',a''>0$,
\begin{eqnarray*}
P_{\theta}V(x) & \leq&\bigl(1-a'/|\Gamma+
\epsilon_{\mathrm{AM}}I_{n_{x}\times
n_{x}}|^{n_{x}/2}\bigr)V(x)+b'
\mathbb{I}\{x\in\mathsf{C}\}
\\
& \leq& \bigl[1-a''/ \bigl(|\epsilon_{\mathrm{AM}}I_{n_{x}\times
n_{x}}|^{n_{x}/2}+w^{n_{x}/2}(
\theta) \bigr) \bigr]V(x)\mathbb{I}\{ x\notin\mathsf{C}\}+b\mathbb{I}\{x\in
\mathsf{C}\},
\end{eqnarray*}
where $b=b'+\sup_{x\in\mathsf{C}}V(x)$. Naturally here $\iota=1$.
Now from Proposition~\ref{propwforAM} we have
\[
P_{\theta,\gamma}w(\theta,x)\leq w(\theta)-\gamma w(\theta)\Delta \biggl(w(
\theta)^{-\epsilon/(2+\epsilon)}+\frac{V^{\beta
}(x)\mathbb{I}\{x\notin\mathsf{C}\}+b^{\beta}\mathbb{I}\{x\in
\mathsf{C}\}}{w(\theta)} \biggr)
\]
with $\Delta(z):=1-C[z+z^{1/(2+\epsilon)}]$. Therefore here $c(\theta
)=w(\theta)^{-\epsilon/(2+\epsilon)}$,
$d(\theta)=w(\theta)^{-\epsilon/(2+\epsilon)}+b^{\beta}/w(\theta)$
and $e(\theta)=w(\theta)$. We note that $p_{\Delta}=1\leq
n_{x}/2+1\leq\iota/\beta$.
The condition $\beta\leq1/(1+n_{x}/2)$ implies \ref
{hypwdriftSA1}\ref{hypenurelationawV}
as for any $\epsilon>0$, on $V^{\beta}(x)/e(\theta)\geq\epsilon$
\[
\frac{a(\theta)w(\theta)e^{-1}(\theta)}{V^{\iota-\beta}(x)}\leq C\frac{w^{n_{x}/2}(\theta)}{V^{\iota-\beta}(x)}\leq\frac
{C\epsilon^{-n_{x}/2}}{V^{\iota-\beta(1+n_{x}/2)}(x)},
\]
for some $C>0$, and we conclude with Theorem~\ref
{propestablishcompositedrift}.
\end{pf*}

\subsubsection{Relaxed tail conditions, univariate scenario: \texorpdfstring{\protect\ref{hyppismallthetasubgeomdrift}}{(A7)}}

Now we draw the same conclusions when $\mathsf{X}=\mathbb{R}$, and
$\pi(\cdot)$ now satisfies less stringent tail conditions.

\begin{pf*}{Proof of Theorem~\ref{thmCVAMalgorithms} under \ref
{hyppismallthetasubgeomdrift}}
Let $\iota,\eta\in(0,1)$ and $\beta\in(0,\iota/2]$. Let
$V(x)\propto\pi^{-\eta}(x)$,
such that $V^{\beta}(x)\geq1+|x|^{2+\epsilon}$ [which is possible
as\break $\pi^{-1}(x)\geq C_{1}\exp(C_{2}|x|^{p})$ for some $C_{1},C_{2}>0$
from Lemma~\ref{lempropertiesUpsilonetc}]. From Theorem~\ref
{thmdriftconditionsimpletailcondition},
there exist $b,R>0$ such that with $\mathsf{C}=B(0,R)$ for any $\theta,x\in\Theta\times\mathsf{X}$,
\[
P_{\theta}V(x)\leq \bigl[V(x)-a^{-1}(\theta)V^{\iota}(x)
\bigr]\mathbb{I}\{x\in\mathsf{C}\}+b\mathbb{I}\{x\in\mathsf{C}\},
\]
with $a^{-1}(\theta)=a_{0}/[(\Gamma+\epsilon_{\mathrm{AM}})^{-1}\vee(\Gamma
+\epsilon_{\mathrm{AM}})^{1/2}]\geq a_{0}/[\epsilon_{\mathrm{AM}}^{-1}\vee
(\epsilon_{\mathrm{AM}}+w(\theta) )]$
and $w(\theta):=1+|\mu|^{2+\epsilon}+|\Gamma|$. From Proposition~\ref{propwforAM} we therefore have
\[
P_{\theta,\gamma}w(\theta,x)\leq w(\theta)-\gamma w(\theta)\Delta \biggl(w(
\theta)^{-\epsilon/(2+\epsilon)}+\frac{V^{\beta
}(x)\mathbb{I}\{x\notin\mathsf{C}\}+b^{\beta}\mathbb{I}\{x\in
\mathsf{C}\}}{w(\theta)} \biggr)
\]
with $\Delta(z):=1-C[z+z^{1/(2+\epsilon)}]$. Therefore here $c(\theta
)=w(\theta)^{-\epsilon/(2+\epsilon)}$,
$d(\theta)=w(\theta)^{-\epsilon/(2+\epsilon)}+b^{\beta}/w(\theta)$
and $e(\theta)=w(\theta)$. Note that we have $p_{\Delta}=1\leq2\leq
\iota/\beta$.
The condition $\beta\leq\iota/2$ implies \ref
{hypwdriftSA1}\ref{hypenurelationawV}
as for any $\epsilon>0$, on $V^{\beta}(x)/e(\theta)\geq\epsilon$
\[
\frac{a(\theta)w(\theta)e^{-1}(\theta)}{V^{\iota-\beta}(x)}\leq C\frac{w(\theta)}{V^{\iota-\beta}(x)}\leq\frac{C\epsilon
^{-1}}{V^{\iota-2\beta}(x)},
\]
for some $C>0$ and we conclude with Theorem~\ref
{propestablishcompositedrift}.
\end{pf*}

\subsection{Stability of the coerced acceptance probability
algorithms}
\label{subStability-of-coerced}

In this subsection we establish the stability of Algorithm \ref{algstandardcoerced}
and Algorithm \ref{algfastcoerced} in a univariate setting. We proceed
as in Section~\ref{subStability-of-AM} and aim to apply Theorem~\ref{propestablishcompositedrift} and its corollary which require
\ref{hypfixedthetaconditiononPtheta} and \ref{hypwdriftSA1}
to be satisfied. A related result has been established in \cite{vihola:2011}
under a more stringent condition on the decay of the tails of the
target density, and not covering constant stepsize sequences $\{\gamma
_{i}\}$.

%
\begin{thm}
\label{thmCVcoercedalgorithms}Consider the controlled MC as defined
by either Algorithm \ref{algstandardcoerced} or Algorithm \ref{algfastcoerced}
for some $\alpha_{*}\in(0,1/2)$. Assume that $\pi(\cdot)$ and
$q(\cdot)$
satisfy \textup{\ref{hyppismallthetasubgeomdrift}} and that the stepsize
sequence $\{\gamma_{i}\}$ is such that
\[
\Bigl(\limsup_{i\rightarrow\infty}\gamma_{i+1}^{-1}-\gamma
_{i}^{-1}\Bigr)+\limsup_{i\rightarrow\infty}
\gamma_{i}<\alpha_{*}\wedge \biggl(\frac{1}{2}-
\alpha_{*}\biggr).
\]
Let $w(\theta)=\exp(|\theta|)$ for Algorithm \ref{algstandardcoerced}
and $w(\theta):=1+|\theta|^{2}$ for Algorithm \ref{algfastcoerced}.
Then there exist $M,R>0$ such that for any $\theta,x\in\Theta\times
\mathsf{X}$,
\[
\mathbb{P}_{\theta,x} \Biggl(\bigcap_{k=0}^{\infty}
\bigcup_{i\geq k}\bigl\{ (\theta_{i},X_{i})\in
\mathcal{W}_{M}\dvtx \times B(0,R)\bigr\} \Biggr)=1,
\]
that is, $\mathbb{P}_{\theta,x}\mbox{-a.s. } \{\theta_{i},X_{i}\}$ visits
$\mathcal{W}_{M}\times B(0,R)$ infinitely often.
\end{thm}

The proofs are given for the two scenarios in the following two subsections.
The following lemma, whose proof can be found in Appendix~\ref
{secProofs-from-Section},
will be useful in both scenarios.

%
\begin{lem}
\label{propsub-expon-w-drift}Assume that $\pi(\cdot)$ is a strictly
positive, differentiable probability density satisfying \ref
{hyppismallthetasubgeomdrift}\ref{subhypsubgeometrictailsofpi}.
Moreover, suppose that $q_{\sigma}(z):=\sigma^{-1}q(z/\sigma)$ where
$q\dvtx \mathsf{Z}\rightarrow[0,\bar{q}]$ for $\bar{q}>0$ is symmetric
and such that it has a finite absolute first order moment. Then,
for any $x\in\mathsf{X}$
\[
\alpha_{\sigma}(x):=\int_{\mathsf{Z}}\min \biggl\{ 1,
\frac{\pi
(x+z)}{\pi(x)} \biggr\} q_{\sigma}(z)\,\mathrm{d}z,
\]
there exist constants $C_{-},C_{+}>0$ such that
\begin{eqnarray*}
\alpha_{\sigma}(x) & \ge&1/2-C_{-}\sigma\qquad  \mbox{for
\ensuremath {\sigma\le1} and \ensuremath{x\in\mathsf{X}}},
\\
\alpha_{\sigma}(x) & \le& C_{+}\frac{(-\log\pi(x))^{1/p}\vee
1}{\sigma} \qquad  \mbox{for
\ensuremath{\sigma\ge1} and \ensuremath {x\in\mathsf{X}}.}
\end{eqnarray*}
\end{lem}

\begin{rem}
Notice from the proof that the moment condition is assumed here in
order to simplify our statement and that more general conditions are
possible. Note that the restriction $\alpha_{*}\in(0,1/2)$ is practically
harmless since this covers relevant values according to the scaling
theory of the RWM \cite{gelman:roberts:gilks:1995}.
\end{rem}

\subsubsection{Proof in the standard scenario: Algorithm \texorpdfstring{\protect\ref{algstandardcoerced}}{3}}
\label{subThe-coerced-standard-scenario}

Before starting the proofs it is worth stressing on the fact that
throughout
\[
P_{\theta}\bigl(x,y;\mathrm{d}x'\times\mathrm{d}y'
\bigr)=q\bigl(x,\mathrm {d}y'\bigr)\bigl[\alpha\bigl(x,y'
\bigr)\delta_{y'}\bigl(\mathrm{d}x'\bigr)+ \bigl(1-\alpha
\bigl(x,y'\bigr) \bigr)\delta_{x}\bigl(
\mathrm{d}x'\bigr)\bigr]
\]
and hence that for any $x,y\in\mathsf{X}$, $P_{\theta}(x,y;\cdot
)=P_{\theta}(x,\cdot)$
and that for notational simplicity the Lyapunov function $V(x)$ should
be understood as being the function $V(x)\times1$ defined on $\mathsf{X}^{2}$.
The following proposition establishes part of \ref{hypwdriftSA1}
under a condition implied by Lemma~\ref{propsub-expon-w-drift}.

%
\begin{prop}
\label{propcoercedconditionacceptprob}Consider the controlled
MC as defined in Algorithm \ref{algstandardcoerced} with $\alpha_{*}\in(0,1/2)$,
$V(x):=c\pi^{-\eta}(x)$ and assume that there exist $C>0$ and $\beta
\in[0,1)$
such that for all $(\theta,x)\in\Theta\times\mathsf{X}$
\[
\sgn(\theta) \bigl(\alpha_{\exp(\theta)}(x)-\alpha_{*}\bigr)\leq-
\bigl[\alpha _{*}\wedge(1/2-\alpha_{*})\bigr]+CV^{\beta}(x)/
\exp\bigl(|\theta|\bigr).
\]
Let $\gamma_{\max}\in(0,\alpha_{*}\wedge(1/2-\alpha_{*}))$. Then
for any $\gamma\in(0,\gamma_{\max}]$ and $\theta,x\in\Theta
\times\mathsf{X}$
and the Lyapunov function $w\dvtx \Theta\to[1,\infty)$ defined by
$w(\theta):=\exp (|\theta| )$
\begin{eqnarray*}
P_{\gamma,\theta}w(\theta,x)&\le& w(\theta)\\
&&{}-\gamma w(\theta)\Delta \biggl(
\mathbb{I}\bigl\{|\theta|\leq\gamma_{\max}\bigr\}C^{-1}\biggl(2+\alpha
_{*}\wedge\biggl(\frac{1}{2}-\alpha_{*}\biggr)
\biggr)+\frac{V^{\beta}(x)}{w(\theta
)} \biggr),
\end{eqnarray*}
with
\[
\Delta (z )= \bigl(\alpha_{*}\wedge\bigl(\tfrac{1}{2}-\alpha
_{*}\bigr) \bigr)-\gamma_{\max}-Cz.
\]
\end{prop}

\begin{pf}
For $|\theta|>\gamma$ and since for all $x,y_{+}\in\mathsf{X}^{2}$
$|\alpha(x,y_{+})-\alpha_{*}|\le1$, one can write (with $\theta
_{+}=\theta+\gamma[\alpha(x,y_{+})-\alpha_{*}]$)
%
\begin{eqnarray*}
w(\theta_{+}) & =&w\bigl(\theta+\gamma\bigl[\alpha (x,y_{+})-
\alpha_{*}\bigr]\bigr)
\\
& =&\exp\bigl(|\theta|+\sgn(\theta)\gamma\bigl[\alpha(x,y_{+})-
\alpha_{*}\bigr]\bigr).
\end{eqnarray*}
Now since $\gamma\le1$, from the inequality $\exp(u)\le1+u+u^{2}$
valid for $|u|\le1$, one obtains
\[
w(\theta_{+})\le w(\theta) \bigl(1+\gamma\sgn(\theta)\bigl[\alpha
(x,y_{+})-\alpha_{*}\bigr]+\gamma^{2}\bigr).
\]
Taking the conditional expectations yields for $|\theta|>\gamma$
and by assumption
\begin{eqnarray*}
P_{\gamma,\theta}w(\theta,x) & \le& w(\theta) \bigl(1+\gamma\sgn(\theta )\bigl[
\alpha_{\exp(\theta)}(x)-\alpha_{*}\bigr]+\gamma^{2}\bigr)
\\
& \leq& w(\theta)-\gamma w(\theta) \bigl(\bigl[\alpha_{*}\wedge (1/2-
\alpha_{*})\bigr]-CV^{\beta}(x)/\exp\bigl(|\theta|\bigr)-
\gamma_{\max
} \bigr).
\end{eqnarray*}
Also notice that for any $\theta\in\Theta$ we have $|\theta_{+}|\le
|\theta|+\gamma$,
whence $w(\theta_{+})\le w(\theta)\exp(\gamma)\le w(\theta
)(1+\gamma+\gamma^{2})\le w(\theta)(1+2\gamma)$
for all $\gamma\leq\gamma_{\mathrm{max}}\le1$. From this inequality
and the display above, we deduce for all $\theta,x\in\Theta\times
\mathsf{X}$
and $\gamma\in(0,\gamma_{\max}]$,
\begin{eqnarray*}
&&P_{\gamma,\theta}w(\theta,x)\le w(\theta)-\gamma w(\theta) \bigl[ \bigl(
\alpha_{*}\wedge\bigl(\tfrac{1}{2}-\alpha_{*}\bigr)
\bigr)-\gamma _{\max}\\
&&\hphantom{P_{\gamma,\theta}w(\theta,x)\le w(\theta)-\gamma w(\theta) \biggl[}{}-\mathbb{I}\bigl\{|\theta|\leq\gamma_{\max}\bigr\}
\bigl(2+\alpha _{*}\wedge\bigl(\tfrac{1}{2}-
\alpha_{*}\bigr)\bigr)\\
&&\hspace*{220pt}{}-CV^{\beta}(x)/w(\theta) \bigr].
\end{eqnarray*}
\upqed
\end{pf}

\begin{pf*}{Proof of Theorem~\ref{thmCVcoercedalgorithms} in the
case of Algorithm \ref{algstandardcoerced}}
First notice that there exists $i_{0}\in\mathbb{N}$ such that $\sup_{i\geq i_{0}}\gamma_{i+1}^{-1}-\gamma_{i}^{-1}<\alpha_{*}\wedge
(\frac{1}{2}-\alpha_{*})-\gamma_{\max}$
with $\gamma_{\max}:=\sup_{i\geq i_{0}}\gamma_{i}$. We show that
\ref{hypfixedthetaconditiononPtheta} and \ref{hypwdriftSA1}
are satisfied and conclude with Theorem~\ref{propestablishcompositedrift}
and Corollary~\ref{corcorollaryrecurrence}, for $i\geq i_{0}$.
Let $\eta,\iota\in(0,1)$, $\beta\in(0,\iota/3]$, and define
$V(x)\propto\pi^{-\eta}(x)$,
such that $V^{\beta}(x)\geq1\vee(-\log\pi(x))^{1/p}$ (which is possible
since for any $a_{1},a_{2},M>0$, $\sup_{0\leq u\leq M}u^{a_{1}}|\log
u|^{a_{2}}<\infty$).
From Theorem~\ref{thmdriftconditionsimpletailcondition}, there exist
$b,R>0$ such that for any $\theta,x\in\Theta\times\mathsf{X}$,
\[
P_{\theta}V(x)\leq \bigl[V(x)-a^{-1}(\theta)V^{\iota}(x)
\bigr]\mathbb{I}\{x\notin\mathsf{C}\}+b\mathbb{I}\{x\in\mathsf{C}\},
\]
with $a(\theta)=[\exp(\theta)\vee\exp(-2\theta)]/a_{0}$ (for some
$a_{0}>0$) and $\mathsf{C}=B(0,R)$. Now with $w(\theta)=\exp(|\theta|)$
from Lemma~\ref{propsub-expon-w-drift} there exists $C>0$ such
that
\begin{eqnarray*}
\alpha_{\exp(\theta)}(x) & \geq&1/2-C\frac{V^{\beta}(x)}{\exp
(-\theta)} \qquad  \mbox{for
\ensuremath{\theta\le0} and \ensuremath {x\in\mathsf{X}}},
\\
\alpha_{\exp(\theta)}(x) & \le& C_{+}\frac{(-\log\pi(x))^{1/p}\vee
1}{\exp(\theta)}\leq C
\frac{V^{\beta}(x)}{\exp(\theta)} \qquad  \mbox{for \ensuremath{\theta\ge0} and \ensuremath{x\in
\mathsf{X}}.}
\end{eqnarray*}
One can apply Proposition~\ref{propcoercedconditionacceptprob},
leading to the existence of $C>0$ such that for any $\theta,x\in
\Theta\times\mathsf{X}$,
\[
P_{\gamma,\theta}w(\theta,x)\le w(\theta)-\gamma w(\theta)\Delta \biggl(
\biggl[c(\theta)+\frac{V^{\beta}(x)}{w(\theta)} \biggr]\mathbb{I}\{x\notin\mathsf{C}\}+d(
\theta)\mathbb{I}\{x\in\mathsf {C}\} \biggr),
\]
with
\begin{eqnarray*}
c(\theta)&=&C^{-1}\mathbb{I}\bigl\{|\theta|\leq\gamma_{\max}\bigr\}
\bigl(2+\alpha _{*}\wedge\bigl(\tfrac{1}{2}-
\alpha_{*}\bigr)\bigr),
\\
d(\theta)&=&\frac{\sup_{x\in\mathsf{C}}V^{\beta}(x)}{w(\theta
)}+C^{-1}\mathbb{I}\bigl\{|\theta|\leq
\gamma_{\max}\bigr\}\biggl(2+\alpha _{*}\wedge\biggl(
\frac{1}{2}-\alpha_{*}\biggr)\biggr)
\end{eqnarray*}
and
\[
\Delta (z )=\alpha_{*}\wedge\bigl(\tfrac{1}{2}-\alpha
_{*}\bigr)-\gamma_{\max}-Cz.
\]
Notice that from our choice of $i_{0} \sup_{i\geq i_{0}}\gamma
_{i+1}^{-1}-\gamma_{i}^{-1}<\Delta(0)$
and that we have $\iota/\beta\geq3>p_{\Delta}=1$ in \ref
{hypwdriftSA1}.
Clearly here $e(\theta)=w (\theta )$. Now the condition
$\beta\leq\iota/3$ implies \ref{hypwdriftSA1}\ref
{hypenurelationawV}
as for any $\epsilon>0$, on $V^{\beta}(x)/e(\theta)\geq\epsilon$
\[
\frac{a(\theta)w(\theta)e^{-1}(\theta)}{V^{\iota-\beta}(x)}=\frac
{\exp(2|\theta|)}{a_{0}V^{\iota-\beta}(x)}\leq\frac{\epsilon
^{-2}a_{0}^{-1}}{V^{\iota-3\beta}(x)}<\infty,
\]
and we conclude.
\end{pf*}

\subsubsection{Proof for the accelerated version: Algorithm \texorpdfstring{\protect\ref{algfastcoerced}}{4}}

The arguments are similar to those of Section~\ref{subThe-coerced-standard-scenario},
but here $w(\theta)$ is here of a different form. The following proposition
is similar to Proposition~\ref{propcoercedconditionacceptprob}
but takes this change of Lyapunov function into account.

%
\begin{prop}
\label{propcoercedconditionacceptprobfast}Consider Algorithm \ref
{algfastcoerced}
with $\alpha_{*}\in(0,1/2)$, let $w(\theta):=1+|\theta|^{2}$ and
assume that there exists $C>0$ and $\beta\in[0,1)$ such that for
any $\theta,x\in\Theta\times\mathsf{X}$
\[
\sgn(\theta) \bigl(\alpha_{\exp(\theta)}(x)-\alpha_{*}\bigr)\leq-
\bigl[\alpha _{*}\wedge(1/2-\alpha_{*})\bigr]+CV^{\beta}(x)/
\exp\bigl(|\theta|\bigr).
\]
Let $\gamma_{\max}\in(0,\alpha_{*}\wedge(1/2-\alpha_{*}))$. Then
there exists $C'>0$ such that for any $\gamma\in(0,\gamma_{\max}]$
and $\theta,x\in\Theta\times\mathsf{X}$,
\begin{eqnarray*}
P_{\theta,\gamma}w(\theta,x)&\leq& w(\theta)\\
&&{}-\gamma w(\theta)\Delta \bigl(
\mathbb{I}\bigl\{|\theta|\leq1\bigr\}C'^{-1} \bigl(
\alpha_{*}\wedge (1/2-\alpha_{*}) \bigr)+V^{\beta}(x)/
\exp\bigl(|\theta|\bigr)\bigr)
\end{eqnarray*}
with
\[
\Delta(z)=2\bigl[\alpha_{*}\wedge(1/2-\alpha_{*})-
\gamma_{\max}-C'z\bigr].
\]
\end{prop}

\begin{pf}
With $\theta_{+}=\theta+\gamma(1+|\theta|)[\alpha(x,y)-\alpha_{*}]$
we have
\begin{eqnarray*}
w(\theta_{+}) & =&w(\theta)+2\gamma|\theta|\bigl(|\theta|+1\bigr)\sgn(\theta )
\bigl(\alpha(x,y)-\alpha_{\ast}\bigr)\\
&&{}+\gamma^{2}\bigl(1+|
\theta|\bigr)^{2}\bigl[\alpha (x,y)-\alpha_{\ast}
\bigr]^{2}
\\
& \leq& w(\theta)+2\gamma|\theta|\bigl(|\theta|+1\bigr)\sgn(\theta) \bigl(\alpha (x,y)-
\alpha_{\ast}\bigr)+2\gamma^{2}w(\theta)
\end{eqnarray*}
so
\begin{eqnarray*}
P_{\theta,\gamma}w(\theta,x)&\leq& w(\theta)+2\gamma|\theta |\bigl(|\theta|+1\bigr) \bigl[-
\bigl[\alpha_{*}\wedge(1/2-\alpha_{*})\bigr]+CV^{\beta
}(x)/
\exp\bigl(|\theta|\bigr) \bigr]\\
&&{}+2\gamma^{2}w(\theta).
\end{eqnarray*}
Notice that for $|\theta|\geq1$ we have $|\theta|(1+|\theta|)\geq
1+|\theta|^{2}$. Consequently for any $|\theta|\geq1$ and $x\in\mathsf{X}$ such
that $-[\alpha_{*}\wedge(1/2-\alpha_{*})]+CV^{\beta}(x)/\exp
(|\theta|)\leq0$
\begin{eqnarray*}
P_{\theta,\gamma}w(\theta,x) & \leq& w(\theta)+2\gamma w(\theta ) \bigl[-\bigl[
\alpha_{*}\wedge(1/2-\alpha_{*})\bigr]+\gamma+CV^{\beta
}(x)/
\exp\bigl(|\theta|\bigr) \bigr]
\\
& \leq& w(\theta)-2\gamma w(\theta) \bigl[\alpha_{*}\wedge (1/2-
\alpha_{*})-\gamma-2CV^{\beta}(x)/\exp\bigl(|\theta|\bigr) \bigr].
\end{eqnarray*}
Notice that for any $\theta\in\Theta$, $|\theta|(1+|\theta|)\leq
(1+|\theta|)^{2}\leq2(1+|\theta|^{2})$.
For the specific case $-[\alpha_{*}\wedge(1/2-\alpha_{*})]+CV^{\beta
}(x)/\exp(|\theta|)\geq0$,
we therefore have
\begin{eqnarray*}
P_{\theta,\gamma}w(\theta,x) & \leq& w(\theta)\\
&&{}+2\gamma2w(\theta ) \bigl[-\bigl[
\alpha_{*}\wedge(1/2-\alpha_{*})\bigr]+\gamma/2+CV^{\beta
}(x)/
\exp\bigl(|\theta|\bigr) \bigr]
\\
& \leq& w(\theta)+2\gamma w(\theta) \bigl[-\bigl[\alpha_{*}\wedge
(1/2-\alpha_{*})\bigr]+\gamma+2CV^{\beta}(x)/\exp\bigl(|\theta|\bigr)
\bigr],
\end{eqnarray*}
and for any $\theta,x\in\Theta\times\mathsf{X}$ one has
\[
P_{\theta,\gamma}w(\theta,x)\leq w(\theta)-2\gamma w(\theta )
\bigl[-2CV^{\beta}(x)/\exp\bigl(|\theta|\bigr)-\gamma\bigr].
\]
We can now combine these intermediate results, yielding for any $\theta,x\in\Theta\times\mathsf{X}$
\begin{eqnarray*}
&&P_{\theta,\gamma}w(\theta,x)\leq w(\theta)-\gamma2w(\theta) \bigl[
\alpha_{*}\wedge(1/2-\alpha_{*})-\gamma_{\max}\\
&&\hphantom{P_{\theta,\gamma}w(\theta,x)\leq w(\theta)-\gamma2w(\theta) \bigl[}{}-
\mathbb{I}\bigl\{ |\theta|\leq1\bigr\} \bigl(\alpha_{*}\wedge(1/2-
\alpha_{*}) \bigr)\\
&&\hspace*{181pt}{}-2CV^{\beta}(x)/\exp\bigl(|\theta|\bigr) \bigr],
\end{eqnarray*}
and we conclude.
\end{pf}

\begin{pf*}{Proof of Theorem~\ref{thmCVcoercedalgorithms} in the
case of Algorithm \ref{algfastcoerced}}
The beginning of the proof is similar to that of Algorithm \ref{algstandardcoerced}
by using Proposition~\ref{propcoercedconditionacceptprobfast}
and Lemma~\ref{propsub-expon-w-drift}, but here we set $\beta\in
(0,\iota/2)$.
This leads to the existence of $C>0$ such that for any $\theta,x\in
\Theta\times\mathsf{X}$,
\[
P_{\gamma,\theta}w(\theta,x)\le w(\theta)-\gamma w(\theta)\Delta \biggl(
\biggl[c(\theta)+\frac{V^{\beta}(x)}{\exp (|\theta
| )} \biggr]\mathbb{I}\{x\notin\mathsf{C}\}+d(
\theta)\mathbb {I}\{x\in\mathsf{C}\} \biggr),
\]
with
\begin{eqnarray*}
c(\theta)&=&C'^{-1}\mathbb{I}\bigl\{|\theta|\leq1\bigr\}\bigl(
\alpha_{*}\wedge \bigl(\tfrac{1}{2}-\alpha_{*}\bigr)
\bigr),
\\
d(\theta)&=&2\frac{\sup_{x\in\mathsf{C}}V^{\beta}(x)}{\exp(|\theta
|)}+C'^{-1}\mathbb{I}\bigl\{|\theta|
\leq1\bigr\}\biggl(\alpha_{*}\wedge\biggl(\frac
{1}{2}-
\alpha_{*}\biggr)\biggr)
\end{eqnarray*}
and
\[
\Delta (z )=2 \bigl[ \bigl(\alpha_{*}\wedge\bigl(
\tfrac
{1}{2}-\alpha_{*}\bigr) \bigr)-\gamma_{\max}-C'z
\bigr].
\]
Notice that from our choice of $i_{0}$, $\sup_{i\geq i_{0}}\gamma
_{i+1}^{-1}-\gamma_{i}^{-1}<\Delta(0)$
and that we have $\iota/\beta\geq2>p_{\Delta}=1$ in \ref
{hypwdriftSA1}.
Now the condition $\beta<\iota/2$ implies \ref
{hypwdriftSA1}\ref{hypenurelationawV}
as for any $\epsilon>0$ there exist $C''\in(0,\infty)$ such that
on $V^{\beta}(x)/e(\theta)\geq\epsilon$ [since here $e(\theta
)=\exp (|\theta| )$]
\[
\frac{a(\theta)w(\theta)e^{-1}(\theta)}{V^{\iota-\beta}(x)}=\frac
{\exp(|\theta|)[1+|\theta|^{2}]}{a_{0}V^{\iota-\beta}(x)}\leq C''
\frac{1+ (\log V(x) )^{2}}{V^{\iota-2\beta}(x)}<\infty,
\]
and we conclude.
\end{pf*}\vspace*{-9pt}

%
\begin{appendix}
%
\section{Appendix for Section~\texorpdfstring{\protect\ref{121212121}}{1}}\vspace*{-6pt}
\label{secAppendix-for-Section1}

\begin{pf*}{Proof of Lemma~\ref{lemgeneral-drift-stability}}
For any $k\geq1$, we introduce the stopping times $\tau(k):=\inf\{
i>k\dvtx (\theta_{i},X_{i})\in\mathcal{C}\}$.
We proceed by contradiction and observe first that if the claim did
not hold, then there would be an integer $i_{w}\leq n<\infty$ such
that with positive probability the stopping time $\tau(n)$ would
be infinite, that is, $\mathbb{P}_{\theta,x}(\tau(n)=\infty)>0$. We
establish
a result similar to \cite{meyn:tweedie:1993}, Proposition~11.3.3, page 266,
but take care of the inhomogeneity and do not require the same precision.
We introduce the following notation for simplicity:
$W_{i}:=W_{i}(\theta_{i},X_{i})$
and for any $m\in\mathbb{N}$, $\tau^{m}:=\tau(n)\wedge m$ (we omit
the dependence on $n$ in order to alleviate notation). Assumption
\eqref{eqgeneric-W-drift} implies that for $i\geq n+1$,
%
\begin{eqnarray*}
&&\mathbb{E}_{\theta,x}\bigl[W_{i+1}\mathbb{I}\bigl\{
\tau^{m}\geq i+1\bigr\}\bigr] \\
&&\qquad  =\mathbb{E}_{\theta,x}
\bigl[W_{i}\mathbb{I}\bigl\{\tau^{m}\geq i+1\bigr\} +
\mathbb{E}_{\theta,x} [W_{i+1}-W_{i}\mid\mathcal
{F}_{i}]\mathbb{I}\bigl\{\tau^{m}\geq i+1\bigr\} \bigr]
\\
&&\qquad  \le\mathbb{E}_{\theta,x}\bigl[W_{i}\mathbb{I}\bigl\{
\tau^{m}\geq i\bigr\} \bigr]-\mathbb{E}_{\theta,x}\bigl[
\delta_{i+1}\mathbb{I}\bigl\{\tau^{m}\geq i+1\bigr\}\bigr],
\end{eqnarray*}
and consequently, we can establish
\[
\mathbb{E}_{\theta,x}\Biggl[\sum_{i=n+1}^{\infty}
\delta _{i+1}\mathbb{I}\bigl\{\tau^{m}-1\geq i\bigr\}\Biggr]
\leq \mathbb{E}_{\theta,x}[W_{n+1}]-\mathbb{E}_{\theta,x}[W_{\tau^{m}}]
\leq\mathbb {E}_{\theta,x}[W_{n+1}].
\]
Now, by using the trivial inequality $\mathbb{E}_{\theta,x}[\mathbb
{I}\{\tau(n)=\infty\}{\sum}_{i=n+1}^{\infty}\delta
_{i+1}\mathbb{I}\{\tau^{m}-1\geq i\}]\leq\mathbb{E}_{\theta,x}[{\sum}_{i=n+1}^{\infty}\delta_{i+1}\mathbb{I}\{\tau
^{m}-1\geq i\}]$
and the monotone convergence theorem (thanks to our assumptions on
$\{\delta_{i}\}$) we obtain the contradictory statement
\[
\mathbb{P}_{\theta,x}\bigl(\tau(n)=\infty\bigr)\sum
_{i=n+1}^{\infty}\delta _{i+1}\leq
\mathbb{E}_{\theta,x}[W_{n+1}]<\infty.
\]
We therefore conclude that for any $i\geq i_{w}$, $\mathbb{P}_{\theta,x}(\tau(i)=\infty)=0$,
and the result follows.
\end{pf*}\vspace*{-9pt}

\section{Appendix for Section~\texorpdfstring{\protect\ref{1313131313}}{3}}
\label{secAppendix-for-Section3}

We state the following result for the reader's convenience.

\begin{thm}[(see \cite{andrieu:moulines:priouret:2005} for a proof)]
\label{theokeyproposition} Assume \textup{\ref{hypALyapunovFunction}}.
For any $M\in(M_{0},M_{1}]$ there exist $\delta_{0}>0$ and $\lambda_{0}>0$
such that, for all $n\geq1$, all $\vartheta_{0}\in\mathcal{W}_{M_{0}}$,
all\vspace*{2pt} sequences $\boldsymbol{\rho}=\{\rho_{k}\}$ of nonnegative real numbers
and all sequences $\{\varsigma_{k}\}\subset\Theta^{\mathbb{N}}$ of
$n_{\theta}$-dimensional vectors satisfying
\[
\sup_{1\leq k\leq n}\rho_{k}\leq\lambda_{0}
\quad \mbox{and}\quad  \sup_{1\leq k\leq n}\Biggl\llvert \sum
_{j=1}^{k}\rho_{j}\varsigma
_{j}\Biggr\rrvert \leq\delta_{0},
\]
we have for $k=1,\ldots,n$, $w(\vartheta_{k})\leq M$, where
$\vartheta_{k}=\vartheta_{k-1}+\rho_{k}h(\vartheta_{k-1})+\rho
_{k}\varsigma_{k}$.
\end{thm}

\section{Proofs from Section~\texorpdfstring{\protect\ref{171717171717}}{5}}
\label{secProofs-from-Section}

Before proving Lemmas \ref{lempropertiesUpsilonetc}--\ref
{lemsubgeomboundonTsthetalargerthanx}
we state and prove an intermediate result.

%
\begin{lem}
\label{lemsmallresultaboutintegralratiopi}Let $c,p>0$ be constants.
Then there exist constants $M=M(c,p)\in(0,\infty)$ and
$x_{0}=x_{0}(c,p)\in(0,\infty)$
such that
\[
\int_{0}^{\infty}\exp \bigl(-c\bigl[(x+z)^{p}-x^{p}
\bigr] \bigr)\,\mathrm{d}z\le Mx^{1-p} \qquad \mbox{for all \ensuremath{x\ge
x_{0}}}.
\]
\end{lem}

\begin{pf}
By a change of variable $u=c(x+z)^{p}$, we obtain
\[
\int_{0}^{\infty}\exp \bigl(-c\bigl[(x+z)^{p}-x^{p}
\bigr] \bigr)\,\mathrm {d}z=\frac{e^{cx^{p}}}{cp}\int_{cx^{p}}^{\infty}e^{-u}u^{\sfrac
{1}{p}-1}\,\mathrm{d}u.
\]
Integration by parts yields
%
\begin{equation}\label{eqgamma-int-by-parts}
\hspace*{15pt}\int_{cx^{p}}^{\infty}e^{-u}u^{\sfrac
{1}{p}-1}\,\mathrm{d}u=e^{-cx^{p}}
\bigl(cx^{p}\bigr)^{\sfrac{1}{p}-1}+ \biggl(\frac
{1}{p}-1 \biggr)
\int_{cx^{p}}^{\infty}e^{-u}u^{\sfrac{1}{p}-2}\,\mathrm{d}u.
\end{equation}
Now, if $p\ge1$, this is enough to yield the claim. Suppose then
$p\in(0,1)$, and fix a constant $\lambda\in(0,1)$. By \eqref
{eqgamma-int-by-parts},
\begin{eqnarray*}
(1-\lambda)\int_{cx^{p}}^{\infty}e^{-u}u^{\sfrac
{1}{p}-1}\,\mathrm{d}u&=&e^{-cx^{p}}
\bigl(cx^{p}\bigr)^{\sfrac{1}{p}-1}\\
&&{}+\int_{cx^{p}}^{\infty
}e^{-u}u^{\sfrac{1}{p}-1}
\biggl[ \biggl(\frac{1}{p}-1 \biggr)\frac
{1}{u}-\lambda \biggr]\,\mathrm{d}u.
\end{eqnarray*}
Now, if $cx^{p}\ge (\frac{1}{p}-1 )\frac{1}{\lambda}$, the
latter integrand is negative. Setting
\[
x_{0}:= \biggl[ \biggl(\frac{1}{p}-1 \biggr)
\frac{1}{\lambda c} \biggr]^{1/p},
\]
we therefore have for $x\ge x_{0}$ the desired bound
\[
\int_{0}^{\infty}\exp \bigl(-c\bigl[(x+z)^{p}-x^{p}
\bigr] \bigr)\,\mathrm{d}z\le \frac{(cx^{p})^{\sfrac{1}{p}-1}}{cp(1-\lambda)}=\frac{c^{\sfrac
{1}{p}-2}}{p(1-\lambda)}x^{1-p}.
\]
We remark that the constant $\lambda\in(0,1)$ can be used to optimize
the value constants $M$ and $x_{0}$.
\end{pf}

\begin{pf*}{Proof of Lemma~\ref{lempropertiesUpsilonetc}}
First from assumption \ref{hyppismallthetasubgeomdrift}\ref
{subhypsubgeometrictailsofpi}
there exist $R_{\ell},C_{\ell}>0$ such that for all $x\in
B^{c}(0,R_{\ell})$,
we have
\[
\ell'(x)\leq-C_{\ell}|x|^{p-1},
\]
and consequently for all $x\in B^{c}(0,R_{\ell})$ and $z\geq0$, we
have
\begin{eqnarray*}
\frac{\pi(x+\sgn(x)z)}{\pi(x)}&=&\exp \biggl(\sgn(x)\int_{0}^{z}
\ell'\bigl(x+\sgn(x)t\bigr)\,\mathrm{d}t \biggr)\\
&\leq&\exp \biggl(-
\frac{C_{\ell
}}{p} \bigl[\bigl\vert |x|+z\bigr\vert ^{p}-|x|^{p}
\bigr] \biggr).
\end{eqnarray*}
Consequently for any $x\in B^{c}(0,R_{\ell})$ we deduce that
%
\begin{equation}\label{eqsimpleboundtailofpi}
\pi(x)\leq\bigl[\pi(-R_{\ell})\vee\pi(R_{\ell})\bigr]\exp
\bigl(C_{\ell
}/p|R_{\ell}|^{p}\bigr)\exp
\bigl(-C_{\ell}/p|x|^{p}\bigr).
\end{equation}
We deduce that there exists $R_{1}\geq R_{\ell}$ such that
\[
\bigl[\pi(-R_{\ell})\vee\pi(R_{\ell})\bigr]\exp
\bigl(C_{\ell}/p|R_{\ell
}|^{p}\bigr)\exp
\bigl(-C_{\ell}/pR_{1}^{p}\bigr)\leq\inf
_{x\in B(0,R_{\ell})}\pi(x),
\]
and from $\pi(\cdot)>0$, its continuity and the fact that it is monotone
on both $(-\infty,-R_{\ell}]$ and $[R_{1},\infty)$ we deduce the
first statement. Now from (\ref{eqsimpleboundtailofpi}) we deduce
that there exists $C_{1}>0$ such that for $x\in[R_{1},\infty)$
\[
\pi(x)=\pi\bigl(\Upsilon(x)\bigr)\leq C_{1}\exp\bigl(-C_{\ell}\bigl|
\Upsilon(x)\bigr|^{p}\bigr)
\]
which implies the existence of $C_{\Upsilon,1},C_{\Upsilon,2},R_{\Upsilon}>0$
such that for any $x\in\mathsf{X}$ such that $|x|\geq R_{\Upsilon}$
%
\begin{equation}\label{eqboundUpsilon-1}
\bigl|\Upsilon(x)\bigr|\vee|x|\leq C_{\Upsilon,1} \bigl(-\log\bigl(\pi
(x)/C_{\Upsilon,2}\bigr) \bigr)^{1/p}.
\end{equation}
From above we have the upper bound
\[
I_{\gamma}(x)\leq\int_{0}^{\infty}\exp
\bigl(-C_{\ell}\gamma/p \bigl[\bigl|x+\sgn(x)z\bigr|^{p}-|x|^{p}
\bigr]\bigr)\,\mathrm{d}z.
\]
We can conclude with the result of Lemma~\ref
{lemsmallresultaboutintegralratiopi}.
We proceed similarly with $J_{\gamma}(x)$ by noticing that $\ell
(x)-\ell(x-\sgn(x)z)=\sgn(x)\int_{-z}^{0}\ell'(x+\sgn(x)t)\,\mathrm{d}t\leq
-C_{\ell}/p[|x|^{p}-|x-\sgn(x)z|^{p}]$
and again conclude with Lemma~\ref{lemsmallresultaboutintegralratiopi}
above.
\end{pf*}

\begin{pf*}{Proof of Lemma~\ref{lemdecompPVoverVsubgeometric}}
Let $\eta\in(0,1)$, and consider
\begin{eqnarray*}
P_{\sigma}V(x) & = & \int_{\mathsf{X}}V(y)\min \biggl\{ 1,
\frac{\pi
(y)}{\pi(x)} \biggr\} q_{\sigma}(x,y)\,\mathrm{d}y\\
&&{}+V(x)\int
_{\mathsf
{X}} \biggl(1-\min \biggl\{ 1,\frac{\pi(y)}{\pi(x)} \biggr\}
\biggr)q_{\sigma}(x,y)\,\mathrm{d}y
\\
& = & \int_{A_{x}}V(y)q_{\sigma}(x,y)\,\mathrm{d}y+\int
_{R_{x}}V(y)\frac{\pi(y)}{\pi(x)}q_{\sigma}(x,y)\,\mathrm
{d}y\\
&&{}+V(x)\int_{R_{x}} \biggl(1-\frac{\pi(y)}{\pi(x)}
\biggr)q_{\sigma}(x,y)\,\mathrm{d}y,
\end{eqnarray*}
where $A_{x}:=\{y\in\mathbb{R}\dvtx \pi(y)\ge\pi(x)\}$ and $R_{x}:=\{
y\in\mathbb{R}\dvtx \pi(y)<\pi(x)\}$
are the regions of (almost) sure acceptance and possible rejection,
respectively. From this expression, we obtain\vspace*{-1pt}
%
\begin{eqnarray}\label{eqPVoverVminus1decompositin}
\frac{P_{\sigma}V(x)}{V(x)}-1 & = & \int_{A_{x}} \biggl(
\frac
{V(y)}{V(x)}-1 \biggr)q_{\sigma}(x,y)\,\mathrm{d}y
\nonumber
\\
& &{} +\int_{R_{x}} \biggl[ \biggl(\frac{V(y)}{V(x)}
\frac{\pi(y)}{\pi
(x)}-1 \biggr)+ \biggl(1-\frac{\pi(y)}{\pi(x)} \biggr)
\biggr]q_{\sigma}(x,y)\,\mathrm{d}y
\nonumber
\\[-8pt]\\[-8pt]
& = & \int_{A_{x}} \biggl[ \biggl(\frac{\pi(y)}{\pi(x)}
\biggr)^{-\eta}-1 \biggr]q_{\sigma}(x,y)\,\mathrm{d}y
\nonumber
\\
& & {}+\int_{R_{x}} \biggl\{ \biggl[ \biggl(\frac{\pi(y)}{\pi(x)}
\biggr)^{1-\eta}-1 \biggr]+ \biggl[1-\frac{\pi(y)}{\pi(x)} \biggr] \biggr\}
q_{\sigma}(x,y)\,\mathrm{d}y.\nonumber
\end{eqnarray}
Notice that thanks to \ref{hyppismallthetasubgeomdrift} and
Lemma~\ref{lempropertiesUpsilonetc} $\lim_{x\rightarrow\infty
}\Upsilon(x)=-\infty$
and that for $R$ sufficiently large, for any $x\geq R$ we have that
$A_{x}=[\Upsilon(x),x]$ and $R_{x}=(-\infty,\Upsilon(x))\cup
(x,\infty)$.
Then with $y=x\pm z$ and by taking into account that the support
of $q_{\sigma}(z)$ is included in $[-\sigma,\sigma]$, we have\vspace*{-1pt}
\begin{eqnarray*}
\frac{P_{\sigma}V(x)}{V(x)}-1&=&\int_{0}^{(x-\Upsilon(x))\wedge
\sigma} \bigl(
\phi_{x,-\eta,-1}(z)-1 \bigr)q_{\sigma}(z)\,\mathrm {d}z\\
&&{}+\int
_{0}^{\sigma} \bigl[\bigl(\phi_{x,1-\eta,1}(z)-1 \bigr)-
\bigl(\phi _{x,1,1}(z)-1 \bigr) \bigr]q_{\sigma}(z)\,\mathrm{d}z
\\
&&{}+\mathbb{I}\bigl\{\sigma\geq x-\Upsilon(x)\bigr\}\\
&&\hphantom{{}+}{}\times\int_{x-\Upsilon
(x)}^{\sigma}
\bigl[\bigl(\phi_{x,1-\eta,-1}(z)-1 \bigr)-\bigl(\phi _{x,1,-1}(z)-1 \bigr)
\bigr]q_{\sigma}(z)\,\mathrm{d}z
\end{eqnarray*}
and therefore, because $x-\Upsilon(x)>x$, we may write
\begin{eqnarray*}
\frac{P_{\sigma}V(x)}{V(x)}-1&=&\int_{0}^{\sigma\wedge x}\psi
_{x}(z)q_{\sigma}(z)\,\mathrm{d}z\\
&&{}+\mathbb{I}\{\sigma\geq x\}\int
_{x}^{\sigma}\bigl[ \bigl(\phi_{x,1-\eta,1}(z)-1 \bigr)-
\bigl(\phi _{x,1,1}(z)-1 \bigr) \bigr]q_{\sigma}(z)\,\mathrm{d}z
\\
&&{}+\mathbb{I}\{\sigma\geq x\}\int_{x\wedge\sigma}^{(x-\Upsilon
(x))\wedge\sigma} \bigl(
\phi_{x,-\eta,-1}(z)-1 \bigr)q_{\sigma
}(z)\,\mathrm{d}z
\\
&&{}+\mathbb{I}\bigl\{\sigma\geq x-\Upsilon(x)\bigr\}\\
&&\hphantom{{}+}{}\times\int_{x-\Upsilon
(x)}^{\sigma}
\bigl[\bigl(\phi_{x,1-\eta,-1}(z)-1 \bigr)-\bigl(\phi _{x,1,-1}(z)-1 \bigr)
\bigr]q_{\sigma}(z)\,\mathrm{d}z,
\end{eqnarray*}
and we conclude by using that $\pi(\Upsilon(x))=\pi(x)$ and the intermediate
change of variable $z'=\Upsilon(x)-x+z$.
\end{pf*}

\begin{pf*}{Proof of Lemma~\ref{lemsubgeomboundpsithetalessthanx}}
Note first that for $s\in\{-1,1\}$, because $\phi_{x,\eta,s}(z):=[\pi(x+sz)/\pi(x)]^{\eta}=\exp[\eta(\ell(x+sz)-\ell(x))]$,
\begin{eqnarray*}
\phi_{x,\eta,s}'(z) & = & \eta s\ell'(x+sz)
\phi_{x,\eta,s}(z)\quad  \mbox{and}
\\
\phi_{x,\eta,s}''(z) & = & \bigl[
\eta^{2}\bigl|\ell'(x+sz)\bigr|^{2}+\eta
\ell''(x+sz) \bigr]\phi_{x,\eta,s}(z).
\end{eqnarray*}
We now prove the desired upper bounds on $\psi_{x}(z)$ by considering
the following three cases: (a) $0\leq z\leq c_{0}/|\ell'(x)|$, (b)
$c_{0}/|\ell'(x)|\leq z\leq C_{0}/|\ell'(x)|$ and (c)~$C_{0}/|\ell
'(x)|\leq z\leq x$
for an appropriate choice of the constants $c_{0},C_{0}>0$ to be
determined.

Case (a) $0\leq z\leq c_{0}/|\ell'(x)|$. We consider a first-order
Taylor expansion of $\psi_{x}(z)$ at $z_{0}=0$ with integral error
form and obtain
\begin{eqnarray*}
\psi_{x}(z)&=&z\eta\ell'(x)+z(1-\eta)
\ell'(x)-z\ell'(x)\\
&&{}+\int_{0}^{z}
\bigl[\phi_{x,-\eta,-1}''(t)+\phi_{x,1-\eta,1}''(t)-
\phi _{x,1,1}''(t)\bigr](z-t)\,\mathrm{d}t
\\
&=&\int_{0}^{z}a_{x,\eta}(t) (z-t)\,
\mathrm{d}t,
\end{eqnarray*}
where
\begin{eqnarray*}
a_{x,\eta}(t)&:=&\eta^{2}\bigl[\bigl|\ell'(x-t)\bigr|^{2}-
\eta\ell''(x-t) \bigr]\phi _{x,-\eta,-1}(t)
\\
&&{}+ \bigl[(1-\eta)^{2}\bigl|\ell'(x+t)\bigr|^{2}+(1-\eta)
\ell''(x+t) \bigr]\phi _{x,1-\eta,1}(t)\\
&&{}- \bigl[\bigl|
\ell'(x+t)\bigr|^{2}+\ell''(x+t)
\bigr]\phi_{x,1,1}(t).
\end{eqnarray*}
We seek to upperbound $a_{x,\eta}(t)$. We choose $\epsilon_{0}\in
(0,\eta(1-\eta))$
and first show that for any $c_{0}\in(0,\epsilon_{0}/2)$,
%
\begin{equation}\label{eqat-upperbound}
\lim_{x\rightarrow\infty}\inf_{0\leq z\leq c_{0}/|\ell'(x)|}\phi
_{x,1,1}(z)>1-\epsilon_{0}/2.
\end{equation}
Indeed, for $0\leq z\leq c_{0}/|\ell'(x)|$ and $x$ large enough
to ensure $\ell'(x)<0$, we have for some $\xi_{x,z}\in[x,x+z]$, the
following Taylor expansion:
\begin{eqnarray*}
\ell(x+z)-\ell(x) & =&\ell'(x)z+\frac{1}{2}z^{2}
\ell''(x+\xi _{x,z})
\\
& \geq&-c_{0}-\frac{c_{0}^{2}}{2}\frac{|\ell''(x+\xi_{x,z})|}{|\ell
'(x)|^{2}}
\end{eqnarray*}
and with \ref{hyppismallthetasubgeomdrift}\ref
{subhypregularityderivatives}
the last term vanishes as $x\rightarrow\infty$, and we conclude by
the assumption that $-c_{0}>-\epsilon_{0}/2$.

Now choose $\epsilon_{1},\epsilon_{2},\epsilon_{3}>0$. From \eqref
{eqat-upperbound}
and \ref{hyppismallthetasubgeomdrift}\ref
{subhypregularityderivatives}
there exists $R>0$ such\linebreak[4]  that for any $x\geq R$, $\inf_{|z|\leq
c_{0}/|\ell'(x)|}\phi_{x,1,1}(z)\geq1-\epsilon_{0}$,
$\sup_{|t|\leq c_{0}/|\ell'(x)|}|\ell''(x+t)|/\allowbreak |\ell'(x+t)|^{2}\le
\epsilon_{1}$,
$\sup_{|t|\leq c_{0}/|\ell'(x)|}|\ell'(x-t)|^{2}/|\ell
'(x+t)|^{2}\leq1+\epsilon_{2}$
and\linebreak[4]  $\sup_{|t|\leq c_{0}/|\ell'(x)|}|\ell''(x+t)|/|\ell
'(x)|^{2}\leq\epsilon_{3}/c_{0}$.
With these, and observing that for the values considered here we have
$0\le\phi_{x,-\eta,-1}(t),\phi_{x,1-\eta,1}(t),\phi_{x,1,1}(t)\le1$,
we obtain the following upper bound:
\begin{eqnarray*}
a_{x,\eta}(t)&\leq&\bigl|\ell'(x)\bigr|^{2}\frac{|\ell'(x+t)|^{2}}{|\ell
'(x)|^{2}}
\bigl[\eta(\eta+\epsilon_{1})(1+\epsilon_{2})\\
&&\hphantom{|\ell'(x)|^{2}\frac{|\ell'(x+t)|^{2}}{|\ell
'(x)|^{2}}
\bigl[}{}+(1-\eta )
(1-\eta+\epsilon_{1})-(1-\epsilon_{1}) (1-
\epsilon_{0}) \bigr].
\end{eqnarray*}
We consider then the case where $\epsilon_{0}$, $\epsilon_{1}$ and
$\epsilon_{3}$ are chosen small enough so that the term in brackets
in the last display is negative. We note now that since for some $\xi
_{x,t}\in[x,x+t]$,
\[
\ell'(x+t)=\ell'(x)+t\ell''(x+
\xi_{x,t}).
\]
Then with $0\leq t\leq c_{0}/|\ell'(x)|^{2}$ we have
\[
\frac{\ell'(x+t)}{\ell'(x)}\geq1-c_{0}\frac{|\ell''(x+\xi
_{x,t})|}{|\ell'(x)|^{2}},
\]
which leads to the following upper bound:
\begin{eqnarray*}
a_{x,\eta}(t)&\leq&\bigl|\ell'(x)\bigr|^{2}(1-
\epsilon_{3}) \bigl[\eta(\eta +\epsilon_{1}) (1+
\epsilon_{2})\\
&&\hphantom{|\ell'(x)|^{2}(1-
\epsilon_{3}) \bigl[}{}+(1-\eta) (1-\eta+\epsilon _{1})-(1-
\epsilon_{1}) (1-\epsilon_{0}) \bigr].
\end{eqnarray*}
Notice that by our choice of $\epsilon_{0}$ above, we have $\eta
^{2}+(1-\eta)^{2}-(1-\epsilon_{0})\leq-\eta(1-\eta)$.
Now since $\epsilon_{1},\epsilon_{2},\epsilon_{3}>0$ can be chosen
arbitrarily small above, we conclude about the existence of $M>0$,
$c_{0}>0$ and $R>0$ such that for any $x\geq R$
\[
\sup_{|t|\leq c_{0}/|\ell'(x)|}a_{x,\eta}(t)\leq-M\bigl|\ell'(x)\bigr|^{2},
\]
and we therefore conclude that in such a case, for $0\leq z\leq
c_{0}/|\ell'(x)|$
\[
\psi_{x}(z)\leq-M\tfrac{1}{2}z^{2}\bigl|
\ell'(x)\bigr|^{2}.
\]

Case (b) $c_{0}/|\ell'(x)|\leq z\leq C_{0}/|\ell'(x)|$.
First notice that $\psi_{x}(0)=0$ and inspect the derivative of this
function and aim to prove that it is negative. For any $x\in\mathsf{X}$
we have
\begin{eqnarray*}
\psi_{x}'(z) & =&\eta\ell'(x-z)
\phi_{x,-\eta,-1}(z)+(1-\eta)\ell '(x+z)\phi_{x,1-\eta,1}(z)\\
&&{}-
\ell'(x+z)\phi_{x,1,1}(z)
\\
& =&\ell'(x+z) \biggl[\eta\frac{\ell'(x-z)}{\ell'(x+z)}\phi
_{x,-\eta,-1}(z)+(1-\eta)\phi_{x,1-\eta,1}(z)-\phi _{x,1,1}(z)
\biggr].
\end{eqnarray*}
Because $\ell'(x+z)<0$ and the two first terms in brackets form a
convex combination, the second line of (\ref{eqboundonpsi}) will
be established for $c_{0}/|\ell'(x)|\leq z\leq C_{0}/|\ell'(x)|$
once we will have shown that for $x\geq0$ sufficiently large,\vspace*{1.5pt}
\[
\phi_{x,1,1}(z)\leq \biggl(\frac{\ell'(x-z)}{\ell'(x+z)}\phi _{x,-\eta,-1}(z)
\biggr)\wedge\phi_{x,1-\eta,1}(z).
\]
Clearly $1\geq\phi_{x,1-\eta,1}(z)=\phi_{x,1,1}^{1-\eta}(z)\geq
\phi_{x,1,1}(z)$, so we are left with showing that $\phi_{x,1,1}(z)\leq\frac{\ell
'(x-z)}{\ell'(x+z)}\phi_{x,-\eta,-1}(z)$,
or equivalently,\vspace*{1.5pt}
\[
\frac{\pi(x+z)}{\pi(x)} \biggl(\frac{\pi(x-z)}{\pi(x)} \biggr)^{\eta}\leq
\frac{\ell'(x-z)}{\ell'(x+z)}.
\]
We consider the following Taylor expansion:\vspace*{1.5pt}
\begin{eqnarray*}
&&\ell(x+z)-\ell(x)+\eta \bigl[\ell(x-z)-\ell(x) \bigr] \\[1.5pt]
&&\qquad  =z\ell
'(x)+\tfrac{1}{2}z^{2}\ell''(x+
\xi_{x,z})+\eta \bigl[-z\ell '(x)+\tfrac{1}{2}z^{2}
\ell''(x+\xi_{x,-z}) \bigr]
\\[1.5pt]
&&\qquad  =(1-\eta)z\ell'(x)+\tfrac{1}{2}z^{2} \bigl[
\ell''(x+\xi _{x,z})+\eta\ell''(x+
\xi_{x,-z}) \bigr]
\end{eqnarray*}
for some $\xi_{x,z}\in[0,z]$ and $\xi_{x,-z}\in[-z,0]$. For now
choose any $C_{0}>c_{0}$ and notice that for $c_{0}/|\ell'(x)|\leq
z\leq C_{0}/|\ell'(x)|$,
we have that\vspace*{1.8pt}
\begin{eqnarray*}
(1-\eta)z\ell'(x) & \leq&-c_{0}(1-\eta),
\\[1.8pt]
z^{2} \bigl[\bigl|\ell''(x+\xi_{x,z})\bigr|+
\eta\bigl|\ell''(x+\xi_{x,-z})\bigr| \bigr] & \leq&
C_{0}^{2}\frac{|\ell''(x+\xi_{x,z})|+\eta|\ell''(x+\xi
_{x,-z})|}{|\ell'(x)|^{2}}.
\end{eqnarray*}
Let $\epsilon_{1}\in(0,c_{0}(1-\eta))$, and choose $\epsilon_{2}>0$
such that $\exp (-c_{0}(1-\eta)+\epsilon_{1} )<1-\epsilon_{2}$.
By \ref{hyppismallthetasubgeomdrift}\ref{subhypregularityderivatives}
we can conclude by letting $x$ be sufficiently large to ensure that
for $c_{0}/|\ell'(x)|\leq z\leq C_{0}/|\ell'(x)|$,\vspace*{1.8pt}
\[
\frac{\pi(x+z)}{\pi(x)} \biggl(\frac{\pi(x-z)}{\pi(x)} \biggr)^{\eta}\leq\exp
\bigl(-c_{0}(1-\eta)+\epsilon_{1} \bigr)<1-
\epsilon_{2}\leq\frac{\ell'(x-z)}{\ell'(x+z)}.
\]
Now using the result of case (a) we conclude that\vspace*{1.8pt}
\[
\psi_{x}(z)\leq\psi_{x} \biggl(\frac{c_{0}}{|\ell'(x)|}
\biggr)\leq -\frac{M}{2}c_{0}^{2}.
\]

Case (c) $C_{0}/|\ell'(x)|\leq z<x$. We have the following
simple bound:\vspace*{1.8pt}
%
\begin{equation}\label
{equpperboundpsicaseC}
\psi_{x}(z)\leq \biggl(\frac{\pi(x)}{\pi(x-z)} \biggr)^{\eta
}-1+
\biggl(\frac{\pi(x+z)}{\pi(x)} \biggr)^{1-\eta}.
\end{equation}
We inspect, for $C_{0}/|\ell'(x)|\leq z\leq x$ and $x$ large enough,
the following difference:
\begin{eqnarray*}
\ell(x+z)-\ell(x) & =&\int_{0}^{z}
\ell'(x+t)\,\mathrm{d}t
\\
& \leq&\int_{0}^{C_{0}/|\ell'(x)|}\ell'(x+t)\,\mathrm{d}t
\\
& \leq&-C_{0}\sup_{0\leq t\leq C_{0}/|\ell'(x)|}\biggl|\frac{\ell
'(x+t)}{\ell'(x)}\biggr|,
\end{eqnarray*}
and we can similarly obtain a bound on
\[
\ell(x)-\ell(x-z)\leq
-C_{0} \sup_{0\leq t\leq C_{0}/|\ell'(x)|}\biggl|\frac{\ell'(x-t)}{\ell'(x)}\biggr|.
\]

From \ref{hyppismallthetasubgeomdrift}\ref
{subhypregularityderivatives}
and the Taylor expansion $\ell'(x+t)=\ell'(x)+z\ell''(x+\xi_{x,t})$,
we conclude that for $C_{0}$ and $x$ sufficiently large enough, we
can ensure that the upper bound in (\ref{equpperboundpsicaseC})
is negative.

The proof is now concluded by choosing $c_{0}$ as in (a), which leads
to the first line of (\ref{eqboundonpsi}), $C_{0}$ as in (c) and
$R$ large enough to cover cases (b) and (c), which imply the second
line of (\ref{eqboundonpsi}).
\end{pf*}

\begin{pf*}{Proof of Lemma~\ref{lemsubgeomboundonTsthetalargerthanx}}
We start with $T_{1}(\sigma,x)+T_{2}(\sigma,x)$, and with the notation
of Lemma~\ref{lempropertiesUpsilonetc}, we obtain
\begin{eqnarray*}
T_{1}(\sigma,x)+T_{2}(\sigma,x)&\leq&\int
_{0}^{x} \biggl[ \biggl(\frac
{\pi(x-z)}{\pi(x)}
\biggr)^{-\eta}-1 \biggr]q_{\sigma}(z)\,\mathrm {d}z\\
&&{}+\int
_{0}^{\infty} \biggl(\frac{\pi(x+z)}{\pi(x)}
\biggr)^{1-\eta}q_{\sigma}(z)\,\mathrm{d}z
\\
&\leq&\frac{\underline{q}}{\sigma} \bigl[J_{\eta}(x)-x \bigr]+\frac
{\bar{q}}{\sigma}I_{1-\eta}(x).
\end{eqnarray*}
For $\sigma\geq x$, because $\phi_{\Upsilon(x),-\eta,-1}(z)\le1$
in the integration domain,
\[
T_{3}(\sigma,x)  \leq0.
\]
For $\sigma\geq x-\Upsilon(x)\geq x$ we have on the one hand
\begin{eqnarray*}
T_{3}(\sigma,x) & =&\int_{\Upsilon(x)}^{0}
\bigl[\phi_{\Upsilon
(x),-\eta,-1}(z)-1 \bigr]q_{\sigma}\bigl(z+x-\Upsilon(x)
\bigr)\,\mathrm{d}z
\\
& \leq&\frac{\underline{q}}{\sigma} \biggl(\Upsilon(x)+\int_{\Upsilon(x)}^{0}
\phi_{\Upsilon(x),-\eta,-1}(z)\,\mathrm{d}z \biggr)
\\
& \leq&\frac{\underline{q}}{\sigma} \biggl(\Upsilon(x)+\int_{0}^{-\Upsilon(x)}
\phi_{\Upsilon(x),-\eta,1}(z)\,\mathrm{d}z \biggr)
\\
& \leq&\frac{\underline{q}}{\sigma} \bigl(\Upsilon(x)+C\bigl|\Upsilon (x)\bigr|^{1-p}
\bigr),
\end{eqnarray*}
where we have used Lemma~\ref{lempropertiesUpsilonetc}. On the
other hand we also have
\begin{eqnarray*}
T_{4}(\sigma,x) & =&\int_{0}^{\sigma-(x-\Upsilon(x))}
\biggl[ \biggl(\frac{\pi(\Upsilon(x)-z)}{\pi(\Upsilon(x))} \biggr)^{1-\eta
}-\frac{\pi(\Upsilon(x)-z)}{\pi(\Upsilon(x))}
\biggr]\\
&&\hphantom{\int_{0}^{\sigma-(x-\Upsilon(x))}}{}\times q_{\sigma
}\bigl(z+x-\Upsilon(x)\bigr)\,\mathrm{d}z
\\
& \leq&\frac{\bar{q}}{\sigma}\int_{0}^{\infty} \biggl(
\frac{\pi
(\Upsilon(x)-z)}{\pi(\Upsilon(x))} \biggr)^{1-\eta}\,\mathrm{d}z
\\
& \leq&\frac{C}{\sigma}\bigl(-\Upsilon(x)\bigr)^{1-p},
\end{eqnarray*}
where we have again used Lemma~\ref{lempropertiesUpsilonetc}. We
now conclude.
\end{pf*}

\begin{pf*}{Proof of Lemma~\ref{propsub-expon-w-drift}}
For any $x\in\mathsf{X}$ let $A_{\mathsf{Z}}(x):=\{z\in\mathsf
{Z}\dvtx \pi(x+z)/\pi(x)\geq1\}$
and $R_{\mathsf{Z}}(x):=A_{\mathsf{Z}}^{c}(x)$ (where the complement
is with respect to $\mathsf{Z}$) and $A(x):=x+A_{\mathsf{Z}}(x)$.
Without loss of generality we focus on the case $x>0$. From Lemma~\ref{lempropertiesUpsilonetc} there exists $R_{1}>0$ such that
for any $x\geq R_{1}$, $R_{\mathsf{Z}}(x)=(-\infty,-x+\Upsilon
(x))\cup(0,\infty)$
and $A_{\mathsf{Z}}(x)=[-x+\Upsilon(x),0]$, where $\Upsilon(x)$
is as in Lemma~\ref{lempropertiesUpsilonetc}. For $x\geq R_{1}$
and $\sigma\leq1$, we have the inequalities
\begin{eqnarray*}
\alpha_{\sigma}(x) & =&\int_{\mathsf{Z}}\min \biggl\{ 1,
\frac{\pi
(x+z)}{\pi(x)} \biggr\} q_{\sigma}(z)\,\mathrm{d}z
\\
& =&1+\int_{R_{\mathsf{Z}}(x)} \biggl[\frac{\pi(x+z)}{\pi
(x)}-1
\biggr]q_{\sigma}(z)\,\mathrm{d}z
\\
& \geq&1-\int_{R_{\mathsf{Z}}(x)}q_{\sigma}(z)\,\mathrm{d}z
\\
& =&\frac{1}{2}-\int_{-\infty}^{(-x+\Upsilon(x))/\sigma}q(z)\,\mathrm
{d}z
\\
& \geq&\frac{1}{2}-\int_{-\infty}^{-x/\sigma}q(z)
\,\mathrm{d}z.
\end{eqnarray*}
Now with $\mu_{1}<\infty$ the first-order moment of $q$, we notice
that from Chebyshev's inequality and for $x\geq R_{1}$,
\[
\int_{x/\sigma}^{\infty}q(z)\,\mathrm{d}z\leq\sigma
R_{1}^{-1}\times \mu_{1}
\]
from which we deduce the first statement for $\sigma\leq1$ and $x\geq R_{1}$.
Now for $x\geq R_{1}$ and $\sigma\geq1$,
\begin{eqnarray*}
\alpha_{\sigma}(x)&\le&\frac{\bar{q}}{\sigma} \biggl(\int_{A_{\mathsf{Z}}(x)\cup R_{\mathsf{Z}}(x)}
\min \biggl\{ 1,\frac{\pi
(x+z)}{\pi(x)} \biggr\} \,\mathrm{d}z \biggr)
\\
&\le&\frac{\bar{q}}{\sigma} \biggl(2C_{\Upsilon,1} \bigl(-\log\bigl(\pi
(x)/C_{\Upsilon,2}\bigr) \bigr)^{1/p}\\
&&\hphantom{\frac{\bar{q}}{\sigma} \biggl(}{}+\int_{0}^{\infty}
\frac{\pi
(\Upsilon(x)-z)}{\pi(\Upsilon(x))}\,\mathrm{d}z+\int_{0}^{\infty
}
\frac{\pi(x+z)}{\pi(x)}\,\mathrm{d}z \biggr)
\\
&\leq& C\frac{ (-\log(\pi(x)/C_{\Upsilon,2})
)^{1/p}}{\sigma},
\end{eqnarray*}
where we have used the results of Lemma~\ref{lempropertiesUpsilonetc}
to upper bound the Lebesgue measure of $A_{\mathsf{Z}}(x)$ and the
last two integrals. We now turn to the case $0\leq x\leq R_{1}$.
Let $M>0$ such that $\int_{M}^{\infty}q(z)\,\mathrm{d}z\leq1/4$ and
$\sigma\leq1$
and with $\phi_{x}(z)=\pi(x+z)/\pi(x)$
\begin{eqnarray*}
\alpha_{\sigma}(x)&=&1+\frac{1}{\sigma}\int_{\mathsf{Z}}
\biggl(1\wedge\frac{\pi(x+z)}{\pi(x)}-1 \biggr)q \biggl(\frac{z}{\sigma
} \biggr)\,
\mathrm{d}z
\\
&\geq&1+\int_{-\infty}^{-M} \biggl(1\wedge
\frac{\pi(x+\sigma z)}{\pi
(x)}-1 \biggr)q(z)\,\mathrm{d}z\\
&&{}+\int_{M}^{\infty}
\biggl(1\wedge\frac
{\pi(x+\sigma z)}{\pi(x)}-1 \biggr)q(z)\,\mathrm{d}z
-\int
_{-M}^{M}\biggl\llvert \frac{\pi(x+\sigma z)}{\pi(x)}-1
\biggr\rrvert q(z)\,\mathrm {d}z
\\
&\ge&1-2\int_{M}^{\infty}q(z)\,\mathrm{d}z-2M\bar{q}\sup
_{x\in
B(0,R_{1}),z\in B(0,M)}\bigl\llvert \phi_{x}'(z)\bigr
\rrvert \sigma,
\end{eqnarray*}
and we deduce the first statement of the lemma. We now consider the
case $0\leq x\leq R_{1}$ and $\sigma\geq1$. There exists (cf. the
proof of Lemma~\ref{lempropertiesUpsilonetc}) $R_{2}>0$ such that
for all $x\leq R_{1}$
\begin{eqnarray*}
\alpha_{\sigma}(x)&\le&\sigma^{-1}\int_{-R_{2}}^{R_{2}}q(z/
\sigma )\,\mathrm{d}z+\sigma^{-1}\int_{-\infty}^{-R_{2}}
\frac{\pi
(x+z)}{\pi(x)}q(z/\sigma)\,\mathrm{d}z\\
&&{}+\sigma^{-1}\int
_{R_{2}}^{\infty}\frac{\pi(x+z)}{\pi(x)}q(z/\sigma)\,\mathrm{d}z.
\end{eqnarray*}
From the proof of Lemma~\ref{lempropertiesUpsilonetc}, we have the
bound $\pi(x+z)/\pi(x)\leq C_{1}\*\exp (-C_{2}|z|^{p} )$ for
some $C_{1},C_{2}>0$ and since $q(z)\leq\bar{q}$, we deduce the existence
of $C>0$ such that for $x\leq R_{1}$ and $\sigma\geq1$ we have
$\alpha_{\sigma}(x)\leq C/\sigma$.\vadjust{\goodbreak}
\end{pf*}
\end{appendix}




\printaddresses

\end{document}